\documentclass[a4paper,11pt, reqno]{amsart}
\usepackage{amssymb,amsmath,latexsym, upgreek,comment}
\usepackage{epsfig,caption}
\usepackage{mathrsfs}
\usepackage{hyperref}
\usepackage{color,graphics,subfigure}
\usepackage{mathtools,textcomp}
 \usepackage{systeme}

\usepackage{epic,eepic,wrapfig,ifthen,enumerate}
\usepackage{url}
\newtheorem{thm}{Theorem}[section]

\newtheorem{prop}[thm]{Proposition}
\newtheorem{defn}[thm]{Definition}

\newtheorem{remark}[thm]{Remark}

\newtheorem{cor}[thm]{Corollary}
\newtheorem{lem}[thm]{Lemma}
\numberwithin{equation}{section}

\allowdisplaybreaks

\def\pf{{\medskip\noindent {\bf Proof. }}}
\def\qed{{\hfill $\Box$ \bigskip}}

\DeclareMathOperator*{\esssup}{ess\,sup}

\DeclareMathOperator{\dist}{dist}

\DeclareFontFamily{U}{mathx}{}
\DeclareFontShape{U}{mathx}{m}{n}{<-> mathx10}{}
\DeclareSymbolFont{mathx}{U}{mathx}{m}{n}
\DeclareMathAccent{\wc}{0}{mathx}{"71}

\newcommand{\cal}[1]{\mathcal{#1}}
  \def\sC {{\cal C}}
 \def\sE {{\cal E}} \def\sF {{\cal F}}

 \def\sN {{\cal N}} 
  
 \def\sT {{\cal T}}

\def\R {{\mathbb R}}
\def\N {{\mathbb N}}
\def\Z {{\mathbb Z}}
\def\P{{\mathbb P}}
\def\E{{\mathbb E}}

\def\EE{{\mathcal E}}
\def\eps{\varepsilon}

\def\FF{{\mathcal F}}
\def\EE{{\mathcal E}}

\def\R{{\mathbb R}}

\def\E{{\mathbb E}}

\def\P{{\mathbb P}}
\def\N{{\mathbb N}}
\def\eps{\varepsilon}
\def\wh{\widehat}
\def\wt{\widetilde}
\def\pf{\noindent{\bf Proof.} }

\def\ub{{\overline{\beta}}}

\def\diam{{\text{\rm diam}}}
\def\1{{\bf 1}}

\def\nn{\nonumber}

\def\qed{{\hfill $\Box$ \bigskip}}

\def\eps{\varepsilon}
\def\wh{\widehat}
\def\wt{\widetilde}

\def\vp{{\varphi}}

\def\As{{\bf (A)}}

\makeatletter
\@addtoreset{equation}{section}

\makeatother

\begin{document}
\bibliographystyle{plain}

\title[HKE for Markov processes with  blowing-up jump kernels]
{ \bf   Heat kernel estimates for Markov processes with  blowing-up jump kernels}

\author{Soobin Cho, \quad Panki Kim, \quad Renming Song \quad and \quad Zoran Vondra\v{c}ek}

\address[Cho]{School of Mathematics and Statistics, Beijing Institute of Technology, Beijing 100081, China}
\curraddr{}
\email{soobincho.math@gmail.com}

\address[Kim]{Department of Mathematical Sciences and Research Institute of Mathematics,
	Seoul National University,	Seoul 08826, Republic of Korea}
\thanks{}
\curraddr{}
\email{pkim@snu.ac.kr}
\thanks{Panki Kim is supported by the National Research Foundation of Korea(NRF) grant funded by the Korea government(MSIP) (No. RS-2026-25518141).} 

\address[Song]{
	Department of Mathematics, University of Illinois Urbana-Champaign, Urbana, IL 61801,
	USA}
\curraddr{}
\thanks{}
\email{rsong@illinois.edu}
\thanks{Research of Renming Song is supported in part by a grant from	the Simons Foundation \#960480.}

\address[Vondra\v{c}ek]
{
	Dr.~Franjo Tu{\dj}man Defense and Security 
	University
	 and Department of Mathematics, Faculty of Science, University of Zagreb, Zagreb, Croatia
}
\curraddr{}
\thanks{} 
\email{vondra@math.hr}
\thanks{Zoran Vondra\v{c}ek is supported in part by the Croatian Science Foundation under the project IP-2025-02-8793.}

\date{}

\begin{abstract} 
In this paper, we establish sharp two-sided heat kernel estimates for a large class of purely discontinuous symmetric Markov processes on 
closed subsets $F$ of $\mathbb{R}^d$, 
 whose jump kernels blow up on a 
Borel subset $\Sigma$ of $F$. 
 We assume that $F\setminus \Sigma$ is  
 a $\kappa$-fat set and is dense in $F$,
  and additionally, in case of unbounded $F$, that the Assouad dimension of $F$ is strictly less than $d$.  
To the best of our knowledge, this is the first work establishing sharp  two-sided heat kernel estimates for jump processes 
whose jump kernels  blow up  on part of the state space.

The jump kernels under consideration take the form 
$J(x,y)=|x-y|^{-d-\alpha}{\mathcal B}(x,y)$, where $\alpha\in (0,2)$ and 
 the function 
${\mathcal B}(x,y)$ blows up at a subset $\Sigma$ of  
$F$. 
A fundamental obstacle is that the tails of the jump measures 
 are not uniformly bounded, and hence 
 standard techniques in heat kernel analysis do not provide a priori off-diagonal estimates. 
To overcome this difficulty, we develop a new approach based on  weighted integral estimates for the heat kernel that are sensitive to both the blow-up behavior of the jump kernel and the geometry of $F\setminus \Sigma$. 

  Examples of processes falling 
 within our general
  framework include traces of isotropic $\alpha$-stable processes in  some unbounded $C^{1,\rm Dini}$ sets,  
  processes in  
  Lipschitz sets arising in connection with the nonlocal Neumann problem, and a large class of resurrected self-similar processes in the closed upper half-space.

  \bigskip
  \medskip
  
  
\end{abstract}

\maketitle

\section{Introduction}\label{s:intro}	
	Analysis of discontinuous Markov processes and their associated nonlocal Dirichlet forms has been one of the most rapidly expanding areas in probability theory and partial differential equations over the last
few decades, see, for example, \cite{BGK09, CKS87, CF12, CKW-jems, CKW-memo, GHH18, GHH24} and the references therein. 
	
	 Following the seminal 
	 work
	 on stable-like processes \cite{CK03}, a rich and robust 
	 framework has been developed 
	 for establishing sharp heat kernel 
	 estimates and regularity theory for discontinuous Markov processes associated with nonlocal Dirichlet forms; 
see for instance, \cite{CKW-jems, CKW-memo, GHH18, GHH24} and the references therein.
	 In a broad range of  settings—ranging from metric measure spaces to 
	subsets of $\R^d$ with smooth or rough boundaries—the 
	  standard assumption has almost invariably been that 
	 the tails of the associated jump measures are uniformly bounded.

	In this paper, we depart from 	the classical framework above
	by allowing jump kernels to blow up on
	 certain parts of the state space. More precisely, we
	study a class of purely discontinuous symmetric Markov processes on closed subsets $F$ of $\mathbb{R}^d$, $d \ge 1$, whose jump kernels take the form:$$J(x,y)=|x-y|^{-d-\alpha}{\mathcal B}(x,y), \quad \alpha\in (0,2)$$ where  ${\mathcal B}(x,y)$ 
	blows up 
		on a subset, with Assouad dimension strictly less than $d$, of $F$.
	
	In our setting, the tails of the jump measures are not uniformly bounded and 
 the parabolic Harnack principle  fails. Thus, the standard techniques 
from the existing literature for 
deriving
heat kernel estimates are not applicable.  Indeed, in our setting, the correct off-diagonal behavior of the heat kernel $p(t,x,y)$ depends on $x$ and $y$ in a way that is not captured solely by $|x-y|$. In particular, for some $x \neq y$, the heat kernel $p(t,x,y)$ may diverge as $t \to 0$; see Remark \ref{r:th}(iv). To overcome 
these difficulties, 
we develop a new method 
for establishing
weighted integral estimates for heat kernels, based on recently developed stability results, which establish the equivalence of weighted functional inequalities, Hölder regularity, and near-diagonal heat kernel lower bounds in \cite{CK25+}. We then derive point-wise heat kernel bounds from these integral estimates via bootstrap arguments based on Meyer's construction.
To the best of our knowledge, this work provides the first sharp two-sided heat kernel estimates for discontinuous Markov processes 
whose jump kernels blow up on part of the state space.

	The present paper substantially extends the  framework developed in \cite{KSV-jfa25} for non-conservative self-similar Markov processes in the half-space.
Our framework accommodates a diverse family of processes. Examples include traces of isotropic $\alpha$-stable processes on 
 some unbounded $C^{1,\rm Dini}$ sets; 
processes in Lipschitz sets 
 arising 
in connection with the nonlocal Neumann problem; and a large class of resurrected self-similar processes in the closed half-space.
	 
The jump kernel blowing up at certain subset of the state space is one way for the Dirichlet form to be degenerate. 
Another type of degeneracy occurs when the jump kernel decays to 0 at certain subset of the state space. The potential theory of Markov processes with jump kernels  vanishing 
at the boundary has been studied systematically in \cite{CKSV-aop, CKSV24, CKSV25, KSV-pota23, KSV-jems24, KSV-ma24}. The papers \cite{ CKSV24, KSV-pota23, KSV-jems24, KSV-ma24} focus on elliptic aspects, such as the boundary Harnack principle and Green function estimates, while \cite{CKSV-aop, CKSV25} are devoted to heat kernel estimates.

We now describe the setting of the current paper and state the main result.
For $E\subset \R^d$, we use $\diam(E)$ to denote the diameter of $E$ 
and ${\rm dim}_{\rm A}(E)$ its Assouad dimension (see Definition \ref{def:Assouad} for the precise  definition).

\begin{defn}\label{d:fat}
\rm Let $\kappa \in (0,1)$. A 
 set $E\subset \R^d$  is said to be \textit{$\kappa$-fat} if there exists a localization constant $R_0 \in (0,\diam(E)]$ such that for any $x \in  E$ and $r\in (0,R_0)$, there exists a point $z \in E$ such that $B(z, \kappa r) \subset  B(x,r)\cap E$. 
\end{defn}

{\it  Throughout  this paper,  
	 $F\subset \R^d$ is a closed set,  
	 $\Sigma\subset F$ is a Borel set and $D:=F\setminus \Sigma$.
	 We assume that $D$  is  
	 a $\kappa$-fat set and is dense in $F$.
	 In case when $F$ is unbounded, we additionally assume that ${\rm dim}_{\rm A}(\Sigma)<d$. 
	 }

\medskip

We note that if $F$ is a closed, bounded set and $\Sigma \subset F$ is a Borel set such that $D:=F\setminus \Sigma$ is a $\kappa$-fat set  dense in $F$, then ${\rm dim}_{\rm A}(\Sigma)<d$. 
However, when $F$ is unbounded, it can happen that ${\rm dim}_{\rm A}(\Sigma)=d$ even if the interior of $F$ is a connected,  smooth open set; see Remark \ref{r:some-examples}. 

Here are some examples that fit into our framework. Let $F$ be the closure of a bounded 
Lipschitz open set in $\R^d$, $\Sigma_{d-1}:=\partial F$, $\Sigma_0$ be a finite set in the interior of $F$, and, for $k=1, \dots, d-2$, let $\Sigma_k$ be a smooth hyper-surface of dimension $k$ in the interior of $F$. If $\Sigma$ is a Borel subset of $\cup^{d-1}_{k=0}\Sigma_k$ and $D:=F\setminus \Sigma$, then 
the assumption above is satisfied.

The constant
\begin{align}\label{e:def-gamma}	\gamma	:= d- 	\text{dim}_{\rm A}(\Sigma) \in (0,d],\end{align}which is the	\textit{Assouad codimension} of $\Sigma$ (see Definition \ref{def:Assouadc} for the precise  definition), will play an important role in this paper.  
When $F$ is the  closure of a 
Lipschitz open set and $\Sigma=\partial F$, we have $\gamma=1$. 

\medskip

We say that $\Phi:[0,\infty) \to [1,\infty)$ is a \textit{blow-up weight function} if it is continuous, increasing, 
satisfies  $\Phi(0)=1$ and 
 the following upper scaling property: there exist constants $\beta\ge 0$ and $C>1$ such that
\begin{align}\label{e:upper-scaling}
	\Phi(r)/\Phi(s) \le C(r/s)^{ \beta} \quad \text{for all $0<s\le r$}.
\end{align}
The upper  Matuszewska index $\overline \beta$ of $\Phi$ is defined as  the infimum of all $ \beta$ 
for which \eqref{e:upper-scaling} holds. 
 We define $\Phi(\infty):=\lim_{a \to \infty}\Phi(a)$.
 Although our results also cover the case $\Phi(\infty)<\infty$, which corresponds to the standard setting, the main interest and novelty of this work lie in the regime $\Phi(\infty)=\infty$. 

Typical examples of blow-up weight functions are  
$\Phi(r)=\log(e+r)$ and $\Phi(r)=1\vee r^{\beta}$ with $\beta\ge 0$.

Consider the bilinear form
\begin{align*}
	\sE(u,v) =\frac12 \int_{F\times F} (u(x)-u(y))(v(x)-v(y)) J(x,y) dxdy,
\end{align*}
where $J(x,y)$ is   a symmetric 
Borel function on 
$F\times F$.

We assume that $J$ satisfies the following assumption:

\medskip

\setlength{\leftskip}{0.17in} 
\setlength{\rightskip}{0.17in}

\noindent 
\As \
 There exist a blow-up weight function $\Phi$  with upper Matuszewska index $\overline \beta$ satisfying
\begin{align}\label{e:beta-range}
\overline \beta<\gamma \wedge \alpha
\end{align} 
and constants $A_0 \in (0,\infty]$ and $C_1>1$ such that for all 
 $x,y \in F$, 
\begin{align}\label{e:JJ}
	&\frac{C_1^{-1}}{|x-y|^{d+\alpha}}\Phi\left( \frac{(|x-y|\wedge A_0)^2}{(\delta_\Sigma(x)\wedge A_0) (\delta_\Sigma(y) \wedge A_0)}\right)\nn\\
	&\le 	J(x,y) \le \frac{C_1}{|x-y|^{d+\alpha}}\Phi\left( \frac{(|x-y|\wedge A_0)^2}{(\delta_\Sigma(x)\wedge A_0) (\delta_\Sigma(y) \wedge A_0)}\right).
\end{align}

\setlength{\leftskip}{0in} 
\setlength{\rightskip}{0in}

Here and throughout this paper, we use the notation $a\wedge b:=\min(a, b)$ and $a\vee b:=\max(a, b)$,
and  write 
 $\delta_A(x):=\mathrm{dist}(x,A)$ for the distance from $x$ to $A$.  
Since $\Phi \ge 1$, \As \ implies that
\begin{align}\label{e:J-standard-lower}
	J(x,y) \ge C_1^{-1} |x-y|^{-d-\alpha} \quad \text{for all $x,y \in 
	F$}.
\end{align}

We remark 
that,  in many examples, one can take the  constant  $A_0$ to be $\diam(F^c)$.
When $A_0=\infty$, \eqref{e:JJ} reduces to 
\begin{align}\label{e:JJ-s}
	\frac{C_1^{-1}}{|x-y|^{d+\alpha}}\Phi\left( \frac{|x-y|^2}{\delta_\Sigma(x) \delta_\Sigma(y)}\right)
	\le 	J(x,y) \le \frac{C_1}{|x-y|^{d+\alpha}}\Phi\left( \frac{|x-y|^2}{\delta_\Sigma(x) \delta_\Sigma(y) }\right).
\end{align}

The motivation for \eqref{e:JJ} and \eqref{e:JJ-s} comes from \cite{KSV-jfa25}, where, in the case of  
half-spaces, 
a general method for constructing jump kernels $J$ satisfying \eqref{e:JJ-s} was developed. In particular, 
the trace process  and the process related to the nonlocal Neumann problem in half-spaces
 satisfy estimates of 
the form \eqref{e:JJ-s}
with 
$\Phi(r)=1\vee r^{\alpha/2}$
and $\Phi(r)=\log(e+r)$, respectively. In Sections \ref{s:appl} and \ref{s:proof-examples} we (i) show that analogous estimates hold for these processes in 
more general (smooth) sets; and (ii) revisit the general kernels in the upper half-space together with 
their detailed construction.

\medskip

Let ${\rm Lip}_c(F)$ be the family of all Lipschitz functions on $F$ with compact supports. 
For  $u \in L^2( F)$, define $\sE_1(u,u):= \sE(u,u) + \lVert u \rVert_{L^2(F)}^2$.   
We will see that, under \As, $\sE(u,u)<\infty$ for all $u\in {\rm Lip}_c(F)$ (see Proposition  \ref{p:regular-Dirichlet-form} below).
Let $\sF$ be the $\sE_1$-closure of ${\rm Lip}_c(F)$ in $L^2(F)$. Then $(\sE,\sF)$ is a regular Dirichlet form on $L^2(F)$ under \As.

We now state the main result of this paper. 

\medskip
\begin{thm}\label{t:main-HKE}
 Let $d\ge 1$, $\alpha \in (0,2)$, 
 $F\subset \R^d$ be closed,   
 $\Sigma\subset F$ be a  Borel set and $D:=F\setminus \Sigma$. 
 Assume that $D$ is dense in $F$ and is $\kappa$-fat  with localization constant $R_0$. 
 In case $F$ is unbounded, assume additionally that ${\rm dim}_{\rm A}(\Sigma)<d$. If \As \ holds, then 
 $(\sE,\sF)$ 
admits  a jointly continuous heat kernel $p(t,x,y)$ on $(0,\infty) \times F\times F$ satisfying the following estimates: 
For every $T \in (0, \infty)$, there exists a constant $C \ge 1$ such that, 	 for all $(t,x,y) \in (0, T \vee R_0^\alpha) \times F\times F$,
	\begin{align}\label{e:main-HKE1}
	&C^{-1} \left[t^{-d/\alpha} \wedge  \left( \frac{ t }{|x-y|^{d+\alpha}}\Phi\left( \frac{(|x-y| \wedge A_0)^2}{((\delta_\Sigma(x) \vee t^{1/\alpha}) \wedge A_0)((\delta_\Sigma(y) \vee t^{1/\alpha}) \wedge A_0) }\right) \right)\right] \nn\\
	&\le p(t,x,y)  \le
	C \left[t^{-d/\alpha} \wedge  \left( \frac{ t }{|x-y|^{d+\alpha}}\Phi\left( \frac{(|x-y| \wedge A_0)^2}{((\delta_\Sigma(x) \vee t^{1/\alpha}) \wedge A_0)((\delta_\Sigma(y) \vee t^{1/\alpha}) \wedge A_0) }\right) \right) 
	\right].
	\end{align}
	\end{thm}

	\begin{remark}\label{r:th} 
	{\rm
(i) When $\Phi\equiv 1$, \eqref{e:main-HKE1} reduces to the classical  heat kernel estimates of stable-like processes obtained in \cite{CK03}.

\noindent
(ii)
For $x\in {F}$ and $t \in (0, T \vee R_0^\alpha)$, let $x(t) \in D$ be a point satisfying
\begin{align}\label{e:def-x(t)}
	|x(t)-x| < c_1t^{1/\alpha} \quad \text{and} \quad 	\delta_\Sigma(x(t))  \ge \delta_\Sigma(x) \vee (\kappa c_1 t^{1/\alpha}),
\end{align}
where 
 $c_1:=1 \wedge (R_0/T^{1/\alpha})$. 
For example, let $x(t)=x$  if $\delta_\Sigma(x) \ge \kappa  c_1t^{1/\alpha}$, and otherwise choose  $x(t) \in B(x,c_1t^{1/\alpha})$ such that $B(x(t), \kappa c_1t^{1/\alpha}) \subset D$.
Using \As \
and \eqref{e:def-x(t)}, 
one can get the following equivalent form of \eqref{e:main-HKE1}:
		\begin{align*}
	C^{-1} \left[t^{-d/\alpha} \wedge \left( t J(x(t), y(t))  \right)\right]	 \le p(t,x,y) \le C \left[t^{-d/\alpha} \wedge \left( t J(x(t), y(t))  \right)\right]. 
	\end{align*} 

	\noindent 
	(iii)
	When 	$D$ is  $\kappa$-fat  with localization constant $R_0=\infty$, then $F$ is unbounded and our heat kernel estimates are 
global in time, 
 i.e., \eqref{e:main-HKE1} holds for all $(t,x,y) \in (0, \infty) \times F\times F$.

	\noindent 
	(iv)
	When $\Phi(r)=1\vee r^{\beta}$ with 
	$\alpha/2 <\beta<\gamma \wedge \alpha$, 
	 \eqref{e:main-HKE1} implies that
	for any $x, y \in \Sigma$ and  $t^{1/\alpha}< A_0$,
	$$p(t,x,y) \ge C^{-1}  \left[t^{-d/\alpha} \wedge  \bigg( \frac{ t^{1- (2\beta/\alpha)}(|x-y| \wedge A_0)^{2\beta}}{|x-y|^{d+\alpha}
	} \bigg) \right].
	 $$
	 Thus, for any $x, y \in \Sigma$,  
	$\lim_{t \to 0} p(t,x,y) =\infty $. 
 }
\end{remark}

We now describe the organization of the paper. In Section \ref{s:prelim} we prove some preliminary results. These are of two types.  Subsection \ref{ss:prelim-1} contains several simple geometric and measure-theoretic consequences of the assumption that 
$D$ is  a $\kappa$-fat set and is dense in $F$. The main result of this subsection is Proposition \ref{p:Aikawa}.
In Subsection \ref{ss:prelim-2}, we prove several estimates involving the function $\Phi$, the jump kernel $J(x,y)$, and various integrals thereof. The proofs are elementary although sometimes tedious.

Section \ref{s:prel-anal-df} is devoted to  
some preliminary analysis of 
the Dirichlet form $(\sE,\sF)$ and its truncated versions, and 
relies on the results obtained in \cite{CK25+}. Subsection \ref{ss:Hardy} begins with a Hardy-type inequality in Proposition \ref{p:Hardy-fractional-Laplacian}, whose proof uses the fact that ${\rm dim}_{\rm A}(\Sigma)<d$. 
The truncated form $\sE^{(\rho)}$ 
with jump kernel $J^{(\rho)}(x,y):=\1_{\{|x-y|<\rho\}}J(x,y)$ 
is introduced next and compared in Proposition \ref{p:truncation-equivalent} with the form $\sE$. This leads to the conclusion that 
$(\sE^{(\rho)}, \sF)$ 
is also a regular Dirichlet form on $L^2(F)$. As 
in \cite{CKW-memo}, truncated Dirichlet forms play a crucial role in 
establishing  the heat kernel estimates in this paper.

In Subsection \ref{ss:reg}, we show that $\Phi(r/\delta_\Sigma(x))$ is an admissible weight function, 
in the sense of  \cite{CK25+}, 
for both $\sE$ and $\sE^{(\rho)}$. 
Then we use \cite[Theorem 12.1]{CK25+}  to establish  parabolic H\"older regularity of caloric functions, a Nash-type inequality, existence of heat kernels for the form $\sE$ and its truncated version  $\sE^{(\rho)}$, as well as near-diagonal lower estimates. At the end of this section, 
we also show that the parabolic Harnack inequality does not hold in 
our setting  when $\Phi(\infty)=\infty$.

In Section \ref{s:further-anal} we continue our analysis of the truncated form. 
We first establish several inequalities between 
 heat kernels $p^{(\rho)}(t,x,y)$ and $p^{(\rho')}(t,x,y)$ of the truncated forms
for different parameters $\rho, \rho'\in (0, \overline{R}]$. 
These results rely on Meyer’s construction; however, standard results in the literature cannot be applied in our setting, since the tail of the jump measure is not uniformly bounded. This necessitates substantial modifications via the Mosco convergence method.
This analysis is carried out in Subsection \ref{ss:Meyer}. In Subsection \ref{ss:4.2}, we prove several off-diagonal estimates for the truncated Dirichlet form. Corollary \ref{c:truncated-off-diagonal} then serves as the basis of a bootstrap argument used to obtain the sharp heat kernel upper bound.

Section \ref{s:hke} forms the core of the paper. In the short Subsection \ref{ss:t-lhke}, we 
first establish the sharp lower bound for the heat kernel when $t\in (0, R_0^{\alpha})$. The proof is rather simple. The near-diagonal lower bound follows from Lemma \ref{l:NDL-modified}. For the off-diagonal lower bound, we use 
the strong Markov property, L\'evy system formula, and judiciously chosen small sets
on which the already obtained near-diagonal estimate can be applied. Proof of the upper bound when $t\in (0, R_0^{\alpha})$ is given in Subsection \ref{ss:t-uhke}. The near-diagonal bound follows from Proposition \ref{p:NDU}. The off-diagonal bound is considerably more difficult and constitutes the main challenge of the paper. 
We first use Corollary \ref{c:truncated-off-diagonal} to prove a bootstrap lemma (Lemma \ref{l:UHK-without-Phi}) and then use it to prove the preliminary upper bounds on $p(t, x, y)$ and $p^{(\rho)}(t, x, y)$ in Proposition \ref{p:UHK-tt}. Next we prove another bootstrap lemma (Lemma \ref{l:UHK-onesided-Phi}) and use it to prove the preliminary upper bounds in Proposition \ref{p:UHK-onesided-Phi-iterated}. 
Finally, we use Proposition \ref{p:UHK-onesided-Phi-iterated} and some detailed analysis  to prove the sharp upper bound in Theorem \ref{t:main-HKE} for $t\in (0, R_0^{\alpha})$. 
In Subsection \ref{ss:proof-main}, we use the semigroup property to  extend both the upper and lower bounds to all times $t<T\vee R^\alpha_0$. 

In Section \ref{s:appl}, we present three examples that fall within our framework: the trace process, the process associated with the nonlocal Neumann problem, 
and the general resurrected processes in the upper half-space introduced in \cite{KSV-jfa25}. 
Proofs of the technical estimates that the corresponding jump kernels satisfy condition \As \ are deferred to Section \ref{s:proof-examples}.

We end this introduction with a few words on notation. 
We use $m_d$ to denote the $d$-dimensional Lebesgue measure. 
 For a Borel set $E\subset \R^d$, $x\in \R^d$ and $r>0$, we use the shorthand notation $B_{E}(x, r):=B(x,r)\cap E$.
We also use the shorthand notation $\overline R : =\diam(F) \in (0,\infty]$.  
Throughout this paper, $A\asymp B$ means that 
 the quotient of $A$ and $B$
 is bounded between two positive constants on the specified region.
The  constants $d$, $\gamma$, $\alpha$, $\kappa$ and $R_0$ are fixed throughout the paper. The notation $C=C(a,b,\ldots)$ indicates that the positive constant $C$ depends on $a, b, \ldots$. The dependence on $d$, $\gamma$, $\alpha$, $\kappa$ and $R_0$  may not be mentioned explicitly. Lower case letters  $c_i$, $i=1,2,  \dots$, are used to denote
positive constants in the proofs and the labeling of these constants starts anew in each proof.

\section{Preliminary results}\label{s:prelim}

{\it 
	Throughout this paper, we assume that  
	$F\subset \R^d$ is a closed set,  
	$\Sigma$ is a Borel subset of $F$, $D:=F\setminus \Sigma$ is $\kappa$-fat  with localization constant $R_0 \in (0,\overline R]$ 
and is dense in $F$, and 	that \As \ holds.  
}

 \subsection{Some elementary geometric and measure-theoretic results}\label{ss:prelim-1}

In this subsection, we recall the definition of the Assouad dimension and present some elementary 
geometric and measure-theoretic consequences of the  assumptions that  $D$ is a dense Borel subset of $F$ and $D$ is  $\kappa$-fat. 

\begin{defn}\label{def:Assouad}
		\rm Let  $E \subset \R^d$ be a non-empty subset. Let $\mathfrak I(E)$ be the set of all $\lambda>0$ satisfying the following: there exists a constant $C=C(\lambda)>0$ such that for all $x\in E$ and  $0<s\le r<\diam(E)$, the set $B(x,r) \cap E$ can be covered by at most  $N\le C(r/s)^\lambda$ open balls of radius $s$ centered in $E$. The \textit{Assouad dimension} of $E$ is defined by
		\begin{align*}
			{\rm dim}_{\rm A}(E):= \inf\{\lambda>0: \lambda \in \mathfrak I(E)\}.
		\end{align*} 
		The Assouad dimension of a singleton $\{x\}$ is defined to be 0.
	\end{defn}

\begin{remark}\label{r:some-examples} {\rm (i)
When $F$ is bounded, since $D$ is dense in $F$ and $\kappa$-fat, $\Sigma=F\setminus D$ is porous in $F$
in the sense of \cite[Section 5]{Lu98}. Therefore, by \cite[Theorem 5.2]{Lu98}  and the implication (1) $\Rightarrow$ (4') in \cite[Theorem A.3]{Lu98}, we have  ${\rm dim}_{\rm A}(\Sigma)<d$. 
(In \cite{Lu98}, the author uses a slightly different definition of the Assouad dimension based on $(C,s)$-homogeneity.  By \cite[Theorem A.3]{Lu98}, it is equivalent to Definition \ref{def:Assouad}.) 
Note that this argument does not apply to unbounded $F$, since porosity in  \cite{Lu98}  is defined with the localization constant equal to the diameter of the set, whereas the localization constant $R_0$ in Definition \ref{d:fat} may be finite for unbounded $F$. 

\noindent (ii) 
If the interior of $F$ is  a bounded Lipschitz set, or  a half-space-like Lipschitz open set (see Definition \ref{d:half-space-like}
for the precise definition) or an exterior Lipschitz open set,  it is easy to see that  ${\rm dim}_{\rm A}(\partial F)=d-1$. 

\noindent  (iii) There exist unbounded, connected, closed sets $F$ with smooth boundary such that ${\rm dim}_{\rm A}(\partial F)=d$. For instance, consider $F:=\R^d \setminus \cup_{k \in \Z^d} {B(k,1/4)}$. Fix $s=1/ 8$. For any $r\ge 1$, the set $B(0,r) \cap \partial F$ contains at least $Cr^d$ distinct  spheres of radius $1/4$ 
for some $C=C(d)>0$,  with the distance between any two of them being at least  $1/2>2s$. Hence, any covering of  $B(0,r) \cap \partial F$ by balls of radius $s$ requires at least $Cr^d$ balls. Letting $r\to \infty$  yields  ${\rm dim}_{\rm A}(\partial F)=d$.
}
\end{remark}

It is  known that the Hausdorff dimension of a subset of $\R^d$ does not exceed its Assouad dimension; see   \cite[Lemma 2.4.3]{Fr21}. Thus, since ${\rm dim}_{\rm A}(\Sigma)<d$, we have $m_d(\Sigma)=0$.
  Since $D$ is dense in $F$  and $\kappa$-fat and $m_d(\Sigma)=0$,
  there exists a constant $C>0$ such that
 \begin{align}\label{e:d-set}
  m_d(B_{F}(x,r)) = m_d(B_D(x,r)) \ge Cr^d \quad \text{ for all $x\in  {F}$ and $r\in (0, R_0)$.}
 \end{align}

	\begin{remark}\label{r:Assouad}
	{\rm	In Definition \ref{def:Assouad}, the restriction 
		 $r<\diam(E)$  may be removed. Indeed,  when $\diam(E)<\infty$, since $\overline E \subset \cup_{x\in E} B(x,\diam( E)/2) $ and $\overline E$ is compact, there are $x_1,\cdots,x_N \in  E$ such that $\overline E \subset  \cup_{i=1}^N B(x_i,\diam( E)/2)  $. Thus, for any $\lambda \in \mathfrak I(E)$ and all $0<s\le r$ with $r\ge \diam( E)$, one can see that $E= \cup_{i=1}^N (B(x_i,\diam( E)/2) \cap  E) $ can be covered by at most $CN(\diam( E)/ (2s))^\lambda \le CN(r/(2s))^\lambda$ balls of radius $s$ centered in $E$ if $s\le \diam(E)/2$, and at most $N \le N(r/s)^\lambda$   balls of radius $s$ centered in $E$  if $s\ge \diam(E)/2$. This shows that the restriction 
		 $r<\diam(E)$ in the definition of  $		{\rm dim}_{\rm A}( E)$ can be omitted.
	}
	\end{remark}

\begin{defn}\label{def:Assouadc}
	\rm Let  $E \subset \R^d$ be a non-empty subset. 
	The
	\textit{Assouad codimension} $	{\rm codim}_{\rm A}(E)$ of $E$ is the supremum of those $q\ge 0$ for which there exists a constant $C=C(q)\ge 1$ such that for all $x\in E$ and $0<s\le r< \diam(E)$,
	\begin{align*}
		\frac{m_d\left( \left\{ y \in B(x,r): \dist(y, E)<s\right\} \right) }{m_d(B(x,r))} \le C \bigg( \frac{s}{r}\bigg)^q.
	\end{align*}
\end{defn}

By  \cite[Lemma 3.4]{KLV13}, it holds  that ${\rm dim}_{\rm A}(E) + {\rm codim}_{\rm A}(E) = d$ for every non-empty subset $E\subset \R^d$. Consequently, we get ${\rm codim}_{\rm A}(\Sigma)= \gamma$.

\begin{prop}\label{p:Aikawa}
	For any $q \in [0,\gamma)$, there exists $C=C(q)>0$ such that
	\begin{align}\label{e:Aikawa}
		m_d \left( \left\{y \in B_{F}(x,r): \delta_\Sigma(y)<s\right\} \right)  \le C s^{q} (s\vee r)^{d-q} \quad \text{ for all $x\in {F}$ and $r,s>0$.} 
	\end{align}
\end{prop}
\pf If $\delta_\Sigma(x) \ge 2(s\vee r)$, then $ \left\{y \in B_{F}(x,r): \delta_\Sigma(y)<s\right\}=\emptyset$. Thus, it suffices to prove \eqref{e:Aikawa} for $\delta_\Sigma(x)<2(s\vee r)$.
 Suppose $\delta_\Sigma(x)<2(s\vee r)$ and let $Q\in \Sigma$ be such that $|x-Q|=\delta_\Sigma(x)$.
 
   If $s\vee r<3^{-1}\diam(\Sigma)$, then, since $q<{\rm codim}_{\rm A}(\Sigma)$,  we have
\begin{align*}
	m_d \left( \left\{y \in B_{F}(x,r): \delta_\Sigma(y)<s\right\} \right) & \le  m_d \left( \left\{y \in B_{F}(Q,3(s\vee r)): \delta_\Sigma(y)<s\right\} \right) \le c_1 s^q (s\vee r)^{d-q}.
\end{align*}

We now assume that $\diam(\Sigma )<\infty$ and $s \vee r\ge 3^{-1}\diam(\Sigma )$. 
\ Since $\Sigma$ is bounded, there are $N=N(\Sigma)\in \N$ and  $Q_1, \cdots, Q_N \in \Sigma $ such that $\{ y \in {\R^d}: \delta_\Sigma(y)\le 2^{-3}\diam(\Sigma )\} \subset \cup_{i=1}^N B(Q_i,2^{-2}\diam(\Sigma ) )$. 
 Thus,  since $q<\gamma\le d$ and \eqref{e:Aikawa} holds for $s\vee r<3^{-1}\diam(\Sigma )$, we get  that  for all   $0<s\le 8^{-1}\diam(\Sigma )$ and $r\ge 3^{-1}\diam(\Sigma )$,
\begin{align*}
	&m_{d}\left( \left\{ y \in B_{F}(x,r): \delta_\Sigma(y)<s \right\} \right) \le 	m_{d}\left( \left\{ y \in {F}: \delta_\Sigma(y)<s \right\} \right) \\
	&\le  \sum_{i=1}^N 	m_{d}\left( \left\{ y \in B_{F}(Q_i, 2^{-2}\diam(\Sigma ) ): \delta_\Sigma(y)<s \right\} \right) \le c_2s^{q} \diam(\Sigma )^{d-q} \le 3^{d-q} c_2 s^{q} r^{d-q} .
\end{align*}
On the other hand, for all $s>8^{-1}\diam(\Sigma )$ and $r>0$,  we have
\begin{align*}
	&m_{d}\left( \left\{ y \in B_{F}(x,r): \delta_\Sigma(y)<s \right\} \right) \le 	m_{d}\left( \left\{ y \in {F}:\delta_\Sigma(y)<s \right\} \right) \\
	&\le 	m_{d}\left( B(Q, \diam(\Sigma )  + s) \right)  \le c_3 s^d \le c_3 s^q (s\vee r)^{d-q}.
\end{align*}
The proof is complete. \qed 

For $a>0$, define
\begin{align*}
	D_a:=\left\{ x\in D: \delta_\Sigma(x)>a\right\}.
\end{align*}Recall that $D$ is $\kappa$-fat for some $\kappa \in (0,1)$. Define
\begin{align}\label{e:def-Lambda}
	\Lambda:=4+5\kappa^{-1}.
\end{align}

\begin{lem}\label{l:D-rho-nonempty}
	For any $z \in D$ and $\rho \in (0,\Lambda^{-1}R_0)$, we have	$(B_D(z,(\Lambda-1) \rho) \cap D_{2\rho}) \setminus B(z,2\rho) \neq \emptyset$.
\end{lem}
\pf Let $z \in D$ and $\rho \in (0,\Lambda^{-1}R_0)$.  Since $D$ is $\kappa$-fat with localization constant $R_0>5\kappa^{-1}\rho$, there exists  $y\in D_{5\rho}$ such that $B(y,5 \rho) \subset  B_D(z, 5\kappa^{-1}\rho)$. 
 Note that  $B(y, 3\rho) \subset B_D(z, (\Lambda-1)\rho) \cap  D_{2\rho}$.  Thus, if  $(B_D(z,(\Lambda-1) \rho) \cap D_{2\rho}) \setminus B(z,2\rho) = \emptyset$, then 
\begin{align*}
	m_d(B(z,2\rho)) \ge m_d(B_D(z,(\Lambda-1) \rho) \cap D_{2\rho})  \ge m_d( B(y, 3\rho)) = m_d(B(z,3\rho)),
\end{align*}
which is a contradiction. Hence, $(B_D(z,(\Lambda-1) \rho) \cap D_{2\rho}) \setminus B(z,2\rho) \neq \emptyset$. \qed

Lemma \ref{l:D-rho-nonempty}  implies in particular that $D_{\rho} \neq \emptyset$ for any $\rho \in (0, 2\Lambda^{-1}R_0)$.

\subsection{Some elementary estimates}\label{ss:prelim-2}

In this subsection, we present some  inequalities involving the blow-up function $\Phi$ and the jump kernel $J$. These results will be used in later sections.

Recall that 
$\Phi$
has 
upper Matuszewska index  $\overline \beta$ satisfying \eqref{e:beta-range}. Set
\begin{align}\label{e:beta1}
	\beta_1:=\frac{(\gamma \wedge \alpha)+ \overline\beta}{2} \in (\overline \beta, \gamma \wedge \alpha).
\end{align}
Then there exists a constant $C>1$ such that
\begin{align}\label{e:Phi-scaling}
	\frac{\Phi(r)}{\Phi(s)} \le  C\bigg( \frac{r}{s}\bigg)^{\beta_1} \quad \text{for all $0<s\le r$}.
\end{align} 
Since $\Phi$ is increasing, it follows that 
\begin{align}\label{e:Phi-scaling-monotone}
	\frac{\Phi(r)}{\Phi(s)} \le 1 + C\bigg( \frac{r}{s}\bigg)^{ \beta_1} \quad \text{for all $r,s>0$}.
\end{align} 

Note that 
\begin{align}
	 \frac{r\wedge A}{s\wedge A} \le \frac{r}{s} \quad \text{for all $A\in (0,\infty]$ and $0<s\le r$}.\label{e:without-A0}
\end{align}

\begin{lem}\label{l:Aikawa-integrals}
(i) For any $q\in [0,\gamma)$,	there exists $C=C(q)>0$ such that
\begin{align*}
	\int_{B_F(x,r)} \frac{dy}{(\delta_\Sigma(y) \wedge A)^{q}}   \le \frac{Cr^{d}}{((\delta_\Sigma(x) \vee r) \wedge A)^{q}}    \quad \text{ for all $A\in (0,\infty]$, $x\in  {F}$ and  $r>0$.}
\end{align*}

\noindent (ii) There exists $C>0$ such that
	\begin{align*}
		\int_{B_F(x,r)} \Phi\bigg( \frac{a}{\delta_\Sigma(y) \wedge A}\bigg)  dy \le Cr^{d}\Phi\bigg( \frac{a}{(\delta_\Sigma(x) \vee r) \wedge A}\bigg)    \quad \text{ for all $A\in (0,\infty]$, $x\in  {F}$ and $a,r>0$.}
	\end{align*}
\end{lem}
\pf (i)  If $\delta_\Sigma(x) \ge 2r$, then $\delta_\Sigma(y) \ge \delta_\Sigma(x)-r \ge \delta_\Sigma(x)/2$ for all $y \in B_F(x,r)$. Thus, in this case, we have
\begin{align*}
		\int_{B_F(x,r)} \frac{dy}{(\delta_\Sigma(y) \wedge A)^{q}} \le 	 \frac{2^q}{(\delta_\Sigma(x) \wedge A)^{q}} \int_{B_F(x,r)}  dy \le \frac{c_1r^d}{(\delta_\Sigma(x) \wedge A)^{q}}.
\end{align*}
 Suppose that $\delta_\Sigma(x)<2r$.  Then we have $\delta_\Sigma(y)<3r$ for all $y \in B_F(x,r)$. Thus, using \eqref{e:without-A0} and Proposition \ref{p:Aikawa} (with $q$ replaced by $(q+\gamma)/2$), we get that
 \begin{align*}
 	&\int_{B_F(x,r)} \frac{dy}{(\delta_\Sigma(y) \wedge A)^{q}}  \le \frac{(3r)^q}{	( (3r) \wedge A)^{q}}\int_{B_F(x,r)}  \frac{dy}{\delta_\Sigma(y)^{q}}\\
 	& \le \frac{(3r)^q}{	( (3r) \wedge A)^{q}} \sum_{n=0}^\infty  \frac{ 3^{nq}}{r^q} \int_{B_F(x,r): \delta_\Sigma(y) \in [3^{-n}r, 3^{1-n}r)}  dy\\
 	&\le \frac{c_2}{	( (3r) \wedge A)^{q}} \sum_{n=0}^\infty  3^{nq} (3^{-n}r)^{(q+\gamma)/2} (3r)^{d-(q+\gamma)/2}= \frac{c_3r^d}{	( (3r) \wedge A)^{q}} \le \frac{c_3r^d}{	( (\delta_\Sigma(x) \vee r) \wedge A)^{q}}.
 \end{align*}

\noindent (ii) Note that $\delta_\Sigma(y) \le \delta_\Sigma(x) + r \le 2(\delta_\Sigma(x) \vee r)$ for all  $y\in B_F(x,r)$. Thus, using \eqref{e:Phi-scaling} and the monotonicity of $\Phi$, we obtain for all $y \in B_F(x,r)$,
\begin{align*}
	\Phi\bigg( \frac{a}{\delta_\Sigma(y) \wedge A}\bigg) &\le c_4 \bigg(\frac{(2(\delta_\Sigma(x) \vee r)) \wedge A}{\delta_\Sigma(y) \wedge A}\bigg)^{ \beta_1 }\Phi\bigg( \frac{a}{(2(\delta_\Sigma(x) \vee r)) \wedge A}\bigg)\nn\\
	&\le 2^{\beta_1}c_4 \bigg(\frac{(\delta_\Sigma(x) \vee r) \wedge A}{\delta_\Sigma(y) \wedge A}\bigg)^{ \beta_1 }\Phi\bigg( \frac{a}{(\delta_\Sigma(x) \vee r) \wedge A}\bigg).
\end{align*}Hence, 
since $0 < \beta_1<\gamma$ by \eqref{e:beta1}, 
applying (i) with $q=\beta_1$, we obtain
	\begin{align*}
	\int_{B_F(x,r)} \Phi\bigg( \frac{a}{\delta_\Sigma(y) \wedge A}\bigg)  dy &\le 2^{\beta_1}c_4\Phi\bigg( \frac{a}{(\delta_\Sigma(x) \vee r) \wedge A}\bigg)	\int_{B_F(x,r)} \bigg(\frac{(\delta_\Sigma(x) \vee r) \wedge A}{\delta_\Sigma(y) \wedge A}\bigg)^{ \beta_1 }dy \nn\\
	&\le c_5r^d\Phi\bigg( \frac{a}{(\delta_\Sigma(x) \vee r) \wedge A}\bigg).
\end{align*}
\qed

\begin{lem}\label{l:jump-density-monotonicty}
	There exists $C>0$ such that for all $a, r>0$ and $A \in (0, \infty]$,  
	\begin{align}\label{e:jump-density-monotonicty-1}
		&	\sup_{R\ge r}	\frac{ 1 }{R^{d+\alpha}}\Phi\left( \frac{(R \wedge A)^2}{a}\right)\le 	\frac{ C }{r^{d+\alpha}}\Phi\left( \frac{(r\wedge A)^2}{a}\right).
	\end{align}
	Consequently, for all $r>0$, $x\in F$ and $y \in F\setminus B(x,r)$, we have
	\begin{align}\label{e:jump-density-monotonicty-2}
		J(x,y)\le  \frac{C}{r^{d+\alpha}} \Phi\bigg( \frac{(r\wedge A_0)^2}{(\delta_\Sigma(x) \wedge A_0)(\delta_\Sigma(y) \wedge A_0)}\bigg).
	\end{align}
\end{lem}
\pf  Using \eqref{e:Phi-scaling}  and \eqref{e:without-A0}, we get that for all $R\ge r$,
\begin{align*}
		\frac{ 1 }{R^{d+\alpha}}\Phi\left( \frac{(R \wedge A)^2}{a}\right)	&\le \frac{c_1}{R^{d+\alpha}} \bigg( \frac{R}{ r }\bigg)^{2\beta_1} \Phi\left( \frac{(r \wedge A)^2}{a }\right)\le \frac{c_1}{r^{d+\alpha}}\Phi\left( \frac{(r \wedge A)^2}{a }\right).
\end{align*}
We used $2\beta_1<2(\gamma\wedge \alpha) \le 2(d\wedge \alpha) \le d+\alpha$ in the last inequality. This proves \eqref{e:jump-density-monotonicty-1}. 

For \eqref{e:jump-density-monotonicty-2}, using \eqref{e:JJ} and  \eqref{e:jump-density-monotonicty-1}, we obtain for all $r>0$, $x\in F$ and $y \in F\setminus B(x,r)$, 
\begin{align*}
	J(x,y)&\le  \frac{c_2}{|x-y|^{d+\alpha}} \Phi\bigg( \frac{(|x-y|\wedge A_0)^2}{(\delta_\Sigma(x) \wedge A_0)(\delta_\Sigma(y) \wedge A_0)}\bigg)\\
	&\le \sup_{R\ge r}  \frac{c_2}{R^{d+\alpha}} \Phi\bigg( \frac{(R\wedge A_0)^2}{(\delta_\Sigma(x) \wedge A_0)(\delta_\Sigma(y) \wedge A_0)}\bigg)\le  \frac{c_3}{r^{d+\alpha}} \Phi\bigg( \frac{(r\wedge A_0)^2}{(\delta_\Sigma(x) \wedge A_0)(\delta_\Sigma(y) \wedge A_0)}\bigg).
\end{align*}
\qed

\begin{lem}\label{l:tail-estimate}
For any 
	$\theta >\ub$,  there exists $C=C(\theta)>0$  such that 
	 for all $x \in  F$, $A\in (0,\infty]$, 
and $r, b>0$,
	\begin{align}
	\label{e:tail-estimate0}
		\int_{F\setminus B(x,r)}	 \frac{ 1 }{|x-y|^{d+\theta}}\Phi\left( \frac{(|x-y| \wedge A)^2}{b(\delta_\Sigma(y) \wedge A) }\right)dy \le \frac{C}{r^{\theta}} \Phi\left( \frac{r \wedge A}{b}\right).
	\end{align}
Consequently,	there exists $C>0$  such that 
	for all 
	$x\in F$ and $r>0$,
	\begin{align}
	\label{e:tail-estimate}
		\int_{F\setminus B(x,r)}	J(x,y)dy \le \frac{C}{r^{\alpha}} \Phi\left( \frac{r \wedge A_0}{\delta_\Sigma(x) \wedge A_0}\right).
	\end{align}
\end{lem}
\pf  Choose $\beta \in (\ub, \theta)$ so that 
\begin{align}\label{e:upper-scaling0}
	\Phi(r)/\Phi(s) \le C(r/s)^{ \beta} \quad \text{for all $0<s\le r$}.
\end{align}
Let $x\in F$ and $r>0$. We consider the following two cases separately.

\smallskip

\noindent \textit{Case 1:} 
 $r<A$.
Let 
$N:=\min\{n\ge 2: 2^{n-1}r \ge A\}$.  
For all $1\le n<N$, 
using  Lemma \ref{l:Aikawa-integrals}(ii) (with $a=(2^nr)^2/b$)  in the second inequality below, 
the monotonicity of $\Phi$ in the third, and \eqref{e:upper-scaling0} in the fourth,  
 we get that
\begin{align}\label{e:tail-estimate-1}
	&	\int_{B_F(x,2^{n}r)\setminus B(x,2^{n-1}r)} 
	\frac{ 1 }{|x-y|^{d+\theta}}\Phi\left( \frac{(|x-y| \wedge A)^2}{b(\delta_\Sigma(y) \wedge A) }\right)
	dy\nn\\
	& \le 
	 \frac{2^{d+\theta}C_1}{2^{n(d+\theta)} r^{d+\theta}} 
		\int_{B_F(x,2^{n}r)\setminus B(x,2^{n-1}r)}
		\Phi\bigg(\frac{(2^n r)^2}{b(\delta_\Sigma(y)\wedge A)}\bigg) 
		dy\nn\\
	& \le 
	\frac{c_1 }{2^{n\theta} r^{\theta}} 
	\Phi\bigg(\frac{(2^n r)^2}{b((\delta_\Sigma(x) \vee 2^n r)\wedge A)}\bigg) \le  
	\frac{c_1 }{2^{n\theta} r^{\theta}} 
	\Phi\bigg(\frac{2^{n+1} r}{b}\bigg)  \le  
	\frac{c_2 }{2^{n(\theta-
	\beta)} r^{\theta}} 
	\Phi\bigg(\frac{ r}{b}\bigg).
\end{align}
Further, using 
Lemma \ref{l:Aikawa-integrals}(ii), 
we obtain for all $n\ge N$, 
\begin{align}\label{e:tail-estimate-2}
	&	\int_{B_F(x,2^{n}r)\setminus B(x,2^{n-1}r)} 
	\frac{ 1 }{|x-y|^{d+\theta}}\Phi\left( \frac{(|x-y| \wedge A)^2}{b(\delta_\Sigma(y) \wedge A) }\right)
	dy\nn\\
	& \le 
	 \frac{2^{d+\theta}C_1}{2^{n(d+\theta)} r^{d+\theta}} 
		\int_{B_F(x,2^{n}r)\setminus B(x,2^{n-1}r)}
		\Phi\bigg(\frac{A^2}{b(\delta_\Sigma(y)\wedge A)}\bigg) dy\nn\\
	& \le 
	 \frac{c_3 }{2^{n\theta} r^{\theta}} 
	\Phi\bigg(\frac{A^2}{b((\delta_\Sigma(x) \vee 2^n r)\wedge A)}\bigg) =
	\frac{c_3 }{2^{n\theta} r^{\theta}} 
	\Phi\bigg(\frac{A}{b}\bigg).
\end{align}
Combining \eqref{e:tail-estimate-1} and \eqref{e:tail-estimate-2}, 
since  $\theta-\beta>0$, 
we arrive at 
\begin{align}\label{e:tail-estimate-3}
	&	\int_{F\setminus B(x,r)} 
	\frac{ 1 }{|x-y|^{d+\theta}}\Phi\left( \frac{(|x-y| \wedge A)^2}{b(\delta_\Sigma(y) \wedge A) }\right)
	dy\nn\\ &\le
	\sum_{n= 1}^\infty \int_{B_F(x,2^{n}r)\setminus B(x,2^{n-1}r)} 
	\frac{ 1 }{|x-y|^{d+\theta}}\Phi\left( \frac{(|x-y| \wedge A)^2}{b(\delta_\Sigma(y) \wedge A) }\right)
	dy\nn\\
	&\le \sum_{n=1}^{N-1} 
	\frac{c_2 }{2^{n(\theta-
\beta	)} r^{\theta}} \Phi\bigg(\frac{ r}{b}\bigg)   + \sum_{n=N}^\infty  \frac{c_3 }{2^{n\theta} r^{\theta}} \Phi\bigg(\frac{A}{b}\bigg) \le 
	\frac{c_{4}}{r^\theta} \Phi\bigg(\frac{r}{b}\bigg) + \frac{c_{5}}{2^{N\theta}
	 r^\theta } \Phi\bigg(\frac{A}{b }\bigg).
\end{align}
Applying 
 \eqref{e:upper-scaling0},  since   $2^Nr\ge A>r$ and $\beta<\theta$, we have that
\begin{align*}
\frac{1}{2^{N\theta}r^\theta } \Phi\bigg(\frac{A}{b }\bigg)  \le  	\frac{1}{A^\theta } \Phi\bigg(\frac{A}{b }\bigg)  \le \frac{c_{6}}{r^\theta}  \Phi\bigg(\frac{r}{b}\bigg).
\end{align*}
Combining this with \eqref{e:tail-estimate-3}, the result follows.

\smallskip

\noindent \textit{Case 2:} 
$r\ge A$. 
Note that \eqref{e:tail-estimate-2} holds for all $n\ge 1$ by 
Lemma \ref{l:Aikawa-integrals}(ii). 
Thus,  we get
\begin{align*}
	&	\int_{F\setminus B(x,r)} 
	\frac{ 1 }{|x-y|^{d+\theta}}\Phi\left( \frac{(|x-y| \wedge A)^2}{b(\delta_\Sigma(y) \wedge A) }\right)
	dy\le \sum_{n= 1}^\infty \frac{c_3}{
		2^{n\theta } 
	r^{\theta}} \Phi\bigg(\frac{A}{b}\bigg)= \frac{c_{7}}{r^{\theta} }\Phi\bigg(\frac{A}{b}\bigg).
\end{align*}
The proof of \eqref{e:tail-estimate0} is complete. 
Now \eqref{e:tail-estimate} follows immediately from \eqref{e:tail-estimate0} and  \eqref{e:JJ}. \qed 

\begin{lem}\label{l:lower-tail-estimate}
	There exists $C>0$ such that for all $z\in D$ and  $\rho \in (0, \Lambda^{-1}  R_0)$,
	\begin{align*}
		\int_{(B_F(z,\Lambda\rho)\cap D_\rho) \setminus B(z,\rho)} J(z,w)dw \ge \frac{C}{\rho^\alpha}\Phi\bigg( \frac{\rho \wedge A_0}{\delta_\Sigma(z) \wedge A_0}\bigg).
	\end{align*}
\end{lem}
\pf Let $z\in D$ and  $\rho \in (0, \Lambda^{-1}  R_0)$. By Lemma \ref{l:D-rho-nonempty}, there exists  $x\in (B_D(z,(\Lambda-1) \rho) \cap D_{2\rho}) \setminus B(z,2\rho)$. Note that $B(x, \rho) \subset (B_D(z,\Lambda\rho) \cap D_\rho) \setminus B(z,\rho)$. Hence, for all $w\in B(x,\rho)$, we have $\rho \le |z-w|<\Lambda\rho$ and $\delta_\Sigma(w)\le \delta_\Sigma(z) + \Lambda\rho$. 
Thus, using \As, \eqref{e:Phi-scaling} and the fact that  $1=\Phi(0)\le \Phi(1)$, we get that for all $w\in B(x,\rho)$,
\begin{align*}
	J(z,w) &\ge \frac{c_1}{(\Lambda\rho)^{d+\alpha}} \Phi \bigg( \frac{(\rho \wedge A_0)^2}{(\delta_\Sigma(z) \wedge A_0) ((\delta_\Sigma(z)  + \Lambda\rho) \wedge A_0)}\bigg)\\
	&\ge \frac{c_2}{\rho^{d+\alpha}} \begin{cases}
		\Phi ( (\rho \wedge A_0)/ (\delta_\Sigma(z) \wedge A_0)) &\mbox{ if $\delta_\Sigma(z) \le \rho$},\\
		1&\mbox{ if $\delta_\Sigma(z) > \rho$}
	\end{cases}\\
	&\ge \frac{c_2}{\Phi(1)\rho^{d+\alpha}}\Phi \bigg( \frac{\rho \wedge A_0}{\delta_\Sigma(z) \wedge A_0}\bigg).
\end{align*} 
Using this, we conclude that 
\begin{align*}
	\int_{(B_F(z,\Lambda\rho)\cap D_\rho) \setminus B(z,\rho)} J(z,w)dw  \ge \frac{c_2}{\Phi(1)\rho^{d+\alpha}}\Phi \bigg( \frac{\rho \wedge A_0}{\delta_\Sigma(z) \wedge A_0}\bigg)	\int_{B(x,\rho)}dw = \frac{c_3}{\rho^{\alpha}}\Phi \bigg( \frac{\rho \wedge A_0}{\delta_\Sigma(z) \wedge A_0}\bigg).
\end{align*}
\qed

\begin{lem}\label{l:tail-estimate-2}
	(i)		There exists $C>0$  such that 
	for all $x,z\in F$ and $r,s>0$,
	\begin{align*}
		\int_{B_F(z,s)\setminus B(x,r)}	J(x,y)dy \le \frac{Cs^d}{r^{d+\alpha}} \Phi\left( \frac{(r \wedge A_0)^2}{(\delta_\Sigma(x) \wedge A_0) ( (\delta_\Sigma(z) \vee s) \wedge A_0) }\right).
	\end{align*}
	
	\noindent (ii) 	There exists $C>0$  such that 
	for all 
	 $a>0$, 
	$A\in (0, \infty]$, 
	$x,z\in F$ and $r>0$ and $0<u\le s$,
	\begin{align*}
		&\int_{B_F(z,s)\setminus B(x,r)}	 \Phi\bigg(\frac{a}{(\delta_\Sigma(y) \vee u) \wedge A}\bigg) J(x,y)dy \\
		&\le  \frac{Cs^d}{r^{d+\alpha}} \bigg[1 + \bigg( \frac{a}{s \wedge A}\bigg)^{\beta_1}\bigg]\Phi\left( \frac{(r \wedge A_0)^2}{(\delta_\Sigma(x) \wedge A_0) (u \wedge A_0) }\right).
	\end{align*}
\end{lem}
\pf	(i) Applying Lemmas \ref{l:jump-density-monotonicty} and  \ref{l:Aikawa-integrals}(ii), we obtain
\begin{align*}
	\int_{B_F(z,s)\setminus B(x,r)}	J(x,y)dy &\le \frac{c_1}{r^{d+\alpha}} 	\int_{B_F(z,s)} \Phi\bigg( \frac{(r\wedge A_0)^2}{(\delta_\Sigma(x) \wedge A_0)(\delta_\Sigma(y) \wedge A_0)}\bigg)dy\\
	&\le \frac{c_2s^d}{r^{d+\alpha}} \Phi\bigg( \frac{(r\wedge A_0)^2}{(\delta_\Sigma(x) \wedge A_0)((\delta_\Sigma(z) \vee s) \wedge A_0)}\bigg).
\end{align*}
(ii) 
 By Lemma \ref{l:jump-density-monotonicty}, we have
\begin{align*}
	&\int_{B_F(z,s)\setminus B(x,r)}	 \Phi\bigg(\frac{a}{(\delta_\Sigma(y) \vee u) \wedge A}\bigg) J(x,y)dy\\
	& \le \frac{c_3}{r^{d+\alpha}}\int_{B_F(z,s)\setminus B(x,r) : \delta_\Sigma(y)\ge u}	 \Phi\bigg(\frac{a}{\delta_\Sigma(y) \wedge A}\bigg) \Phi\bigg( \frac{(r\wedge A_0)^2}{(\delta_\Sigma(x) \wedge A_0)(\delta_\Sigma(y) \wedge A_0)}\bigg)dy\\
	&\quad +  \frac{c_3}{r^{d+\alpha}} \int_{B_F(z,s)\setminus B(x,r) : \delta_\Sigma(y)<u}	 \Phi\bigg(\frac{a}{u \wedge A}\bigg) \Phi\bigg( \frac{(r\wedge A_0)^2}{(\delta_\Sigma(x) \wedge A_0)(\delta_\Sigma(y) \wedge A_0)}\bigg)dy\\
	&=: \frac{c_3}{r^{d+\alpha}} \left( I_1+I_2\right).
\end{align*}
Using Lemma \ref{l:Aikawa-integrals}(ii) and \eqref{e:Phi-scaling-monotone}, we get that
\begin{align*}
	I_1 &\le \Phi\bigg( \frac{(r\wedge A_0)^2}{(\delta_\Sigma(x) \wedge A_0)(u \wedge A_0)}\bigg) \int_{B_F(z,s)\setminus B(x,r) }	 \Phi\bigg(\frac{a}{\delta_\Sigma(y)  \wedge A}\bigg) dy\\
	&\le c_4 s^d \Phi\bigg(\frac{a}{s \wedge A}\bigg) \Phi\bigg( \frac{(r\wedge A_0)^2}{(\delta_\Sigma(x) \wedge A_0)(u \wedge A_0)}\bigg)  \\
	&\le c_5 s^d \bigg[1 + \bigg( \frac{a}{s \wedge A}\bigg)^{\beta_1}\bigg] \Phi\bigg( \frac{(r\wedge A_0)^2}{(\delta_\Sigma(x) \wedge A_0)(u \wedge A_0)}\bigg) .
\end{align*}
For $I_2$, 
since $\beta_1<\gamma$ by \eqref{e:beta1}, using Proposition \ref{p:Aikawa} with $q=(\beta_1+\gamma)/2 
<\gamma
$  in the second inequality below, \eqref{e:Phi-scaling} and \eqref{e:Phi-scaling-monotone} in the third, and  \eqref{e:without-A0} in the fourth,  we obtain
\begin{align*}
	I_2& \le c_6\Phi\bigg(\frac{a}{u   \wedge A}\bigg)  \sum_{m=1}^\infty \Phi\bigg( \frac{(r\wedge A_0)^2}{(\delta_\Sigma(x) \wedge A_0)(2^{-m}u \wedge A_0)}\bigg) \int_{B_F(z,s)\setminus B(x,r) : \delta_\Sigma(y)\in [ 2^{-m} u, 2^{1-m}u)}	 dy\\
	& \le c_7  \Phi\bigg(\frac{a}{u   \wedge A}\bigg) \sum_{m=1}^\infty (2^{1-m}u)^{(\beta_1+\gamma)/2} s^{d-(\beta_1+\gamma)/2} \Phi\bigg( \frac{(r\wedge A_0)^2}{(\delta_\Sigma(x) \wedge A_0)(2^{-m}u \wedge A_0)}\bigg) \\
	& \le c_8u^{(\beta_1+\gamma)/2} s^{d-(\beta_1+\gamma)/2}\bigg[1 + \bigg( \frac{a}{u \wedge A}\bigg)^{\beta_1}\bigg]    \Phi\bigg( \frac{(r\wedge A_0)^2}{(\delta_\Sigma(x) \wedge A_0)(u \wedge A_0)}\bigg)\sum_{m=1}^\infty 2^{-m(\beta_1/2 + \gamma/2 - \beta_1) }   \\
	&= c_9 u^{(\beta_1+\gamma)/2} s^{d-(\beta_1+\gamma)/2}\bigg[1 + \bigg( \frac{a}{u \wedge A}\bigg)^{\beta_1}\bigg]    \Phi\bigg( \frac{(r\wedge A_0)^2}{(\delta_\Sigma(x) \wedge A_0)(u \wedge A_0)}\bigg)\\
	&\le c_9 u^{(\beta_1+\gamma)/2} s^{d-(\beta_1+\gamma)/2} \bigg( \frac{s}{u}\bigg) ^{\beta_1 }\bigg[1 + \bigg( \frac{a}{s \wedge A}\bigg)^{\beta_1}\bigg]    \Phi\bigg( \frac{(r\wedge A_0)^2}{(\delta_\Sigma(x) \wedge A_0)(u \wedge A_0)}\bigg) \\
	&\le  c_9  s^d \bigg[1 + \bigg( \frac{a}{s \wedge A}\bigg)^{\beta_1}\bigg]    \Phi\bigg( \frac{(r\wedge A_0)^2}{(\delta_\Sigma(x) \wedge A_0)(u \wedge A_0)}\bigg).
\end{align*}
Here we also used  $\beta_1<\gamma$ in the equality and the last inequality. The proof is complete. 	\qed

\section{Preliminary analysis of Dirichlet forms}\label{s:prel-anal-df} 

In this section, we study the Dirichlet form $\sE$ and its truncated versions. In Subsection \ref{ss:Hardy} we prove a Hardy-type inequality and present some relationship between $\sE$ and
its truncated version. 
The fact $\text{dim}_{\rm A}(\Sigma )<d$  will play a crucial role in the proof of the Hardy-type inequality in Proposition \ref{p:Hardy-fractional-Laplacian}.
 In Subsection \ref{ss:reg}, we prove the parabolic H\"{o}lder regularity and near-diagonal heat kernel bounds. We start with the following result. 

\begin{prop}\label{p:regular-Dirichlet-form}
	$\sE(u,u)<\infty$ for all $u \in {\rm Lip}_c({F})$.
\end{prop}
\pf The proof is essentially the same as that of \cite[Proposition 5.1]{CK25+}, and we give the proof for the reader's convenience.    Let  $u \in {\rm Lip}_c({F})$, $x_0\in F$ and $r>0$ 
be such that $\text{supp}[u] \subset B_{F}(x_0,r)$.   
Since $u \in {\rm Lip}_c({F})$, there exists $c_1>0$ such that $|u(x)-u(y)|\le c_1(|x-y| \wedge 1)$ for all $x,y \in F$. Hence, by Lemma \ref{l:tail-estimate}, we have  for all $x\in B_{F}(x_0,r)$,
\begin{align*}
	&\int_F (u(x)-u(y))^2 J(x,y)dy \le c_1^2\int_F  ( |x-y| \wedge 1)^2J(x,y)dy\\
	& \le c_1^2\sum_{n=0 }^\infty \int_{B_F(x,  2^{-n}) \setminus B(x,2^{-n-1})} 
	2^{-2n} J(x,y)dy  + c_1^2 \int_{F\setminus B( x, 1)} J(x,y)dy\\
	&\le  c_2\sum_{n=0}^\infty 
	2^{-2n+(n+1)\alpha} \Phi\bigg( \frac{2^{-n-1} \wedge A_0}{\delta_\Sigma(x) \wedge A_0 }\bigg) + c_3 \Phi\bigg( \frac{1 \wedge A_0}{\delta_\Sigma(x) \wedge A_0 }\bigg) \\
	&\le  c_4\Phi\bigg( \frac{1 \wedge A_0}{\delta_\Sigma(x) \wedge A_0}\bigg) \sum_{n=0}^\infty 2^{-n(2-\alpha)}  = c_5\Phi\bigg( \frac{1 \wedge A_0}{\delta_\Sigma(x) \wedge A_0}\bigg).
\end{align*}
It follows that 
\begin{align*}
	\sE(u,u) &\le 2  \int_{B_F(x_0,r)} \int_F (u(x)-u(y))^2 J(x,y) dy\,dx\le 
	c_5 \int_{B_F(x_0,r)} \Phi\bigg( \frac{1 \wedge A_0}{\delta_\Sigma(x) \wedge A_0}\bigg) dx  <\infty,
\end{align*}
where we used Lemma \ref{l:Aikawa-integrals}(ii)  in the last inequality. \qed

Let $\sF$ be the $\sE_1$-closure of ${\rm Lip}_c(F)$ in $L^2(F)$, which is well defined by Proposition \ref{p:regular-Dirichlet-form}. Then $(\sE,\sF)$ is a regular Dirichlet form on $L^2(F)$. 
Let $(P_t)_{t\ge 0}$ be the semigroup associated with $(\sE,\sF)$. Since $(\sE,\sF)$ is a regular Dirichlet form on $L^2(F)$, there exists a symmetric Hunt process  $X=(X_t, t\ge 0 ; \P^x, x\in F\setminus \sN)$ properly associated with $(\sE,\sF)$ in the sense that for every $f\in L^2(F)$ and $t>0$, the map $x\mapsto \E^x f(X_t)$ is a quasi-continuous version of $P_tf$. Here $\sN$ is a properly exceptional set for $X$.

\subsection{Hardy-type inequality and  truncated Dirichlet forms}\label{ss:Hardy}

Using the results of \cite{DV14}, we can get the following  fractional Hardy-type inequality with localization constant $R_0$. When $D$ is $\kappa$-fat with $R_0=\diam(D)$, this was established in  \cite[Proposition 4.8]{CK25+}. 
\begin{prop}\label{p:Hardy-fractional-Laplacian}
	For any $q \in (0,\gamma \wedge \alpha)$, there exists  $C=C(q)>0$ such that for all $x\in F$, $r\in (0, R_0)$ and $u \in C_c(F) \cap \sF$,
	\begin{align}\label{e:Hardy-fractional-Laplacian}
		r^q\!	\int_{B_F(x,r)} \!\frac{u(y)^2}{\delta_\Sigma(y)^{q}}  dy \le C \bigg( r^\alpha\! \int_{B_F(x,2r)\times B_F(x,2r)}\! \frac{(u(z)-u(y))^2}{|z-y|^{d+\alpha}} dzdy + 	\int_{B_F(x,2r)} \!u(y)^2  dy \bigg).
	\end{align}
\end{prop}
\pf Let $x \in F$, $r\in (0, R_0)$ and $u\in C_c(F)\cap \sF$. If $\delta_\Sigma(x)\ge 2r$, then, since $\delta_\Sigma(y)\ge r$ for all $y \in B_F(x,r)$, we have
\begin{align*}
	r^q	\int_{B_F(x,r)} \frac{u(y)^2}{\delta_\Sigma(y)^{q}}  dy \le \int_{B_F(x,r)} u(y)^2dy,
\end{align*}
implying \eqref{e:Hardy-fractional-Laplacian}.

 Suppose  $\delta_\Sigma(x)<2r$. Set $M:=\R^d \setminus (B(x,2r) \cap F^c)$ and let $\wt D:=M \setminus (B(x,2r) \cap \Sigma )$. Since $\delta_\Sigma(x)<2r$, $\wt D$ is a proper open subset of $M$ in the relative topology. Let $\partial_M \wt D:=B(x,2r) \cap \Sigma $ and $\wt q:= (q + (\gamma \wedge \alpha))/2$.  We claim that the following two conditions hold:

\medskip

(T1) There exists  $c_1>0$ such that for all $0<a<b$ and $Q \in \partial_M \wt D$, there is a cover of $B(Q,b) \cap \partial_M \wt D$ by  balls $B(Q_j, a) \cap M$ with $Q_j\in  \partial_M \wt D$, $j=1, \cdots, N$, where $N\le c_1(b/a)^{d-\wt q}$.

(T2) For any $Q \in \partial_M \wt D$ and $a>0$, there exists $z\in B(Q,a)$ such that $B(z, \kappa a/5) \subset \wt D$.

\medskip

For (T1), since $d-\wt q>d-\gamma= \text{dim}_A(\Sigma )$, 
by Remark \ref{r:Assouad}, 
there exists $c_2>0$ such that for all $0<a<b$ and $Q \in \partial_M \wt D$, there is a cover of $B(Q,b) \cap \Sigma $ by  balls $B(Q_j, a/2)$ with $Q_j\in  \partial  D$, $j=1, \cdots, N$, where $N\le c_2(2b/a)^{d-\wt q}$. Let $\mathfrak I:=\{j : Q_j \in \partial_M\wt D\}$. If $B(Q,b) \cap \partial_M\wt D \subset \cup_{j\in \mathfrak I} (B(Q_j, a/2)\cap M)$, then we are done.   Assume $(B(Q,b) \cap \partial_M\wt D) \setminus \cup_{j\in \mathfrak I} (B(Q_j, a/2)\cap M) \neq \emptyset$.  For any $y \in (B(Q,b) \cap \partial_M\wt D) \setminus \cup_{j\in \mathfrak I} (B(Q_j, a/2)\cap M) $, since $y \in B(Q,b) \cap  \Sigma $, there exists $j_y\in \{1,\cdots, N\}\setminus \mathfrak I$ such that  $y \in B(Q_{j_y}, a/2)$. Let $\mathfrak I':=\{j \notin \mathfrak I : j = j_y \text{ for some $y \in  (B(Q,b) \cap \partial_M\wt D) \setminus \cup_{j\in \mathfrak I } (B(Q_j, a/2)\cap M)$}\}$. 
For each $j\in \mathfrak I'$, let us fix $y(j) \in  (B(Q,b) \cap \partial_M\wt D) \setminus \cup_{j\in \mathfrak I} (B(Q_j, a/2)\cap M)$ such that $j=j_{y(j)}$. For every $j\in \mathfrak I'$, it holds that $B(Q_{j},a/2) \subset B(y(j), a)$. Hence, we have
\begin{align*}
	B(Q,b) \cap \partial_M \wt D =	B(Q,b) \cap \partial_M \wt D  \cap M &\subset \cup_{j\in \mathfrak I} (B(Q_j, a/2)\cap M) \, \cup \, \cup_{j\in \mathfrak I'} (B(Q_j, a/2)\cap M)\\&\subset \cup_{j\in\mathfrak I} (B(Q_j, a)\cap M) \, \cup \, \cup_{j\in \mathfrak I'} (B(y(j), a)\cap M).
\end{align*}
Since  
 $|\mathfrak I| + |\mathfrak I'|\le N\le 2^{d-q}c_2(b/a)^{d-\wt q}$, this proves (T1).

For (T2),  
if $a<5R_0$, then, since $D$ is $\kappa$-fat with localization constant $R_0$,  there exists $z\in B_D(Q,a/5)$ such that  $B(z,\kappa a/5) \subset B_D(Q,a/5) \subset \wt D$.   For
$a\ge 5R_0>5r$ and $Q\in \partial_M \wt D \subset B(x,2r)$, there exists $z \in B(Q,a)$ such that 
$|z-x| = |x-Q| + 3a/5$.
Since $|z-x| \ge 
3a/5  > \kappa a/5 + 2r$, we get that $B(z, \kappa a/5) \subset \R^d \setminus B(x, 2r) \subset \wt D$. This completes the proof for (T2).

Using (T1) and (T2), we can follow the argument in \cite[Proposition 10]{DV14} to deduce that  $\wt D \subset M$ satisfies  condition DC$(2,d-\wt q,d)$ of \cite[Definition 3.6]{DV14}, with  constants independent of $x$ and $r$. Applying \cite[Theorem 5]{DV14} (with $\phi(t)=t^q$, $\eta = q$, $\gamma=d-\wt q$ and  $p=2$),  we get that there exist constants $c_3,c_4>0$  independent of $x$ and $r$ such that
\begin{align}\label{e:Hardy-wtD}
	\int_{\wt D} \frac{f(y)^2}{\dist(y,\partial_M \wt D)^q}dy &\le c_3 \int_{\wt D} \int_{\wt D \cap B(y, c_4\, \dist(y,\partial_M \wt D))} \frac{(f(y)-f(z))^2}{\dist(y,\partial_M \wt D)^{d+q}}dz\,dy,
\end{align}
for all Borel functions $f$ on $\wt D$ 
for which the left-hand side is finite.

Define $\vp(y):= ((2-(|x-y|/r)) \wedge 1) \vee 0.$
Since $|\vp(y)-\vp(z)| \le (|y-z|/r) \wedge 1$ for all $y,z\in \R^d$, 
we have for all $y\in \R^d$ that 
\begin{align}\label{e:cut-off}
	\int_{\R^d} \frac{(\vp(y)-\vp(z))^2}{|y-z|^{d+q}} dz \le \frac{1}{r^2}\int_{B(y,r)}  \frac{dz}{|y-z|^{d+q-2}}  + \int_{\R^d\setminus B(y,r)} \frac{dz}{|y-z|^{d+q}} = \frac{c_5}{r^q}.
\end{align}
For any $y\in B_D(x,r)$, we have either $\delta_\Sigma(y) =  \dist(y,\partial_M \wt D)$ or $\delta_\Sigma(y) = \dist(y, \Sigma  \setminus B(x,2r))\ge r$. Thus, for any $y\in B_D(x,r)$,  since $ \dist(y,\partial_M \wt D) \le  \dist(x,\partial_M \wt D) + r <3r$, we have  $\delta_\Sigma(y) \ge  \dist(y,\partial_M \wt D)/3$. Using this, the fact that  $\vp=1$ on $B(x,r)$, and \eqref{e:Hardy-wtD}, we obtain
\begin{align*}
	&	r^q	\int_{B_D(x,r)} \frac{u(y)^2}{\delta_\Sigma(y)^{q}}  dy \le  3^qr^q	\int_{B_D(x,r)} \frac{u(y)^2}{\dist(y,\partial_M \wt D)^q}dy\le   3^qr^q	\int_{\wt D} \frac{(u(y) \vp(y))^2}{\dist(y,\partial_M \wt D)^q}dy  \\
	&  \le c_6 r^q \int_{\wt D} \int_{\wt D \cap B(y, c_4\, \dist(y,\partial_M \wt D))} \frac{(u(y)\vp(y)-u(z)\vp(z))^2}{\dist(y,\partial_M \wt D)^{d+q}}dz\,dy\\
	&  \le c_4^{d+q}c_6 r^q \int_{\wt D} \int_{\wt D} \frac{(u(y)\vp(y)-u(z)\vp(z))^2}{|y-z|^{d+q}}dz\,dy=:c_4^{d+q}c_6I.
\end{align*}
Since $\wt D \cap B(x,2r) = B_D(x,2r)$,  $\vp=0$ on $B(x,2r)^c$ and  $\vp^2\le 1$ on $\wt D$,
using the Cauchy-Schwarz inequality and \eqref{e:cut-off},   we get that
\begin{align*}
	I&= r^q \int_{B_D(x,2r)}\int_{B_D(x,2r)}\! \! \!\!\frac{(u(y)\vp(y)-u(z)\vp(z))^2}{|y-z|^{d+q}}dz\,dy 
	+ 2r^q\int_{B_D(x,2r)}\int_{\wt D \setminus B(x,2r)} \! \!\frac{(u(y)\vp(y))^2}{|y-z|^{d+q}}dz\,dy\\
	&\le  2 r^q\int_{B_D(x,2r)} \int_{B_D(x,2r)} \frac{\vp(z)^2(u(y)-u(z))^2}{|y-z|^{d+q}}dz\,dy \\
	&\quad + 2r^q\int_{B_D(x,2r)} \int_{B_D(x,2r)} \frac{u(y)^2(\vp(y)-\vp(z))^2}{|y-z|^{d+q}}dz\,dy \\
	&\quad + 2r^q\int_{B_D(x,2r)}\int_{\wt D \setminus B(x,2r)} \frac{u(y)^2(\vp(y)-\vp(z))^2}{|y-z|^{d+q}}dz\,dy \\
	&\le  2r^q\int_{B_D(x,2r)} \int_{B_D(x,2r)} \frac{(u(y)-u(z))^2}{|y-z|^{d+q}}dz\,dy  + 
	4c_5 \int_{B_D(x,2r)} u(y)^2 dy\\
	&\le  2^{1+2(\alpha-q)}r^\alpha\int_{B_D(x,2r)} \int_{B_D(x,2r)} \frac{(u(y)-u(z))^2}{|y-z|^{d+\alpha}}dz\,dy  + 4c_5 \int_{B_D(x,2r)} u(y)^2 dy,
\end{align*}
where we used $q<\alpha$ in the last inequality. 
Since $m_d(\Sigma)=0$, this yields the desired result.
\qed

By \eqref{e:J-standard-lower},   Proposition \ref{p:Hardy-fractional-Laplacian}  immediately gives
\begin{cor}\label{c:Hardy}
	For any $q \in (0,\gamma \wedge \alpha)$, there exists  $C=C(q)>0$ such that for all $x\in F$, $r\in (0, R_0)$ and $u \in C_c(F)\cap \sF$,
	\begin{align*}
		r^q	\int_{B_F(x,r)} \frac{u(y)^2}{\delta_\Sigma(y)^{q}}  dy \le C \bigg( r^\alpha \int_{B_F(x,2r)\times B_F(x,2r)} (u(z)-u(y))^2J(z,y) dzdy + 	\int_{B_F(x,2r)} u(y)^2  dy \bigg).
	\end{align*}
\end{cor}

Recall that we write $\overline R:=\diam(F)$.
For $\rho\in (0, \overline R]$, define
\begin{align*}
	J^{(\rho)}(x,y) :=  \1_{\{
		|x-y| < \rho\}} J(x,y).
\end{align*} 
The \textit{$\rho$-truncated  form} $\sE^{(\rho)}$ of $\sE$  is defined by
\begin{align}\label{def:truncated-forms}
	\sE^{(\rho)}(u,v) = 
	\frac12 \int_{F\times F} (u(x)-u(y))(v(x)-v(y))  \,J^{(\rho)}(x,y)dxdy.
\end{align}
Note that $\sE^{(\overline R)}=\sE$.

\begin{prop}\label{p:truncation-equivalent}
	There exists $C>0$ such that for all $\rho\in (0,\overline R)$ and $u \in L^2(F)\cap \sF$,
	\begin{align*}
		\sE(u,u) \le C \left(  \sE^{(\rho)}(u,u) + (\rho \wedge R_0)^{-\alpha} \lVert u\rVert_{L^2(F)}^2\right) .
	\end{align*}
\end{prop}
\pf Since $\sE^{(\rho)}(u,u)\le \sE^{(\rho')}(u,u)$ for all $0<\rho<\rho'<\overline R$ and $u \in L^2(F) \cap \sF$, it suffices to prove 
the result  for $\rho<R_0$.
Let  $\rho\in (0,R_0)$. By Vitali's covering lemma, there exists a family of open  balls $\{B(z_i,\rho/8)\}_{i\ge 1}$ such that $\R^d=\cup_{i\ge 1} B(z_i,\rho/8)$ and $B(z_i,\rho/40) \cap B(z_j,\rho/40)=\emptyset$ for all $i\neq j$. It follows that
$\sum_{i\ge 1} \1_{B(z_i,\rho)}\le N$ on $\R^d$ for some constant $N=N(d)\ge 1$.  
Let $\mathfrak I:=\{i\ge 1: B(z_i, \rho/8) \cap D \neq \emptyset\}$. For each $i \in \mathfrak I$, fix any point $x_i\in  B(z_i, \rho/8) \cap D$. Then  $D\subset  \cup_{i\in \mathfrak I} B(z_i,\rho/8) \subset \cup_{i\in \mathfrak I} B(x_i,\rho/4)$, and thus, $D=  \cup_{i\in \mathfrak I} B_D(x_i,\rho/4)$. Further, it holds that 
\begin{align}\label{e:finite-intersection}
	\sum_{i\in \mathfrak I} \1_{B_D(x_i,\rho/2)}\le \sum_{i\ge 1} \1_{B(z_i,\rho)}\le  N \quad \text{on $D$}.
\end{align}
 For all $u \in L^2(F)\cap \sF$, using 
the fact that $m_d(\Sigma)=0$ in the first inequality below,
  Lemma \ref{l:tail-estimate} and \eqref{e:without-A0} in the third inequality, \eqref{e:Phi-scaling-monotone} in the fourth, 
 Corollary \ref{c:Hardy} (with $q=\beta_1$) in the fifth, and \eqref{e:finite-intersection} and 
 the fact that $m_d(\Sigma)=0$ in the sixth,  
 we obtain
\begin{align*}
	&	\sE(u,u)-	\sE^{(\rho)}(u,u) \le  4\int_{D} u(x)^2  \int_{D\setminus B(x, \rho)}J(x,y) dy dx \\
	&\le  4\sum_{i\in I} \int_{B_D(x_i,\rho/4)} u(x)^2  \int_{D\setminus B(x, \rho/4)}J(x,y) dy dx\\
	&\le  \frac{c_1}{\rho^\alpha} \sum_i \int_{B_D(x_i,\rho/4)} u(x)^2 \Phi \left( \frac{\rho/4}{\delta_\Sigma(x)}\right) dx\\
	&\le  \frac{c_2\Phi(1)}{\rho^\alpha} \sum_i \left[  \int_{B_D(x_i,\rho/4)} u(x)^2   dx +  (\rho/4)^{\beta_1 }\int_{B_D(x_i,\rho/4)} \frac{u(x)^2}{\delta_\Sigma(x)^{\beta_1}}dx \right] \\
	&\le  \frac{c_3}{\rho^\alpha} \sum_i \left[ \rho^{\alpha} \int_{B_D(x_i,\rho/2)\times B_D(x_i,\rho/2)} (u(x)-u(y))^2 J(x,y)dxdy +  \int_{B_D(x_i,\rho/2)} u(x)^2   dx  \right] \\
	&\le Nc_3\sE^{(\rho)}(u,u)  + N c_3 \rho^{-\alpha} \lVert u \rVert_{L^2(F)}^2.
\end{align*}
The proof is complete. \qed

For $\rho\in (0, \overline R)$, let $\sF^{(\rho)}$ be the $\sE^{(\rho)}_1$-closure of ${\rm Lip}_c(F)$ in $L^2(F)$. 
Then $\sF^{(\rho)}=\sF$ by Proposition \ref{p:truncation-equivalent}, implying that $(\sE^{(\rho)},\sF)$ is  a regular Dirichlet form on $L^2(F)$. 
 Let $(P^{(\rho)}_t)_{t\ge 0}$ and  
 $X^{(\rho)}=(X^{(\rho)}_t, t\ge 0; \P^x, x\in F\setminus \sN_\rho)$   be the semigroup and the  symmetric Hunt process    associated with  $(\sE^{(\rho)},\sF)$, respectively, where $\sN_\rho$ is a properly exceptional set for $X^{(\rho)}$.

  \subsection{Parabolic H\"older regularity and  near-diagonal  heat kernel bounds}\label{ss:reg}

Consider the bilinear form $(\sE^{(\alpha)}, \sF^{(\alpha)})$ on $L^2(F;dx)$  defined by
 \begin{align*}
 	\sE^{(\alpha)}(f,g)&:= \int_{F \times F} \frac{(f(x)-f(y))(g(x)-g(y))}{ |x-y|^{d+\alpha}} dxdy, \quad		\sF^{(\alpha)}:=\big\{f \in L^2(F): \sE^{(\alpha)}(f,f)<\infty\big\}.
 \end{align*}
Using  
\eqref{e:d-set}
 together with \cite[Theorem 1.1]{CK03} and the stability results in \cite[Theorem 1.18 
and Corollary 1.3]{CKW-jems} (see also \cite[Remark 1.19]{CKW-jems}), the Poincar\'e inequality  holds for $(\sE^{(\alpha)}, \sF^{(\alpha)})$. Thus, by \eqref{e:J-standard-lower}, we get the following   Poincar\'e inequality  for $(\sE,\sF)$.
 
 \begin{lem}\label{l:PI}There exist  $C>0$ and $\delta \in (0,1]$ such that for all $x\in {F}$, $r\in (0, \delta  R_0)$ and $f \in L^1(B_F(x,r))\cap \sF$,
 	\begin{align*}
 		\int_{B_F(x,r)} (f(y) - \overline f_{B_F(x,r)} )^2  dy\le Cr^\alpha	\int_{B_F(x, r/\delta) \times B_F(x, r/\delta)} (f(z)-f(y))^2 J(z,y)dzdy,
 	\end{align*}
 	where $\overline f_{B_F(x,r)} := m_d(B_F(x,r))^{-1} \int_{B_F(x,r)} f(z)dz.$
 \end{lem}

 Throughout this paper, an open set in $F$ is understood with respect to the topology on $F$.  
 For an open subset $U$ of $F$, let $\sF^U$ be the $\sE_1$-closure of $C_c(U) \cap \sF$.
 Then $(\sE,\sF^U)$ is a regular Dirichlet form on $L^2(U)$.  For every  $\rho \in (0,\overline R)$, since $\sE_1$ 
 and $\sE_1^{(\rho)}$ induce  equivalent norms on $C_c(U)$ by Proposition \ref{p:truncation-equivalent},  $(\sE^{(\rho)},\sF^U)$ is also a regular Dirichlet form on $L^2(U)$. 
  Denote by $(P^U_t)_{t\ge 0}$  and $(P^{(\rho),U}_t)_{t\ge 0}$  the semigroups associated 
  with $(\sE,\sF^U)$ and $(\sE^{(\rho)},\sF^U)$, respectively.

Let $Z:=(T_s,X_s)_{s \ge 0}$ be the time-space process corresponding to $X$ where $T_s:=T_0-s$. The law of the time-space process $s \mapsto Z_s$ starting from $(t,x)$ will be denoted by $\P^{(t,x)}$.  For an open subset $U$ of $[0,\infty) \times F$, define $\wh \tau_U:=\inf\{t>0:Z_t \notin U\}$.

\begin{defn}
	{\rm  A 
		 Borel 
		function $q$ on $[0,\infty) \times F$ is said to be \textit{caloric} in $(a,b] \times B_{F}(x_0,r)$ (with respect to $X$), if there is a properly exceptional set 
		$\sN'\supset \sN$  
		such that for any relatively compact open subset $U \subset (a,b] \times B_{F}(x_0,r)$, it holds that 
		$\E^{(t,x)} |q(Z_{\wh \tau_U})|<\infty$  and 
		 $q(t,x)=\E^{(t,x)} q(Z_{\wh \tau_U})$ for all 
		 $(t,x) \in U \cap ([0,\infty) \times ({F} \setminus \sN'))$.	
	}
\end{defn}

For $\rho \in (0, \overline R)$, caloric functions with respect to $X^{(\rho)}$ are defined analogously.

The following notion of an admissible weight function was introduced in  \cite[Definition 2.3]{CK25+}.
\begin{defn}\label{d:weight}
	\rm A function $\Theta: 
	 F \times (0,\infty)\to [1,\infty]$ is called an \textit{admissible weight function} for $\sE$   
with scale function $r^\alpha$ and  localization radius $r_0$,
	if   the following  hold:
	
	\medskip

\setlength{\leftskip}{0.18in}

\noindent (1)  For any $x \in 
	 F$, the map $r \mapsto \Theta(x,r)$ is  non-decreasing and  continuous.

\noindent	 	(2) There exist $\gamma_0 \in (0,\alpha)$  and $C>1$ such that 
	\begin{align*}
		\frac{\Theta(x,R)}{\Theta(x,r)}     \le C  \bigg( \frac{R}{r} \bigg)^{\gamma_0} \quad \text{for all $x\in 
	 F$ and $0<r\le R<\infty$.} 
	\end{align*}
	
	\noindent 	(3) There exist $\eps>0$ and $C>1$ such that for all $x\in
	 F$ and $r\in (0,r_0)$,
	\begin{align}\label{e:blow-up-integral-general}
		\int_{B_F(x,r)} \Theta(y,r)^{1+\eps}dy\le C r^d.
	\end{align}
	
\noindent	 (4) There exist  $\eps>0$ and $C>1$ such that for all $x\in
	 F$, $r\in (0,r_0)$ and any bounded Borel function $u$ with compact support in $F$,
	\begin{align}\label{e:blow-up-Hardy}
		\int_{B_F(x,r)} u(y)^2 \Theta(y,r)^{1+\eps}dy&\le C\big( r^\alpha\sE(u,u)+ \lVert u \rVert_{L^2(F)}^2 \big).
	\end{align}

	\setlength{\leftskip}{0in}
	
\end{defn}

\smallskip

Define 
$\Theta_\Phi(x,r):=\Phi(r/\delta_\Sigma(x))$. 
 Then $\Theta_\Phi$ is an admissible weight function for $\sE$ 
with scale function $r^\alpha$ and localization radius $ R_0$. 
Indeed, property (1) is clear, and property  (2)   follows from \eqref{e:Phi-scaling}. For properties (3) and (4), note that $\Theta_\Phi(x,r)\le 1+ C(r/\delta_\Sigma(x))^{\beta_1}$ for all 
$x\in F$
and  $r>0$  by \eqref{e:Phi-scaling-monotone}. Thus, since $\beta_1\in (0,\gamma\wedge \alpha)$, by Lemma \ref{l:Aikawa-integrals}(i) and  Corollary \ref{c:Hardy},  we deduce that $\Theta_\Phi$ satisfies  \eqref{e:blow-up-integral-general} and \eqref{e:blow-up-Hardy}  with $\eps = ((\gamma \wedge \alpha) - \beta_1)/(2\beta_1)$. Using Proposition \ref{p:truncation-equivalent}, we also see that $\Theta_\Phi$ is an admissible weight function for $\sE^{(\rho)}$ for every $\rho \in (0,\overline R)$, 
with scale function $r^\alpha$ and localization radius $(\rho \wedge R_0)/4$.

 Moreover, by Lemmas \ref{l:tail-estimate} and \ref{l:PI}, conditions TJ$^\Phi_\le(r^\alpha)$ and PI$(r^\alpha)$ in \cite{CK25+} hold for $(\sE,\sF)$ with localization radius $R_0$,  and for $(\sE^{(\rho)},\sF)$, uniformly in $\rho \in (0,\overline R)$, with localization radius $(\rho \wedge R_0)/4$.  Consequently, by \cite[Theorem 12.1]{CK25+},  we obtain the following two propositions, that establish parabolic H\"older regularity  and   near-diagonal lower bounds for the Dirichlet heat kernels.

\begin{prop}\label{p:PHR}
(i)	There exist  $C,\theta>0$ and $\delta \in (0,1)$ such that for all $x_0 \in {F}$, $R \in (0, R_0)$,  $t_0 \ge 0$ and any bounded Borel function $q$ on $[t_0,t_0+R^\alpha] \times {F}$ which is caloric with respect to $X$ in $[t_0,t_0+R^\alpha] \times B_{F}(x_0,R)$, 
	\begin{align}\label{e:PHR}
		|q(s,x)-q(t,y)| \le C \bigg(\frac{|s-t|+|x-y|^\alpha}{
 R^\alpha		} \bigg)^\theta \esssup_{[t_0, t_0+R^\alpha] \times {F}} |q|, 
	\end{align}
	for all $s,t \in [t_0+(1-\delta) R^\alpha, t_0 +  R^\alpha]$ and 
	a.e. $x,y \in B_{F}(x_0, \delta R)$.
	
\noindent (ii) There exist  $C,\theta>0$ and $\delta \in (0,1)$ such that for all $\rho \in (0,\overline R)$, $x_0 \in {F}$, $R \in (0,(\rho \wedge R_0)/4]$,  $t_0 \ge 0$ and any bounded Borel function $q$ on $[t_0,t_0+R^\alpha] \times {F}$ which is caloric with respect to $X^{(\rho)}$ in $[t_0,t_0+R^\alpha] \times B_{F}(x_0,R)$,  
\eqref{e:PHR} holds for all $s,t \in [t_0+(1-\delta) R^\alpha, t_0 +  R^\alpha]$ and a.e. $x,y \in B_F(x_0, \delta R)$.
\end{prop}

\begin{prop}\label{p:NDL}
(i)	There exists  $\eta_0 \in (0,1)$ such that for all $x_0 \in {F}$ and $r\in (0, \eta_0R_0)$, the semigroup $(P^{B_{F}(x_0,r)}_{t})_{t\ge 0}$ has a jointly continuous heat kernel $p^{B_{F}(x_0,r)}(t,x,y)$. Moreover, there exist  $\delta_1,\delta_2\in (0,1)$ and $C>0$ independent of $x_0$ and $r$ such that
	\begin{align*}
		p^{B_{F}(x_0,r)}(t,x,y) \ge C t^{-d/\alpha} \quad \text{for all $t\le  \delta_1r^\alpha$ and $x,y \in 
		B_{F}(x_0, \delta_2 t^{1/\alpha})$
		}.
	\end{align*}
	
	\noindent (ii) 	There exists  $\eta_0 \in (0,1)$ such that for all $\rho \in (0,\overline R)$,  $x_0 \in {F}$ and $r\in (0, \eta_0 (\rho \wedge R_0)/4]$, the semigroup $(P^{(\rho), B_{F}(x_0,r)}_{t})_{t\ge 0}$ has a jointly continuous heat kernel $p^{(\rho), B_{F}(x_0,r)}(t,x,y)$. Moreover, there exist  $\delta_1,\delta_2\in (0,1)$ and $C>0$ independent of $\rho$, $x_0$ and $r$ such that
	\begin{align*}
		p^{(\rho), B_{F}(x_0,r)}(t,x,y) \ge C t^{-d/\alpha} \quad \text{for all $t\le  \delta_1r^\alpha$ and $x,y \in 
		B_{F}(x_0, \delta_2 t^{1/\alpha})$
		}.
	\end{align*}
\end{prop}

We next establish a Nash-type inequality for $(\sE,\sF)$.
\begin{prop}\label{p:Nash}
	There exists  $C>0$ such that for all $u \in \sF$ with $\lVert u \rVert_{L^1(D)}=1$,
	\begin{align*}
		\lVert u \rVert_{L^2(F)}^{2(1+\alpha/d)}  \le  C \left(  \sE(u,u) + R_0^{-\alpha} \lVert u \rVert_{L^2(F)}^2 \right).
	\end{align*}
\end{prop}
\pf 
Recall that \eqref{e:d-set} holds. Hence, when $R_0<\infty$, the result follows from  the third display on page 41 of  \cite{CK03}. When $R_0=\infty$,  the result follows from Lemma \ref{l:PI} and  \cite[Proposition 7.6 and the implication (1) $\Rightarrow$ (2) in Proposition 7.3]{CKW-memo}.
\qed

Combining Proposition \ref{p:Nash} with Proposition \ref{p:truncation-equivalent}, we obtain the following result.
\begin{prop}\label{p:Nash-truncated}
	There exists  $C>0$ such that for all $\rho \in (0,\overline R)$ and  $u \in \sF$ with $\lVert u \rVert_{L^1(D)}=1$,
	\begin{align*}
		\lVert u \rVert_{L^2(F)}^{2(1+\alpha/d)}  \le  C \left(  \sE^{(\rho)}(u,u) + (\rho \wedge R_0)^{-\alpha} \lVert u \rVert_{L^2(F)}^2 \right).
	\end{align*}
\end{prop}

As a standard consequence of Propositions \ref{p:Nash} and \ref{p:Nash-truncated}, we obtain the following ultracontractivity of the semigroups $(P_t)_{t\ge 0}$ and $(P^{(\rho)}_t)_{t\ge 0}$ (see \cite[Theorem 2.1]{CKS87}): There exists $C>1$ such that for all $\rho \in (0,\overline R)$ and $t>0$, 
\begin{align}\label{e:Ultracontractivity}
	\lVert P_t \rVert_{L^1(D) \to L^\infty(D)} \le C t^{-d/\alpha} e^{t/ R_0^\alpha}, \quad 
	\lVert P^{(\rho)}_t \rVert_{L^1(D) \to L^\infty(D)} \le C t^{-d/\alpha} e^{t/(\rho \wedge R_0)^\alpha}.
\end{align}

Combining  \eqref{e:Ultracontractivity}  with parabolic H\"older regularity (Proposition \ref{p:PHR}), we obtain the  next result. See \cite[Subsection 5.5]{GHH18} for details.
\begin{prop}\label{p:NDU}
	For any non-empty open set $U\subset F$ and    $\rho\in (0,\overline R)$,	the semigroups $(P^U_t)_{t\ge 0}$ and $(P_t^{(\rho),U})_{t\ge 0}$ admit
	 jointly continuous heat kernels $p^U(t,x,y)$ and $p^{(\rho),U}(t,x,y)$ on 
	$(0,\infty) \times U \times U$. 
	 Moreover, there exists  $C>0$ independent of $U$ and $\rho$ such that for all 
	 $t>0$ and $x,y\in U$,
	\begin{align*}
		p^U(t,x,y) \le Ct^{-d/\alpha}e^{t/ R_0^\alpha} \quad 
		\text{and} \quad p^{(\rho),U}(t,x,y) \le Ct^{-d/\alpha} e^{t/(\rho \wedge R_0)^\alpha}.
	\end{align*}
\end{prop}

The heat kernels $p^{F}(t,x,y)$ and $p^{(\rho),F}(t,x,y)$ will be denoted by $p(t,x,y)$ and $p^{(\rho)}(t,x,y)$, respectively.

 Since $X$ and $X^{(\rho)}$, $\rho \in (0,\overline R)$, admit jointly continuous heat kernels on 
  $(0,\infty) \times F \times F$, 
 they can be refined to start from every point in  $F$, and the properly exceptional sets $\sN$ and $\sN_\rho$ can be taken to be empty.

In the remainder of this paper, we  take $\sN=\sN_\rho=\emptyset$.
 Further, we use the notations $X^{(\overline R)}_t:=X_t$, $P^{(\overline R)}_t:=P_t$, $p^{(\overline R)}(t,\cdot,\cdot):=p(t,\cdot,\cdot)$,  $P^{(\overline R),U}_t:=P^U_t$ and $p^{(\overline R),U}(t,\cdot,\cdot):=p^U(t,\cdot,\cdot)$ for $t>0$ and open sets $U\subset F$.

By the semigroup property, we can obtain the following refinement of Proposition \ref{p:NDL}.

\begin{lem}\label{l:NDL-modified}
	There exist  $\delta_3\in (0,1/4)$, $\mu>0$ and  $C>0$ such that 
	for  any $\rho \in (0,\overline R]$, $x_0 \in {F}$ and $r>0$, we have
	\begin{align}\label{e:NDL-modified-1}
		p^{(\rho),B_{{F}}(x_0,r)}(t,x,y) \ge C (t^{1/\alpha} \wedge \rho \wedge r\wedge R_0)^{-d}e^{-\mu t/(\rho \wedge r \wedge R_0)^\alpha}
	\end{align}
	for all $t>0$,  $x\in B_{{F}}(x_0, r/4)$  and $y \in B_{F}(x,\delta_3(t^{1/\alpha} \wedge \rho \wedge r\wedge R_0))$.
	Consequently,	it holds that
	\begin{align}   \label{e:NDL-modified-2}	\inf_{x\in B_{{F}}(x_0, r/4)}  \int_{B_{{F}}(x_0,r)} p^{(\rho),B_{{F}}(x_0,r)}(t,x,y) dy \ge C'e^{-\mu t/(\rho \wedge r \wedge R_0)^\alpha} \quad \text{for all $t>0$.}
	\end{align} 
\end{lem}
\pf   Set $r_0:=\rho \wedge r \wedge R_0$ and $\delta_3:= 2^{-2-1/\alpha}\eta_0\delta_1^{1/\alpha}\delta_2   $, where    $\eta_0, \delta_1,\delta_2\in (0,1)$ are the constants in Proposition \ref{p:NDL}. Let $x\in B_{{F}}(x_0, r/4)$, $t>0$ and $y \in B_{F}(x, \delta_3(t^{1/\alpha} \wedge r_0))$.
Note that $B_{F}(x, \eta_0r_0/4)\subset B_{F}(x_0,r)$ and $y \in B_{F}(x,\delta_2t^{1/\alpha})$.
If  $t\le \delta_1(\eta_0r_0/4)^\alpha$, then by Proposition \ref{p:NDL}(ii),
\begin{align*}
	& p^{(\rho),B_{{F}}(x_0,r)}(t,x,y)  \ge p^{(\rho),B_{{F}}(x,\eta_0r_0/4)}(t,x,y)    \ge c_1 t^{-d/\alpha} \ge c_2 (t^{1/\alpha} \wedge r_0)^{-d}.
\end{align*} 
Suppose  $t> \delta_1(\eta_0r_0/4)^\alpha$ and let $N\ge2$ be such that $ (N-1)\delta_1(\eta_0r_0/4)^\alpha < t \le N \delta_1(\eta_0r_0/4)^\alpha$.  Note that $y \in B_{F}(x, 2^{-2-1/\alpha}\eta_0\delta_1^{1/\alpha}\delta_2 r_0) \subset B_{F}(x, \delta_2(t/N)^{1/\alpha})$.
Hence, by the semigroup property, Proposition \ref{p:NDL}(ii) and  \eqref{e:d-set}, 
we get that
\begin{align*}
	&p^{(\rho),B_{{F}}(x_0,r)}(t,x,y)   \ge  p^{(\rho),B_{{F}}(x,\eta_0r_0/4)}(t,x,y) \\
	&\ge 	\int_{B_{{F}}(x, \delta_2 (t/N)^{1/\alpha})} 	\cdots	\int_{B_{{F}}(x, \delta_2 (t/N)^{1/\alpha})}  p^{(\rho),B_{{F}}(x,\eta_0r_0/4)} (t/N,x,z_1) \\
	&\hspace{2.8in} \cdots  p^{(\rho),B_{{F}}(x,\eta_0r_0/4)}(t/N,z_{N-1},y)  dz_1 \cdots dz_{N-1} \\
	&\ge 	(c_3(t/N)^{-d/\alpha})^{N} m_d(B_{{F}}(x, \delta_2 (t/N)^{1/\alpha}))^{N-1}\ge  c_4(t/N)^{-d/\alpha} e^{-c_5(N-1)} \ge  c_6 r_0^{-d} e^{-c_5 t /(\delta_1(\eta_0r_0/4)^\alpha)}.
\end{align*} 
This completes the proof for \eqref{e:NDL-modified-1}. For \eqref{e:NDL-modified-2}, using  \eqref{e:NDL-modified-1} and \eqref{e:d-set}, 
we get that for all $t>0$,
\begin{align*}
	&\inf_{x\in B_{{F}}(x_0, r/4)}  \int_{B_{{F}}(x_0,r)} p^{(\rho),B_{{F}}(x_0,r)}(t,x,y) dy \ge \inf_{x\in B_{{F}}(x_0, r/4)}  \int_{B_{F}(x, \delta_3(t^{1/\alpha} \wedge r_0))} p^{(\rho),B_{{F}}(x_0,r)}(t,x,y) dy\\
	&\ge c_7 (t^{1/\alpha}\wedge r_0)^{-d}e^{-\mu t/r_0^\alpha} \inf_{x\in F} m_d(B_{F}(x, \delta_3(t^{1/\alpha} \wedge r_0)) ) \ge c_8e^{-\mu t/r_0^\alpha}.
\end{align*}
\qed

We end this section with the  observation that when $\Phi(\infty)=\infty$,  the parabolic Harnack inequality  with  space--time scaling  $r^\alpha$  \textit{fails} for $X$; see  \cite[Definition 1.1]{CKW-jems} for the definition of the parabolic Harnack inequality for Hunt processes. 
 Indeed, it is  known that the parabolic Harnack inequality is equivalent to the combination of the near-diagonal lower heat kernel bound in Proposition \ref{p:NDL} and the following condition, referred to as \textbf{(UJS)} in the literature (see \cite[Theorem 1.18]{CKW-jems}):  There  exists $C>0$ such that  
 \begin{align}\label{e:UJS}
 	J(x,y)\le \frac{C}{r^d}\int_{	B_{F}(x,r)} J(z,y)dz 
 	\quad  
 	\text{ for all } x, y\in F  \,\, \text{whenever} \;\, 0< r \le \frac12 ( |x-y| \wedge \wh R).
 \end{align}
Since $|z-y| \asymp |x-y|$ for $z\in B(x,r)$ and  $0< r \le  |x-y|/2$, 
by the monotonicity of $\Phi$, \eqref{e:Phi-scaling}, \eqref{e:JJ} and Lemma \ref{l:Aikawa-integrals}, we have that, for $0< r \le \frac12 ( |x-y| \wedge \wh R),$
\begin{align*}	& \frac{1}{r^d}\int_{	B_{F}(x,r)} J(z,y)dz \asymp \frac{1}{r^d|x-y|^{d+\alpha}} \int_{	B_{F}(x,r)} 
	\Phi\left( \frac{(|x-y|\wedge A_0)^2}{(\delta_\Sigma(z)\wedge A_0) (\delta_\Sigma(y) \wedge A_0)}\right)dz\\
	&\le 
	\frac{c}{|x-y|^{d+\alpha}} 
	\Phi\left( \frac{(|x-y|\wedge A_0)^2}{((\delta_\Sigma(x) \vee r) \wedge A_0) (\delta_\Sigma(y) \wedge A_0)} 
	\right).
\end{align*}
Thus, when $\Phi(\infty)=\infty$,   $J(x,y) = \infty$ for $ x \in \Sigma$ and  $y \in D$, but the right-hand side of  \eqref{e:UJS} is finite, so \eqref{e:UJS} fails. Therefore, the parabolic Harnack inequality does not hold.

\section{Further analysis of truncated Dirichlet forms}\label{s:further-anal}

Recall that the truncated Dirichlet form $\sE^{(\rho)}$ is defined in \eqref{def:truncated-forms}, 
and $\sE^{(\overline R)}=\sE$.
In Subsection \ref{ss:Meyer}, we use Meyer's construction to prove some inequalities between
the heat kernels $p^{(\rho)}(t,x,y)$ and $p^{(\rho')}(t,x,y)$ of the truncated forms
for different parameters $\rho, \rho'\in (0, \overline{R}]$. 
In Subsection \ref{ss:4.2}, we prove some off-diagonal heat kernel estimates for the truncated forms.

The Dirichlet form $(\sE^{(\rho)},\sF)$ on $L^2(F)$ is called \textit{conservative} if $P^{(\rho)}_t\1_F=\1_F$ for all $t>0$.

\begin{prop}\label{p:conservative}
	For any $\rho \in (0,\overline R]$,  the Dirichlet form $(\sE^{(\rho)},\sF)$ is conservative.
\end{prop}
\pf By Lemma \ref{l:NDL-modified},   we obtain for all $\rho \in (0,\overline R]$,  $x_0 \in {F}$, $r\in (0, \rho \wedge R_0)$ and $t\in (0,  r^\alpha]$,	\begin{align*}
	\inf_{x\in B_{{F}}(x_0, r/4)}  \int_{B_{{F}}(x_0,r)} p^{(\rho),B_{{F}}(x_0,r)}(t,x,y) dy \ge c_1.	\end{align*}
Thus, 
 since $(\sE^{(\rho)},\sF)$ is a regular Dirichlet form without killing,  the result follows from \cite[Lemma 4.6]{GHL17}.  \qed

 \subsection{Meyer's decomposition}\label{ss:Meyer}

For $x\in F$, $0<\rho<\rho'\le \overline R$ and 
open set $U\subset D$ with $\overline U\subset D$, define
\begin{align*}
	\sT^{(\rho,\rho')}(x):=\int_{B_F(x,\rho')\setminus B(x, \rho)} J(x,y)dy, \qquad 	 \sT^{(\rho,\rho')}(U):=\sup_{x\in U}\sT^{(\rho,\rho')}(x).
\end{align*}
By Lemma \ref{l:tail-estimate}, we have  for all $0<\rho<\rho'\le \overline R$,
\begin{align} \label{e:finite-jumps-added}	\sT^{(\rho,\rho')}(U)\le \frac{c_1}{\rho^\alpha}  \sup_{x\in U}\Phi\bigg(\frac{ \rho \wedge A_0}{\delta_\Sigma(x) \wedge A_0}\bigg) \le \frac{c_1}{\rho^\alpha} \Phi\bigg(\frac{ \rho \wedge A_0}{\dist(\overline U, \Sigma) \wedge A_0}\bigg)<\infty.
	\end{align}

Now we
establish relations among the heat kernels $p^{(\rho)}(t,x,y)$, $\rho \in (0,\overline R]$, using Meyer's construction. A similar result was obtained in \cite{BGK09}. However, in our  setting, the integral of the additional jumps $\sT^{(\rho,\rho')}(x)$ is  not bounded on $D$ in general. Thus some modification is needed.

 Let  $0<\rho<\rho'\le \overline R$. For $n\ge 1$, define 
 \begin{align*}
 K^{(\rho,\rho')}_n(x,y)&:=\1_{D_{1/n}\times D_{1/n}}(x,y) (J^{(\rho')}(x,y)-J^{(\rho)}(x,y) ), \quad \;\;x,y \in F,\\
 	J^{(\rho,\rho')}_n(x,y)&:=J^{(\rho)}(x,y) + K^{(\rho,\rho')}_n(x,y), \quad\;\; x,y \in F,
 \end{align*}
\begin{align*}
	\sE^{(\rho,\rho',n)}(u,v) 
	:=\frac12 \int_{F\times F} (u(x)-u(y))(v(x)-v(y)) 	J^{(\rho,\rho')}_n(x,y) dxdy, \quad\;\; u,v\in \sF.
\end{align*} 
By Proposition \ref{p:truncation-equivalent}, there exists $C>1$ such that $\sE_1^{(\rho)}(u,u) \le \sE_1^{(\rho,\rho',n)}(u,u) \le \sE_1(u,u) \le C\sE_1^{(\rho)}(u,u) $ for all $n\ge 1$ and $u\in \sF$. Thus $(\sE^{(\rho,\rho',n)},\sF)$ is a regular Dirichlet form on $L^2(F)$.

\begin{lem}\label{l:Mosco-convergence}
	$( \sE^{(\rho,\rho',n)},\sF)$  converges to $( \sE^{(\rho')},\sF)$ in the sense of Mosco, i.e., the following  hold.

	\noindent (i) For any sequence $(u_n)_{n\ge 1}$ in $L^2(F)$ that converges weakly to $u$ in $L^2(F)$,
	\begin{align*}
		\liminf_{n\to \infty} \sE^{(\rho,\rho',n)}(u_n,u_n) \ge \sE^{(\rho')}(u,u).
	\end{align*}
		\noindent (ii) For any $u\in L^2(F)$, there is a sequence $(u_n)_{n\ge 1}$ in $L^2(F)$ that  converges strongly to $u$ in $L^2(F)$ such that
	\begin{align*}
		\limsup_{n\to \infty} \sE^{(\rho,\rho',n)}(u_n,u_n) \le \sE^{(\rho')}(u,u).
	\end{align*}
\end{lem}
\pf (i) Let $(u_n)_{n\ge 1}$ be a sequence in $L^2(F)$ that  converges weakly to $u$ in $L^2(F)$. 
By taking a subsequence if necessary, 
we may assume that $\sE^{(\rho,\rho',n)}(u_n,u_n)<\infty$ for all $n\ge 1$ and  the limit $A:=\lim_{n\to \infty} \sE^{(\rho,\rho',n)}(u_n,u_n)$ exists and is finite. By the uniform boundedness principle, we have $\sup_{n\ge 1}  \sE^{(\rho,\rho',n)}_1(u_n,u_n)<\infty$. Since  $\sup_{n\ge 1}  \sE^{(\rho,\rho',1)}_1(u_{n},u_{n})\le \sup_{n\ge 1}  \sE^{(\rho,\rho',n)}_1(u_n,u_n)<\infty$, by the Banach-Saks theorem, there exists a subsequence $(u_{a_n})_{n\ge 1}$ such that the Ces\'aro 
mean
$v_n:=\sum_{k=1}^n u_{a_k}/n$ is $ \sE^{(\rho,\rho',1)}_1$-convergent to some $v\in \sF$.  Since $v_n$ converges to $v$ in $L^2(F)$, we have $v=u$ a.e. on $F$, and there exists a subsequence $(v_{b_n})_{n\ge 1}$ converging pointwise to $u$ a.e. on $F$. By Fatou's lemma, we get that for all $n\ge 2$,
\begin{align}\label{e:Mosco-1}
	\liminf_{m\to \infty} 	 \sE^{(\rho,\rho',a_{b_n})}(v_{b_m},v_{b_m}) \ge  \sE^{(\rho,\rho',a_{b_n})}(u,u).
\end{align}
On the other hand, using the Cauchy-Schwarz inequality, we see that for all $n\ge 2$,
\begin{align}\label{e:Mosco-2}
\limsup_{m\to \infty} 	 \sE^{(\rho,\rho',a_{b_n})}(v_{b_m},v_{b_m})& 
\le \limsup_{m\to \infty} \frac{1}{b_m} \bigg(  \sum_{k=1}^{b_n-1}  \sE^{(\rho,\rho',a_{b_n})}(u_{a_k},u_{a_k}) +  \sum_{k=b_n}^{b_m} \sE^{(\rho,\rho',a_{b_n})}(u_{a_k},u_{a_k}) \bigg)  \nn\\	  
&\le \sup_{k\ge b_n}\sE^{(\rho,\rho',a_{b_n})}(u_{a_k},u_{a_k}) .
\end{align}
Combining \eqref{e:Mosco-1} and \eqref{e:Mosco-2}, we arrive at
\begin{align*}
A &=  \lim_{n\to \infty} \sup_{k\ge b_n}\sE^{(\rho,\rho',a_k)}(u_{a_k},u_{a_k})  \ge \lim_{n\to \infty} \sup_{k\ge b_n}\sE^{(\rho,\rho',a_{b_n})}(u_{a_k},u_{a_k})  \\&\ge \lim_{n\to \infty} \limsup_{m\to \infty}  \sE^{(\rho,\rho',a_{b_n})}(v_{b_m},v_{b_m}) \ge \lim_{n\to \infty}  \sE^{(\rho,\rho',a_{b_n})}(u,u) = \sE^{(\rho')}(u,u),
\end{align*}
where we used the monotone convergence theorem in the last equality.

\noindent (ii) The result follows by taking $u_n=u$ for all $n\ge 1$. \qed

Let $X^{(\rho,\rho',n)}$ and   $(P^{(\rho,\rho',n)}_t)_{t\ge 0}$ be the symmetric Hunt process and the  semigroup associated with  $( \sE^{(\rho,\rho',n)},\sF)$. 
 Combining  Lemma \ref{l:Mosco-convergence} with \cite[Theorem 2.4.1 and Corollary 2.6.1]{Mosco}, we get 
\begin{cor}\label{c:Mosco-semigroup}
For all $0<\rho<\rho'\le \overline R$, 
 $t>0$ and $f\in L^2(F)$, $P^{(\rho,\rho',n)}_t f$ converges to $P^{(\rho')}_t f$ in $L^2(F)$ as $n\to \infty$.
\end{cor}

For every $n\ge 1$, since  $ \sE^{(\rho,\rho',n)}(u,u)\ge \sE^{(\rho)}(u,u)  $ for all $u\in \sF$ and $J_n^{(\rho,\rho')}(x,y)\le J(x,y)$, by repeating the argument leading to Proposition  \ref{p:PHR}, we get the next result.
\begin{prop}\label{p:PHR-n}
	There exist $C,\theta>0$ and $\delta \in (0,1)$ such that for all $0<\rho<\rho'\le \overline R$, $n\ge 1$, $x_0 \in {F}$,  $R \in (0,(\rho \wedge R_0)/4]$,   $t_0 \ge 0$ and any bounded Borel function $q$ on $[t_0,t_0+R^\alpha] \times {F}$ which is caloric with respect to $X^{(\rho,\rho',n)}$ in $[t_0,t_0+R^\alpha] \times B_{F}(x_0,R)$, 
	\begin{align*}
		|q(s,x)-q(t,y)| \le C \bigg(\frac{|s-t|+|x-y|^\alpha}{
		R^\alpha} \bigg)^\theta \esssup_{[t_0, t_0+R^\alpha] \times {F}} |q|, 
	\end{align*}
	for all $s,t \in [t_0+(1-\delta) R^\alpha, t_0 +  R^\alpha]$ and a.e. $x,y \in B_{F}(x_0, \delta R)$.
\end{prop}
Moreover, by  Proposition \ref{p:Nash-truncated}, 	there exists $C>0$ such that for all $n\ge 1$ and  $u \in \sF$ with $\lVert u \rVert_{L^1(F)}=1$,
\begin{align*}
	\lVert u \rVert_{L^2(F)}^{2(1+\alpha/d)}  \le  C \left(  \sE^{(\rho,\rho',n)}(u,u) + (\rho \wedge R_0)^{-\alpha} \lVert u \rVert_{L^2(F)}^2 \right).
\end{align*}
Hence, repeating the argument leading to Proposition \ref{p:NDU}, we obtain
\begin{prop}\label{p:NDU-n}
	For all $0<\rho<\rho'\le \overline R$ and $n\ge 1$,	the semigroup $(P_t^{(\rho,\rho',n)})_{t\ge 0}$ has a jointly continuous heat kernel $p^{(\rho,\rho',n)}(t,x,y)$  on $(0,\infty) \times {F} \times {F}$. Moreover, there exists  $C>0$ independent of $\rho,\rho'$ and $n$ such that $	p^{(\rho,\rho',n)}(t,x,y) \le Ct^{-d/\alpha} e^{t/(\rho \wedge R_0)^\alpha}$ for all $t>0$ and $x,y\in {F}$.
\end{prop}
In the following, for all $0<\rho<\rho'\le \overline R$, $n\ge1$,  $t>0$ and $f\in L^2(F)$, we always use the regularized versions of $P^{(\rho,\rho',n)}_tf$ and $P^{(\rho')}_tf$, defined by $P_t^{(\rho,\rho',n)}f(x)=\int_F p^{(\rho,\rho',n)}(t,x,y) f(y)dy$ and  $P_t^{(\rho')}f(x)=\int_F p^{(\rho')}(t,x,y) f(y)dy$ for  $x\in F$.

\begin{lem}\label{l:Mosco-semigroup-pointwise}
For all  $0<\rho<\rho'\le \overline R$,  $t>0$ and $f \in L^2(F) \cap L^\infty(F)$, we have
\begin{align*}
	\limsup_{n\to \infty} P^{(\rho,\rho',n)}_t f(x) =  P^{(\rho')}_t f(x) \quad \text{for all $x\in F$}.
\end{align*}
\end{lem}
\pf Let  $(a_n)_{n\ge 0}$ be an arbitrary increasing sequence of positive integers.  By Corollary \ref{c:Mosco-semigroup}, there exists a subsequence  $(\wt a_n)_{n\ge 1}$ such that $\lim_{n\to \infty} P^{(\rho,\rho',\wt a_n)}_t f(y) =  P^{(\rho')}_t f(y)$ for a.e. $y\in F$. By Propositions \ref{p:PHR} and \ref{p:PHR-n}, since $P^{(\rho,\rho',\wt a_n)}_t$ and $P^{(\rho')}_t$ are $L^\infty$-contractive,  there exist  $\delta,\theta>0$ and $c_1>0$ such that
\begin{align*}
	|P^{(\rho,\rho',\wt a_n)}_t f(y) - P^{(\rho,\rho',\wt a_n)}_t f(z)| +  
	|P^{(\rho')}_t f(y) - P^{(\rho')}_t f(z)| \le c_1|y-z|^\theta  \lVert f \rVert_{L^\infty(F)}
\end{align*}
for all $n\ge 1$, $x\in F$ and  $y,z \in B_{F}(x,\delta)$. Using these, one can easily show that for any $x\in F$, $\lim_{n\to \infty} P^{(\rho,\rho',\wt a_n)}_t f(x) =  P^{(\rho')}_t f(x)$. Since $(a_n)_{n\ge 0}$  is arbitrary,  the proof is complete. \qed

For $n\ge 1$, define
\begin{align*}
 N^{(\rho,\rho')}_n(x):=\int_F 		K^{(\rho,\rho')}_n(x,y) dy, \quad \;\; x\in F.
\end{align*} 
Note that $ N^{(\rho,\rho')}_n(x)=0$ for all $x\notin D_{1/n}$ and 
\begin{align}\label{e:additional-jumps}
N^{(\rho,\rho')}_n(x)&= \int_{D_{1/n} }  (J^{(\rho')}(x,y)-J^{(\rho)}(x,y) )\,dy  \le\int_{F \cap (B(x, \rho')\setminus B(x,\rho))} J(x,y)dy  = \sT^{(\rho,\rho')}(x)
\end{align}
for all $x\in  D_{1/n}$. Hence, by \eqref{e:finite-jumps-added}, it holds that 
\begin{align*}
	\sup_{x\in F} N^{(\rho,\rho')}_n(x) \le \sup_{x\in D_{1/n}} \sT^{(\rho,\rho')}(x) \le \frac{C}{\rho^\alpha} \Phi\bigg( \frac{\rho \wedge A_0}{n^{-1} \wedge A_0}\bigg) <\infty.
\end{align*}
 By Meyer's construction \cite{Meyer} or the Ikeda-Nagasawa-Watanabe piecing together procedure \cite{INW66}, we have that  for  all $t>0$, $x\in F$ and  Borel  $E\subset F$,
\begin{align}\label{e:Meyer-decomposition}
\P^x \big(  X^{(\rho,\rho',n)}_t \in E\big) 	& =  	\E^x \left[ \1_E( X^{(\rho)}_t)  e^{ - H^{(\rho,\rho',n)}_t} \right]\nn\\
&\quad + \E^x \left[ \int_0^t e^{-H^{(\rho,\rho',n)}_s}   \int_{D}  \P^w \big(  X^{(\rho,\rho',n)}_{t-s} \in E\big)  K^{(\rho,\rho')}_n(X^{(\rho)}_s,w) \, dw\, ds\right],
\end{align}
where\begin{align*}	H^{(\rho,\rho',n)}_t:=  \int_0^t  N^{(\rho, \rho')}_n(X^{(\rho)}_u)du.\end{align*}
See \cite[Lemma 3.1(b)]{BGK09} and \cite{BGK09-correction} for details. 

For $\rho \in (0,\overline  R]$ and an open set $U\subset F$, let $\tau^{(\rho)}_U:=\inf\{t>0: X^{(\rho)}_t \notin U\}$ denote the first exit time from $U$ by $X^{(\rho)}$. The subprocess $X^{(\rho), U}$  of $X^{(\rho)}$ killed upon leaving  $U$ is defined as $X^{(\rho), U}_t = X^{(\rho)}_t$ for $t<\tau^{(\rho)}_U$ and $X^{(\rho), U}_t = \partial$ for $t\ge \tau^{(\rho)}_U$, where $\partial$ is a cemetery state.

\begin{lem}\label{l:comparison-Meyer-n}
	Let  $0<\rho<\rho'\le \overline R$ and $n\ge 1$. The following hold.

	\noindent (i) 	For any open subset  $U\subset D_{1/n}$, all  $t>0$, $x\in U$ and Borel $E\subset F$, we have
	\begin{align*}
&\P^x \big(  X^{(\rho,\rho',n)}_t \in E\big)\\
& \ge \int_0^t \int_U e^{-s \sT^{(\rho,\rho')}(U)} p^{(\rho),U}(s, x, z) \int_{D_{1/n} \cap (B(z,\rho')\setminus B(z,\rho))}  J(z,w)   \P^w \big(  X^{(\rho,\rho',n)}_{t-s} \in E\big) \,dw\, dz\, ds.
	\end{align*}

	\noindent (ii)	For all  $t>0$, $x\in F$ and Borel $E\subset F$, we have
	\begin{align*}
		\P^x \big(  X^{(\rho,\rho',n)}_t \in E\big) &\le  \P^x \big(  X^{(\rho)}_t \in E\big)\\
		&\!\!\!\! +\int_0^t \int_{D_{1/n}}  p^{(\rho)}(s, x, z) \int_{D_{1/n} \cap( B(z,\rho')\setminus B(z,\rho))}  J(z,w)   \P^w \big(  X^{(\rho,\rho',n)}_{t-s} \in E\big) \,dw\, dz\, ds.
	\end{align*}

\end{lem}
\pf  (i) Let  $U\subset D_{1/n}$ be an open subset. For all $s<\tau^{(\rho)}_U$, using \eqref{e:additional-jumps}, we see that
\begin{align*}
	H^{(\rho,\rho',n)}_s\le  \int_0^s  \sT^{(\rho, \rho')}(X^{(\rho),U}_u)du \le s \sT^{(\rho,\rho')}(U).
\end{align*}
Hence, using \eqref{e:Meyer-decomposition}, we obtain 
\begin{align*}
&	\P^x \big(  X^{(\rho,\rho',n)}_t \in E\big)\ge \E^x  \int_0^t  \1_{\{s<\tau^{(\rho)}_U\}} e^{-H^{(\rho,\rho',n)}_s}   \int_{D}  \P^w \big(  X^{(\rho,\rho',n)}_{t-s} \in E\big)  K^{(\rho,\rho')}_n(X^{(\rho)}_s,w) \, dw\, ds \nn\\
		&\ge \E^x  \int_0^t  \1_{\{s<\tau^{(\rho)}_U\}} e^{-s\sT^{(\rho,\rho')}(U)}   \int_{D}  \P^w \big(  X^{(\rho,\rho',n)}_{t-s} \in E\big)  K^{(\rho,\rho')}_n(X^{(\rho),U}_s,w) \, dw\, ds\nn\\
		&= \int_0^t \int_U e^{-s \sT^{(\rho,\rho')}(U)} p^{(\rho),U}(s, x, z) \int_{D_{1/n} \cap (B(z,\rho')\setminus B(z,\rho))}  J(z,w)   \P^w \big(  X^{(\rho,\rho',n)}_{t-s} \in E\big) \,dw\, dz\, ds.
\end{align*}

\noindent (ii) The result follows immediately from \eqref{e:Meyer-decomposition}.
\qed

\begin{prop}\label{p:comparison-Meyer-lower-bound}
	Let  $0<\rho<\rho'\le \overline R$ and let  $U\subset D$ be an open subset with $\overline U \subset D$.  For all  $t>0$, $x\in U$ and $y \in F$, we have
	\begin{equation*}
		p^{(\rho')}(t,x,y)\ge \int_0^t \int_U e^{-s \sT^{(\rho,\rho')}(U)} p^{(\rho),U}(s, x, z) \int_{B_F(z,\rho')\setminus B(z,\rho)}  J(z,w)  p^{(\rho')}(t-s,w,y) \,dw\, dz\, ds.
	\end{equation*}
\end{prop}
\pf Since $p^{(\rho')}(t,x,y)$   is jointly continuous, using  Lemma 
\ref{l:Mosco-semigroup-pointwise},    Lemma \ref{l:comparison-Meyer-n}(i) and Fatou's lemma, we obtain 
	\begin{align*}
	&p^{(\rho')}(t,x,y) = \lim_{\delta\to0} \frac{\P^x \big(  X^{(\rho')}_t \in B_{F}(y,\delta)\big)}{m_d(B_{F}(y,\delta))} = \lim_{\delta\to0}\lim_{n\to \infty} \frac{ \P^x \big(  X^{(\rho,\rho',n)}_t \in B_{F}(y,\delta)\big)}{m_d(B_{F}(y,\delta))} \\
	& \ge \liminf_{\delta\to 0} \liminf_{n\to \infty} \int_0^t \int_U e^{-s \sT^{(\rho,\rho')}(U)} p^{(\rho),U}(s, x, z) \\
	&\qquad \qquad \qquad \qquad \times \int_{D_{1/n} \cap (B(z,\rho')\setminus B(z,\rho))}  J(z,w) \frac{  \P^w \big(  X^{(\rho,\rho',n)}_{t-s} \in B_{F}(y,\delta)\big)}{m_d(B_{F}(y,\delta)) } \,dw\, dz\, ds\\
	& \ge \int_0^t \int_U e^{-s \sT^{(\rho,\rho')}(U)} p^{(\rho),U}(s, x, z) \int_{ B_F(z,\rho')\setminus B(z,\rho)}  J(z,w)  p^{(\rho')}(t-s,w,y)\,dw\, dz\, ds.
\end{align*}
\qed

Recall that the constant $\Lambda$ is defined in \eqref{e:def-Lambda}.

\begin{prop}\label{p:blow-up-integrated}
		There exists $C>0$ such that for all  $\rho \in (0,\overline R]$, $x\in F$ and $t\in (0, ((\rho \wedge R_0)/\Lambda)^\alpha)$,
	\begin{align*}
		\int_0^t \int_F p^{(\rho)}(s, x, z) \Phi \bigg( \frac{t^{1/\alpha} \wedge A_0}{\delta_\Sigma(z)\wedge A_0}\bigg)   dz\, ds \le Ct.
	\end{align*}
\end{prop}
\pf 
Let  $\rho \in (0,\overline R]$, $x\in  F$ and $t\in (0, ((\rho \wedge R_0)/\Lambda)^\alpha)$. Set $\sigma:=t^{1/\alpha}$. Recall that $(\sE^{(\rho)},\sF)$ is conservative. Using  symmetry, Proposition \ref{p:comparison-Meyer-lower-bound} and  Fubini's theorem, we obtain
\begin{align}\label{e:blow-up-integrated-1}
&1=\int_F p^{(\rho)}(t,x,y)dy=\int_F	p^{(\rho)}(t,y,x)dy\nn\\
&\ge \int_{D_{\sigma/2}} \int_0^t \int_{D_{\sigma/2}} e^{-s \sT^{(\sigma,\rho)}(D_{\sigma/2})} p^{(\sigma),D_{\sigma/2}}(s, y, w) \int_{B_F(w,\rho)\setminus B(w,\sigma)} \!\!\!\!\!\!\!\!\!\!\!\!J(w,z)  p^{(\rho)}(t-s,x,z) \,dz\, dw\, ds\,dy\nn\\
&= \int_0^t  \int_F p^{(\rho)}(t-s,x,z) e^{-s \sT^{(\sigma,\rho)}(D_{\sigma/2})}\nn\\
&\qquad\qquad \times \int_{(B_F( z, \rho) \cap D_{\sigma/2}) \setminus B_F(z,\sigma)} \int_{D_{\sigma/2}} p^{(\sigma),D_{\sigma/2}}(s, w, y) dy\, J(z,w)  \,dw\, dz\, ds.
\end{align}
 Using \eqref{e:finite-jumps-added} with  $\dist(\overline{D_{\sigma/2}}, 
 \Sigma) \ge \sigma/2$, \eqref{e:Phi-scaling} and \eqref{e:without-A0},  we get that for all $s\in (0,t)$,
\begin{align}\label{e:blow-up-integrated-2}
	s \sT^{(\sigma,\rho)}(D_{\sigma/2}) \le  \frac{c_1s}{\sigma^\alpha} \Phi\bigg( \frac{\sigma \wedge A_0}{(\sigma/2)\wedge A_0}\bigg) \le \frac{c_2s}{\sigma^\alpha} \le c_2.
\end{align}
 Furthermore, applying Lemma \ref{l:NDL-modified}, we obtain for all $s\in (0,t)$ and $w\in  D_{\sigma}$, \begin{align}\label{e:blow-up-integrated-3}	\int_{D_{\sigma/2}} p^{(\sigma),D_{\sigma/2}}(s,w,y) dy \ge 	\int_{B(w,\sigma/2)} p^{(\sigma),B(w,\sigma/2)}(s,w,y) dy \ge c_3 e^{-c_4s/ \sigma^\alpha} \ge c_3 e^{-c_4}.\end{align} 
Since $\Lambda\sigma<\rho \wedge R_0$, using \eqref{e:blow-up-integrated-3} and Lemma  \ref{l:lower-tail-estimate},  we obtain for all $s\in (0,t)$,
\begin{align}\label{e:blow-up-integrated-4}
&\int_{(B_D( z, \rho) \cap D_{\sigma/2}) \setminus B(z,\sigma)} \int_{D_{\sigma/2}} p^{(\sigma),D_{\sigma/2}}(s, w, y) dy\, J(z,w)  \,dw\nn\\
&\ge \int_{(B_D( z, \Lambda\sigma) \cap D_{\sigma}) \setminus B(z,\sigma)} \int_{D_{\sigma/2}} p^{(\sigma),D_{\sigma/2}}(s, w, y) dy\, J(z,w)  \,dw\nn\\
& \ge c_5\int_{(B_D( z, \Lambda\sigma) \cap D_{\sigma}) \setminus B(z,\sigma)} J(z,w)  \,dw   \ge \frac{c_6}{\sigma^\alpha} \Phi\bigg(\frac{\sigma\wedge A_0}{\delta_\Sigma(z) \wedge A_0}\bigg).
\end{align}
Combining \eqref{e:blow-up-integrated-1}, \eqref{e:blow-up-integrated-2} and \eqref{e:blow-up-integrated-4} with 
the conservativeness, we obtain the desired result. \qed

\begin{lem}\label{l:blow-up-integrated-finite}

	 (i) For all   $\rho' \in (0,\overline R]$, $\rho \in (0,\overline R)$ and $t>0$, 
	\begin{align*}
		\sup_{x\in F}		\int_0^t \int_{F}  p^{(\rho')}(s, x, z) \int_{F\setminus B(z,\rho)}  J(z,w)   \,dw\, dz\, ds <\infty.
	\end{align*}
	
	\noindent (ii) For all   $\rho' \in (0,\overline R]$ and $\rho \in (0,\overline R)$, 
	 \begin{align*}
	 	\lim_{t\to 0} \sup_{x\in F}		\int_0^t \int_{F}  p^{(\rho')}(s, x, z) \int_{F\setminus B(z,\rho)}  J(z,w)   \,dw\, dz\, ds  = 0.
	 \end{align*}
\end{lem}
\pf Let $\rho' \in (0,\overline R]$ and $\rho \in (0,\overline R)$. By Lemma \ref{l:tail-estimate},  \eqref{e:Phi-scaling-monotone} and \eqref{e:without-A0},  for all $t>0$ and  $x\in F$,
	\begin{align}\label{e:blow-up-integrated-finite-1}
			&\int_0^t \int_{F}  p^{(\rho')}(s, x, z) \int_{F\setminus B(z,\rho)}  J(z,w)   \,dw\, dz\, ds \le  \frac{c_1}{\rho^\alpha} \int_0^t \int_F p^{(\rho')}(s, x, z) \Phi \bigg( \frac{\rho \wedge A_0}{\delta_\Sigma(z)\wedge A_0}\bigg)   dz\, ds \nn\\
			& \le 
			\frac{c_2}{\rho^\alpha} \bigg(1 +  \frac{\rho^{\beta_1}}{t^{\beta_1/\alpha}}\bigg)\int_0^t \int_F p^{(\rho')}(s, x, z) \Phi \bigg( \frac{t^{1/\alpha} \wedge A_0}{\delta_\Sigma(z)\wedge A_0}\bigg)   dz\, ds .
\end{align}

\noindent (i) By \eqref{e:blow-up-integrated-finite-1}, it suffices to show that for all $t>0$,
\begin{align}\label{e:blow-up-integrated-finite}
\sup_{x\in F}	\int_0^t \int_F p^{(\rho')}(s, x, z) \Phi \bigg( \frac{t^{1/\alpha} \wedge A_0}{\delta_\Sigma(z)\wedge A_0}\bigg)   dz\, ds  <\infty.
\end{align}
Set $T_0:=((\rho' \wedge R_0)/\Lambda)^\alpha$. By Proposition \ref{p:blow-up-integrated}, \eqref{e:blow-up-integrated-finite} holds for all  $t\in (0,T_0)$. Suppose \eqref{e:blow-up-integrated-finite} holds for all $t\in (0,nT_0)$ with $n\ge 1$. For all $t\in [nT_0,(n+1)T_0)$ and $x\in F$, using the semigroup property, Fubini's theorem  and the induction hypothesis, we obtain
	\begin{align*}
&	\int_{nT_0}^t \int_F p^{(\rho')}(s, x, z) \Phi \bigg( \frac{t^{1/\alpha} \wedge A_0}{\delta_\Sigma(z)\wedge A_0}\bigg)   dz\, ds\\ &= 	  \int_F p^{(\rho')}(nT_0, x, w)  \int_{nT_0}^t \int_F p^{(\rho')}(s -nT_0,w,z) \Phi \bigg( \frac{t^{1/\alpha}  \wedge A_0}{\delta_\Sigma(z)\wedge A_0}\bigg)   dz\, ds\, dw\\
&\le \left(\sup_{w\in F}  \int_0^{t-nT_0} \int_F p^{(\rho')}(s,w,z) \Phi \bigg( \frac{t^{1/\alpha}  \wedge A_0}{\delta_\Sigma(z)\wedge A_0}\bigg)   dz\, ds\right) \int_F p^{(\rho')}(nT_0, x,w )dw\\
&\le \sup_{w\in F}  \int_0^{t-nT_0} \int_F p^{(\rho')}(s,w,z) \Phi \bigg( \frac{t^{1/\alpha}  \wedge A_0}{\delta_\Sigma(z)\wedge A_0}\bigg)   dz\, ds<\infty.
\end{align*}
Hence, \eqref{e:blow-up-integrated-finite} holds for all $t\in (0,(n+1)T_0)$. The result now follows by induction.

\noindent (ii) Combining \eqref{e:blow-up-integrated-finite-1} with Proposition \ref{p:blow-up-integrated},   since 
$\beta_1<\alpha$ by \eqref{e:beta1}, 
we obtain
\begin{align*}
	\lim_{t\to0} \sup_{x\in F}		\int_0^t \int_{F}  p^{(\rho')}(s, x, z) \int_{F\setminus B(z,\rho)}  J(z,w)   \,dw\, dz\, ds \le  \lim_{t\to 0} 
	\frac{c_3t}{\rho^\alpha}\bigg(1 +  \frac{\rho^{\beta_1}}{t^{\beta_1/\alpha}}\bigg) =0.
\end{align*}
 \qed

\begin{prop}\label{p:comparison-Meyer}
	For all  $0<\rho<\rho'\le \overline R$, $t>0$ and $x,y\in F$, we have
	\begin{equation*}
		p^{(\rho')}(t,x,y)\le p^{(\rho)}(t,x,y) +  \int_0^t \int_F p^{(\rho)}(s, x, z) \int_{B_F(z,\rho')\setminus B(z,\rho)} J(z,w)  p^{(\rho')}(t-s,w,y) \,dw\, dz\, ds.
	\end{equation*}
\end{prop}
\pf 
Let  $0<\rho<\rho'\le \overline R$, $t>0$, $x\in F$, and let $E\subset F$ be a Borel subset with $m_d(E)<\infty$.  Note that $\1_E \in L^2(F) \cap L^\infty(F)$.
By Lemma \ref{l:Mosco-semigroup-pointwise} and  Lemma \ref{l:comparison-Meyer-n}(ii), we have
\begin{align*}	&\P^x \big(  X^{(\rho')}_t \in E\big)= \lim_{n\to \infty}	\P^x \big(  X^{(\rho,\rho',
		n)}_t \in E\big) \\	&\le  \P^x \big(  X^{(\rho)}_t \in E\big) +\lim_{n\to \infty} \int_0^t \int_{F}  p^{(\rho)}(s, x, z) \int_{ B_F(z,\rho')\setminus B(z,\rho)}  J(z,w)   \P^w \big(  X^{(\rho,\rho',
		n)}_{t-s} \in E\big) \,dw\, dz\, ds\\	&= \P^x \big(  X^{(\rho)}_t \in E\big) + \int_0^t \int_{F}  p^{(\rho)}(s, x, z) \int_{ B_F(z,\rho')\setminus B(z,\rho)}  J(z,w)   \P^w \big(  X^{(\rho')}_{t-s} \in E\big) \,dw\, dz\, ds,\end{align*}where the second equality follows from the dominated convergence theorem, which is applicable   by Lemma \ref{l:blow-up-integrated-finite}(i).  
Since $	p^{(\rho')}(t,x,\cdot)$ and $p^{(\rho)}(t,x,\cdot) $ are  continuous, we have for all $y \in F$, 
\begin{align*}	p^{(\rho')}(t,x,y) = \lim_{\delta \to0} \frac{\P^x \big(  X^{(\rho')}_t \in B_{F}(y,\delta)\big)}{m_d(B_{F}(y,\delta))}&\le p^{(\rho)}(t,x,y)  + \limsup_{\delta\to 0}I(\delta),
\end{align*}
where 
\begin{align*}
&I(\delta)=I(\delta,t,x,y):=	\int_0^t   \int_{F}  p^{(\rho)}(s, x, z)\int_{ B_F(z,\rho')\setminus B(z,\rho)} \!\!\! \!\!\!\!\!\!J(z,w) \int_{B_{F}(y,\delta)}\!\!\! \frac{p^{(\rho')}(t-s, w, u)}{m_d(B_{F}(y,\delta))}   \,du  \,dw\, dz \, ds  .
\end{align*}
For $\eps \in (0,t/2)$, 
we write $I(\delta)=I_1(\delta,\eps) + I_2(\delta,\eps)$, where
\begin{align*}
	I_1(\delta,\eps):= \int_{0}^{t-\eps} \cdots ds \quad \text{and} \quad 
	I_2(\delta,\eps):= \int_{t-\eps}^t \cdots ds.
\end{align*} 
For all $\delta>0$, $\eps \in (0,t/2)$ and $s\in (0,t-\eps)$, by Proposition \ref{p:NDU}, we have
\begin{align*}
	\int_{B_{F}(y,\delta)} \frac{p^{(\rho')}(t-s, w, u)}{m_d(B_{F}(y,\delta))}   \,du  \le c_1\eps^{-d/\alpha} e^{t/(\rho' \wedge R_0)^\alpha}.
\end{align*}
Hence, applying the dominated convergence theorem, we obtain for all $\eps \in (0,t/2)$,
\begin{align*}
	\limsup_{\delta \to0} I_1(\delta,\eps)= \int_0^{t-\eps} \int_F p^{(\rho)}(s, x, z) \int_{B_F(z,\rho')\setminus B(z,\rho)} J(z,w)  p^{(\rho')}(t-s,w,y) \,dw\, dz\, ds.
\end{align*}
On the other hand, for all $\delta>0$ and $\eps \in (0,t/2)$, using Proposition \ref{p:NDU} and symmetry, we obtain
\begin{align*}
	I_2(\delta,\eps)&\le  \frac{c_2(t/2)^{-d/\alpha} e^{t/(\rho \wedge R_0)^{\alpha}}}{m_d(B_{F}(y,\delta)) }\int_{t-\eps}^t  \int_{B_{F}(y,\delta)}   \int_{F}    p^{(\rho')}(t-s,  u,w) \int_{ F\setminus B(w,\rho)}  J(w,z)   \,dz  \,dw\, du\, ds \\
	&\le c_2(t/2)^{-d/\alpha} e^{t/(\rho \wedge R_0)^{\alpha}} \sup_{u \in F}\int_{t-\eps}^t  \int_{F}    p^{(\rho')}(t-s,  u,w) \int_{ F\setminus B(w,\rho)}  J(w,z)   \,dz  \,dw\,  ds=:\wt I_2(\eps).
\end{align*}
By Lemma \ref{l:blow-up-integrated-finite}(ii), we have $\lim_{\eps \to 0} \wt I_2(\eps)=0$. Now, we arrive at
\begin{align*}
\limsup_{\delta\to 0} I(\delta)&\le \limsup_{\eps \to 0}  \left( 	\limsup_{\delta \to0} I_1(\delta,\eps) + \wt I_2(\eps)\right)\\
& = \int_0^{t} \int_F p^{(\rho)}(s, x, z) \int_{B_F(z,\rho')\setminus B(z,\rho)} J(z,w)  p^{(\rho')}(t-s,w,y) \,dw\, dz\, ds.
\end{align*} 
The proof is complete. \qed

\subsection{Off-diagonal estimates for truncated Dirichlet forms} \label{ss:4.2}

The following  result for  
truncated Dirichlet 
forms was established in \cite[Theorem 3.1]{GHL14}.

\begin{prop}\label{p:truncated-self-improvement}
	Let $\rho\in (0,\overline R)$ and let 
	$f(\cdot, \cdot): (0, \infty)\times (0,\infty)\to (0, 1]$ be such that $f(R,\cdot)$ is a non-decreasing function for any $R\in (0,\infty)$.
	Suppose that there exists a constant $T\in (0,\infty]$ such that for all $x_0 \in {F}$, $R>0$ and $t\in (0,T)$,
	\begin{align}\label{e:truncated-self-improvement-ass}
		\inf_{x\in B_{{F}}(x_0,R/4)}	\int_{B_{{F}}(x_0,R)} p^{(\rho),B_{{F}}(x_0,R)}(t,x,y) dy \ge 1- f(R,t).
	\end{align}
	Then for any $x\in F$, $R>\rho$, $t\in (0,T)$ and $n\ge 1$,
	\begin{align*}
		\int_{F\setminus B(x,nR)} p^{(\rho)}(t,x,y) dy \le f(R-\rho, t)^{n-1}.
	\end{align*}
\end{prop} 

Combining Proposition \ref{p:truncated-self-improvement} with Lemma \ref{l:NDL-modified}, we can obtain the next result, cf. \cite[Corollary 3.2]{GHL14}.

\begin{lem}\label{l:truncated-self-improvement-1}
	There exist $C>1$ and $a>0$ such that for any $\rho \in (0,\overline R)$ and $r>0$, we have
	\begin{align}\label{e:truncated-self-improvement-1}
		\sup_{x\in F}	\int_{F\setminus B(x,r/2)} p^{(\rho)}(t,x,y)dy \le C e^{-a r/\rho} \quad \text{for all $t\in (0, (\rho \wedge R_0)^\alpha]$}.
	\end{align}
\end{lem}
\pf If $r< 4\rho$, then by taking $C$ larger than $e^{4a}$, we obtain \eqref{e:truncated-self-improvement-1}. 

Suppose $r\ge 4\rho$ and let $N:=\sup\{ n\ge 1: 2n\rho \le r/2\}$. By Lemma \ref{l:NDL-modified}, there exists 
$c_1\in (0, 1)$
independent of $\rho$ such that for all $x_0\in F$, $R>0$ and $t\in (0, (\rho \wedge  R \wedge R_0)^\alpha]$,
\begin{align}\label{e:truncated-off-diagonal-0}  	\inf_{x\in B_{{F}}(x_0,R/4)} \int_{B_{{F}}(x_0,R)} p^{(\rho),B_{{F}}(x_0,R)}(t,x,y) dy \ge c_1.
\end{align}
Define 
\begin{align*}
	f(R,t):= \begin{cases}
		1-c_1 &\mbox{ if $t \le (\rho \wedge R \wedge R_0)^\alpha$},\\
		1 &\mbox{ otherwise}.
	\end{cases}
\end{align*}
By \eqref{e:truncated-off-diagonal-0}, \eqref{e:truncated-self-improvement-ass} holds with $T=\infty$. Thus, by Proposition \ref{p:truncated-self-improvement},  we have for all $t\in (0,(\rho \wedge R_0)^\alpha]$,
\begin{align*}
	&\sup_{x\in F}	\int_{F\setminus B(x,r/2)} p^{(\rho)}(t,x,y) dy \\
	&\le  \sup_{x\in F}	\int_{F\setminus B(x,2N\rho)} p^{(\rho)}(t,x,y) dy  \le f(\rho, t)^{N-1}  = (1-c_1)^{N-1} \le (1-c_1)^{r/(4\rho)-2}.
\end{align*}
This completes the proof. \qed

\begin{prop}\label{p:truncated-off-diagonal}
	There exist  $a_1,a_2>0$  such that  for any $\rho \in  (0,\overline R)$, 
	\begin{align*}
		p^{(\rho)}(t,x,y)  \le  a_1 t^{-d/\alpha}e^{-a_2|x-y|/\rho} \quad \text{for all $x,y\in {F}$ and $t\in (0,(\rho \wedge R_0)^\alpha]$.} 
	\end{align*}
\end{prop}
\pf   Let  $\rho \in  (0,\overline R)$, $x,y\in {F}$ and $r:=|x-y|$. 
Using the semigroup property, symmetry, Proposition \ref{p:NDU} and Lemma \ref{l:truncated-self-improvement-1}, we obtain for all $t\in (0,(\rho \wedge R_0)^\alpha]$,
\begin{align*}
	&p^{(\rho)}(t,x,y) \\
	&\le \int_{F\setminus B(y,r/2)} p^{(\rho)}(t/2,x,z) p^{(\rho)}(t/2,y,z)  dz + \int_{F\setminus B(x,r/2)} p^{(\rho)}(t/2,x,z) p^{(\rho)}(t/2,z,y)  dz\\
	&\le c_1(t/2)^{-d/\alpha} e^{t/(2(\rho \wedge R_0)^\alpha)} \bigg( \int_{F\setminus B(y,r/2)}  p^{(\rho)}(t/2,y,z)  dz + \int_{F\setminus B(x,r/2)} p^{(\rho)}(t/2,x,z)  dz \bigg) \\
	&\le c_2 t^{-d/\alpha}e^{-c_3r/\rho}.
\end{align*}
\qed

\begin{cor}\label{c:truncated-off-diagonal}
	For any $\delta\in (0,1)$ and $\lambda>0$,	there exists $C=C(\delta/\lambda)>0$   such that  for all $\rho \in  (0,\overline R)$,  $x,y\in {F}$ and $t\in (0,(\rho \wedge R_0)^\alpha]$ satisfying 
$\rho \le t^{\delta/\alpha}|x-y|^{1-\delta}$, 
	\begin{align*}
		p^{(\rho)}(t,x,y)  \le  C t^{-d/\alpha} \bigg( 1 \wedge \frac{t^{1/\alpha}}{|x-y|}\bigg)^\lambda.
	\end{align*}
\end{cor}
\pf 
Note that there exists  $c_1=c_1(\lambda/\delta)>0$ such that $(1\wedge a)^{\lambda/\delta} e^{1/a} \ge c_1$ for all $a>0$. Thus, by Proposition \ref{p:truncated-off-diagonal}, we have
\begin{align*}
	p^{(\rho)}(t,x,y)  \le  c_2 t^{-d/\alpha}e^{-c_3|x-y|/\rho} \le c_4 t^{-d/\alpha} \bigg(1\wedge \frac{\rho}{c_3|x-y|}\bigg)^{\lambda/\delta} \le c_5 t^{-d/\alpha} \bigg(1\wedge \frac{t^{1/\alpha}}{|x-y|}\bigg)^{\lambda} .
\end{align*}
\qed

Recall that $(\sE^{(\rho)},\sF)$ is conservative  for all $\rho \in (0,\overline R]$ (Proposition \ref{p:conservative}). Thus, we obtain the   following result from \cite[Lemma 6.6]{Ch25+}. 
\begin{lem}\label{l:truncated-self-improvement-reverse}
	Let $\rho \in (0,\overline R)$ and let $g$ be a positive function in $(0,\infty)$. Suppose that there exist  $k,T>0$  such that for all $x_0 \in {F}$, $R \ge k\rho$ and $t\in (0,T]$,
	\begin{align*}
		\sup_{x\in B_{F}(x_0,R/4)} \int_{F\setminus B(x_0,R)} p^{(\rho)}(t,x,y) dy  \le g(R).
	\end{align*}
	Then for any $x \in {F}$, $R\ge 3k\rho$ and $t\in (0,T]$, 
	\begin{align}\label{e:truncated-self-improvement-reverse}
		\int_{B_{F}(x,R)}  p^{(\rho),B_{F}(x,R)} (t,x,y) dy \ge 1- 2g(R/3).
	\end{align}
\end{lem}

\begin{remark}\label{r:truncated-self-improvement-reverse}
{\rm	 Since $B_{F}(x,3R/4)\subset B_{F}(x_0,R)$ for all $x_0\in {F}$ and $x\in B_{F}(x_0,R/4)$,  \eqref{e:truncated-self-improvement-reverse} implies that for all $x_0 \in {F}$, $R\ge 4k\rho$ and $t\in (0,T]$, 
	\begin{align*}
		\inf_{x\in B_{F}(x_0,R/4)}	\int_{B_{F}(x_0,R)}  p^{(\rho),B_{F}(x_0,R)} (t,x,y) dy \ge 1- 2g(R/4).
	\end{align*}
	}
\end{remark}

\begin{cor}\label{c:truncated-self-improvement-reverse}
	Let $\rho \in (0,\overline R)$ and let 
	$g(\cdot,\cdot):(0, \infty)\times (0,\infty) \to (0,1/2]$ be such that $g(R,\cdot)$ is 
	a non-decreasing function for any $R\in (0,\infty)$. 
	Suppose that there exists a constant  $k\ge 1$ such that  for all $x_0 \in {F}$, $R \ge k\rho$ and $t\in (0,(\rho \wedge R_0)^\alpha)$,
	\begin{align}\label{e:truncated-self-improvement-reverse-ass-fixed-t}
		\sup_{x\in B_{F}(x_0,R/4)} \int_{F\setminus B(x_0,R)} p^{(\rho)}(t,x,y) dy  \le g(R,t).
	\end{align}
	Then there exists $C>0$ independent of $\rho$ and $g$ such that for all $t\in (0,(\rho \wedge R_0)^\alpha)$,   and $x,y \in F$ with $|x-y|\ge 10k\rho$,
	\begin{align*}
		p^{(\rho)} (t,x,y) \le  C t^{-d/\alpha}(2g(k\rho,t/2))^{(|x-y|/(10k\rho)) -2} .
	\end{align*}
\end{cor}
\pf Let $t\in (0,(\rho \wedge R_0)^\alpha)$.
Since $g(R,\cdot)$ is non-decreasing, by \eqref{e:truncated-self-improvement-reverse-ass-fixed-t}, we have for all $x_0 \in F$, $R\ge k\rho$ and $s\in (0,t/2]$,
\begin{align*}
	\sup_{x\in B_{F}(x_0,R/4)} \int_{F\setminus B(x_0,R)} p^{(\rho)}(s,x,y) dy  \le g(R,s)\le g(R,t/2).
\end{align*}
Hence, by Lemma \ref{l:truncated-self-improvement-reverse} and Remark \ref{r:truncated-self-improvement-reverse}, we get that for all  $x_0 \in {F}$, $R\ge 4k\rho$ and $s\in (0,t/2]$, 
\begin{align}\label{e:self-improvment-reverse-1}
	\inf_{x\in B_{F}(x_0,R/4)}	\int_{B_{F}(x_0,R)}  p^{(\rho),B_{F}(x_0,R)} (s,x,y)dy \ge 1- 2g(R/4,t/2).
\end{align}
Define $f:(0,\infty) \times(0,\infty) \to (0,\infty)$ by
\begin{align*}
	f(R,s):= \begin{cases}
		2g(R/4,t/2)  &\mbox{ if $R\ge 4k\rho$ and $s \le t/2$},\\
		1 &\mbox{ otherwise}.
	\end{cases}
\end{align*}
By \eqref{e:self-improvment-reverse-1}, 
the assumption \eqref{e:truncated-self-improvement-ass} of Proposition \ref{p:truncated-self-improvement} holds  with $T=\infty$.  Thus, applying Proposition \ref{p:truncated-self-improvement} (with $R=(4k+1)\rho$), we obtain for all  $x\in {F}$, $s\in (0,t/2]$ and $n\ge 1$,
\begin{align}\label{e:truncated-self-improvement-g}
	\int_{F\setminus B(x,(4k+1)n\rho)} p^{(\rho)}(s,x,y) dy \le f( 4k\rho, s)^{n-1} = (2g(k\rho, t/2))^{n-1}.
\end{align} 

Now, let $x,y\in F$ with $|x-y|\ge 10k\rho$ and let $n\ge 1$ satisfy $10kn\rho\le |x-y|<10k(n+1)\rho$. Following the proof of Proposition \ref{p:truncated-off-diagonal} with Lemma \ref{l:truncated-self-improvement-1} replaced by  \eqref{e:truncated-self-improvement-g}, since $|x-y|/2 \ge 5kn\rho \ge (4k+1)n\rho$, we obtain 
\begin{align*}
	&	p^{(\rho)}(t,x,y) \\
	&\le c_1(t/2)^{-d/\alpha} e^{t/(2(\rho \wedge R_0)^\alpha)} \bigg( \int_{F\setminus B(y,|x-y|/2)}  p^{(\rho)}(t/2,y,z)  dz + \int_{F\setminus B(x,|x-y|/2)} p^{(\rho)}(t/2,x,z)  dz \bigg) \\
	&\le c_1(t/2)^{-d/\alpha} e^{t/(2(\rho \wedge R_0)^\alpha)} \bigg( \int_{F\setminus B(y,(4k+1)n\rho)}  p^{(\rho)}(t/2,y,z)  dz + \int_{F\setminus B(x,(4k+1)n\rho)} p^{(\rho)}(t/2,x,z)  dz \bigg) \\
	&\le c_2 t^{-d/\alpha}(2g(k\rho,t/2))^{n-1} \le c_2 t^{-d/\alpha}(2g(k\rho,t/2))^{(|x-y|/(10k\rho)) -2}.
\end{align*}
The proof is complete. \qed 

\section{Sharp heat kernel estimates}\label{s:hke}
\subsection{Lower estimate for 
 $t\in (0, R_0^\alpha)$
}\label{ss:t-lhke}
In this subsection, 
we prove the lower estimates in Theorem \ref{t:main-HKE} for $t<R_0^\alpha$.
Throughout this section, $\delta_3\in (0,1)$ stands for the constant in Lemma \ref{l:NDL-modified}.
\begin{lem}\label{l:HKE-lower-1}
	There exists $C>0$ such that  for all  $t\in (0, R_0^\alpha)$ and  $x,y\in F$  with $|x-y|< \delta_3t^{1/\alpha}$,
	\begin{align*}
		p(t,  x, y) \ge Ct^{-d/\alpha}.
	\end{align*}
\end{lem}
\pf  By \eqref{e:NDL-modified-1}, we have  for all  $t\in (0, R_0^\alpha)$ and  $x,y\in F$  with $|x-y|< \delta_3t^{1/\alpha}$,
\begin{align*}
	p(t, x, y) =	p^{(\overline R), B_{F}(x, 2\overline R)}(t, x, y)  \ge c_1 t^{-d/\alpha}.
\end{align*}
\qed

The jump kernel of $(\sE,\sF)$ has a density $J(x,y)$ with respect to the Lebesgue measure. 
Hence, repeating the argument in  
the paragraph
above \cite[Lemma 4.7]{CK03} and the proof of  
\cite[Lemma 4.7]{CK03}, the 
process $X$ satisfies  the following  L\'evy system formula on $F\times F$: for  any stopping time $\sigma$ and    any  Borel function $f: F \times  F \to [0,\infty]$ vanishing on the set $\{(x,x):x \in F\}$,  
\begin{align}\label{e:levy-system}
	\E^x\sum_{s\in (0, \sigma]}f(X_{s-},X_s)= \E^x\int^{\sigma}_0\int_{F}f(X_s, y) J(x,y)   dy ds, \quad x\in F.
\end{align}

  Define for $x\in F$ and $t>0$,
\begin{align}\label{e:def-deltaDt}
	\delta_\Sigma(x,t)&:=\delta_\Sigma(x) \vee t^{1/\alpha}.
\end{align}

\begin{lem}\label{l:HKE-lower-2}
	There exists $C>0$ such that  for all  $t\in (0, R_0^\alpha)$ and  $x,y\in F$  with $|x-y|\ge \delta_3t^{1/\alpha}$,
	\begin{align*}
		p(t, x, y) \ge\frac{C t }{|x-y|^{d+\alpha}}\Phi\left( \frac{(|x-y| \wedge A_0)^2}{(\delta_\Sigma(x,t) \wedge A_0)(\delta_\Sigma(y,t)\wedge A_0) }\right).
	\end{align*}
\end{lem}
\pf Set $\eps_0:=2^{-2-1/\alpha}\delta_3$.  
By Lemma \ref{l:HKE-lower-1}, we have
\begin{align}\label{e:lower-offdiagonal-1} 
	\inf_{s\in [t/2,t],\, w \in B_{F}(y,  \eps_0t^{1/\alpha})}	p(s,  w, y)  \ge c_1 t^{-d/\alpha}. 
\end{align} 
Note that for all $z\in B_{F}(x,  \eps_0t^{1/\alpha})$ and $w\in B_{F}(y, \eps_0t^{1/\alpha})$,
\begin{align*}
	& 	|z-w|  \ge |x-y| - 2 \eps_0t^{1/\alpha} \ge (1/2)|x-y|,\\
&	|z-w|  \le |x-y| +2\eps_0 t^{1/\alpha} \le (3/2)|x-y|,\\
&	\delta_\Sigma(z)\le \delta_\Sigma(x)  + \eps_0t^{1/\alpha}  \le 2\delta_\Sigma(x,t),\quad \;\;
		\delta_\Sigma(w) \le  2\delta_\Sigma(y,t).
\end{align*}
Thus, by \eqref{e:JJ}, \eqref{e:Phi-scaling} and \eqref{e:without-A0},  we have  for all $z\in B_{F}(x, \eps_0 t^{1/\alpha})$ and $w\in B_{F}(y,  \eps_0t^{1/\alpha})$,
\begin{align}\label{e:lower-offdiagonal-2}
	J(z,w)& \ge \frac{C_1^{-1}}{((3/2)|x-y|)^{d+\alpha}} \Phi  \bigg( \frac{( (1/2)|x-y| \wedge A_0 )^2}{( 2\delta_\Sigma(x,t) \wedge A_0)( 2\delta_\Sigma(y,t) \wedge A_0)} \bigg)  \nn\\
	&\ge \frac{c_2}{|x-y|^{d+\alpha}} \Phi  \bigg( \frac{( |x-y| \wedge A_0 )^2}{( \delta_\Sigma(x,t) \wedge A_0)(  \delta_\Sigma(y,t) \wedge A_0)} \bigg).
\end{align}

Using the strong Markov property in the first line below, \eqref{e:lower-offdiagonal-1} in the second,  the L\'evy system formula \eqref{e:levy-system} in the third,    \eqref{e:lower-offdiagonal-2} in the fourth, and  
\eqref{e:d-set} in the fifth, we get
\begin{align}\label{e:lower-offdiagonal-main}
	&p(t,x,y) \nn\\
	&\ge \E^x \Big[ p(t- \tau_{B_{F}(x, \eps_0t^{1/\alpha})}, X_{\tau_{B_{F}(x,\eps_0t^{1/\alpha})}}, y);\,  \tau_{B_{F}(x, \eps_0t^{1/\alpha})}  \le t/2, \, X_{\tau_{B_{F}(x,\eps_0t^{1/\alpha})}} \in B_{F}(y,\eps_0t^{1/\alpha})  \Big] \nn\\
	&\ge c_1t^{-d/\alpha} \P^x \Big( \tau_{B_{F}(x,\eps_0t^{1/\alpha})} \le t/2, \, X_{\tau_{B_{F}(x,\eps_0t^{1/\alpha})}} \in B_{F}(y,\eps_0t^{1/\alpha})  \Big)\nn\\
	&\ge c_1t^{-d/\alpha} \E^x  \int_{0}^{\tau_{B_{F}(x,\eps_0t^{1/\alpha})} \wedge (t/2)} \int_{B_{F}(y,\eps_0t^{1/\alpha})} J(X_s, w) dw\,ds \nn\\
	&\ge  \frac{c_3 t^{-d/\alpha}}{|x-y|^{d+\alpha}} \Phi  \bigg( \frac{( |x-y| \wedge A_0 )^2}{( \delta_\Sigma(x,t) \wedge A_0)(  \delta_\Sigma(y,t) \wedge A_0)} \bigg) \E^x [\tau_{B_{F}(x,\eps_0t^{1/\alpha})} \wedge (t/2)] m_d(B_{F}(y,\eps_0t^{1/\alpha}) )\nn\\
	&\ge  \frac{c_4}{|x-y|^{d+\alpha}} \Phi  \bigg( \frac{( |x-y| \wedge A_0 )^2}{( \delta_\Sigma(x,t) \wedge A_0)(  \delta_\Sigma(y,t) \wedge A_0)} \bigg) \E^x [\tau_{B_{F}(x,\eps_0t^{1/\alpha})} \wedge (t/2)] .
\end{align}
It follows from  Lemma \ref{l:NDL-modified}  that
\begin{align}\label{e:lower-offdiagonal-main-2}
	&\E^x [\tau_{B_{F}(x,\eps_0t^{1/\alpha})} \wedge (t/2)] \ge (t/2) \P^x(\tau_{B_{F}(x,\eps_0t^{1/\alpha})} \ge t/2 ) \ge (t/2) \P^x\big(X^{B_{F}(x,\eps_0t^{1/\alpha})}_{t/2} \in B_{F}(x,\eps_0t^{1/\alpha}) \big) \nn\\
	&= (t/2) \int_{B_{F}(x,\eps_0t^{1/\alpha})} p^{B_{F}(x,\eps_0t^{1/\alpha})}(t/2, x, z)dz \ge  c_5 t e^{-c_6t / (2(\eps_0t^{1/\alpha})^\alpha)} =c_7t.
\end{align} 
Combining \eqref{e:lower-offdiagonal-main} with \eqref{e:lower-offdiagonal-main-2}, we obtain the desired result. \qed

The lower bound in Theorem \ref{t:main-HKE}   for $t\in (0, R_0^\alpha)$ now follows from Lemmas \ref{l:HKE-lower-1} and \ref{l:HKE-lower-2}.

\subsection{Upper estimate  for $t\in (0, R_0^\alpha)$}\label{ss:t-uhke}
In this subsection, we 
prove the upper estimates in Theorem \ref{t:main-HKE}  for $t\in (0, R_0^\alpha)$. 
	From \eqref{e:beta1}  and the fact $\text{dim}_{\rm A}(\Sigma )<d$, we have $d-\beta_1 \ge (\gamma\wedge \alpha)-\beta_1>0$.
	
	\begin{lem}\label{l:Meyer-integral-estimate-1}
		For any   $\lambda \in [0, d-\beta_1)$, there exists $C=C(\lambda)>0$ such that for all $0<t^{1/\alpha}<\rho<\overline R$, $x\in F$  and $w\in F$, 
		\begin{align*}
			\int_{F\setminus B(w,\rho)}  t^{-d/\alpha} \bigg(1 \wedge \frac{t^{1/\alpha}}{|x-z|}\bigg)^\lambda J(w,z) dz  \le \frac{Ct^{(\lambda-d)/\alpha}}{\rho^{\lambda+\alpha}}  \Phi\bigg( \frac{\rho \wedge A_0}{\delta_\Sigma(w) \wedge A_0}\bigg).
		\end{align*}
	\end{lem}
	\pf It follows from  Lemma \ref{l:tail-estimate}  that 
	\begin{align}\label{e:Meyer-integral-estimate-1}
		&	 \int_{F \setminus (B(x,\rho) \cup B(w,\rho))}  t^{-d/\alpha} \bigg(1 \wedge \frac{t^{1/\alpha}}{|x-z|}\bigg)^\lambda J(w,z) dz\nn\\
		& \le  t^{-d/\alpha} \bigg( \frac{t^{1/\alpha}}{\rho}\bigg)^\lambda\int_{F \setminus B(w,\rho)} J(w,z)dz \le \frac{c_1 t^{(\lambda-d)/\alpha}}{\rho^{\lambda+\alpha}} \Phi\bigg( \frac{\rho \wedge A_0}{\delta_\Sigma(w) \wedge A_0}\bigg).
	\end{align}
	On the other hand, using   Lemma \ref{l:tail-estimate-2}(i), we get that
	\begin{align}\label{e:Meyer-integral-estimate-2}
		&	 \int_{B_F(x,t^{1/\alpha}) \setminus B(w,\rho)}  t^{-d/\alpha} \bigg(1 \wedge \frac{t^{1/\alpha}}{|x-z|}\bigg)^\lambda J(w,z) dz\nn\\
		& \le  t^{-d/\alpha}\int_{B_F(x,t^{1/\alpha}) \setminus B(w,\rho)} J(w,z)dz \le \frac{c_2}{\rho^{d+\alpha}} \Phi\bigg( \frac{(\rho \wedge A_0)^2}{(\delta_\Sigma(w) \wedge A_0)(	\delta_\Sigma(x,t) \wedge A_0)}\bigg),
	\end{align}
	and that for all $n\ge 1$,
	\begin{align}\label{e:Meyer-integral-estimate-3}
		&	 \int_{B_F(x,2^nt^{1/\alpha}) \setminus (B(x,2^{n-1}t^{1/\alpha}) \cup B(w,\rho))}  t^{-d/\alpha} \bigg(1 \wedge \frac{t^{1/\alpha}}{|x-z|}\bigg)^\lambda  J(w,z) dz\nn\\
		& \le 2^{-(n-1)\lambda} t^{-d/\alpha}\int_{B_F(x,2^nt^{1/\alpha}) \setminus B(w,\rho)} J(w,z)dz \nn\\
		&\le \frac{c_32^{n(d-\lambda)} }{\rho^{d+\alpha}}  \Phi\bigg( \frac{(\rho \wedge A_0)^2}{(\delta_\Sigma(w) \wedge A_0)(	\delta_\Sigma(x,2^{n\alpha}t) \wedge A_0)}\bigg).
	\end{align}
	Let $N:=\min\{n\ge 1: 2^{n}t^{1/\alpha} \ge \rho\}$.  Using  \eqref{e:Meyer-integral-estimate-2} and \eqref{e:Meyer-integral-estimate-3}  in the second inequality below, \eqref{e:Phi-scaling}, \eqref{e:without-A0} and 
	 $2^{N-1}t^{1/\alpha}<\rho$ in the third, $d-\lambda-\beta_1>0$ in the fourth, and $2^{N+1}<4\rho/t^{1/\alpha}$  in the fifth, we obtain 
	\begin{align*}
		&\int_{B_F(x,\rho) \setminus  B(w,\rho)}  t^{-d/\alpha} \bigg(1 \wedge \frac{t^{1/\alpha}}{|x-z|}\bigg)^\lambda J(w,z) dz\nn\\
		&\le  \int_{B_F(x,t^{1/\alpha}) \setminus B(w,\rho)}  t^{-d/\alpha} \bigg(1 \wedge \frac{t^{1/\alpha}}{|x-z|}\bigg)^\lambda J(w,z) dz\nn\\
		&\quad +  \sum_{n=1}^N  \int_{B_F(x,2^nt^{1/\alpha}) \setminus (B(x,2^{n-1}t^{1/\alpha}) \cup B(w,\rho))}  t^{-d/\alpha} \bigg(1 \wedge \frac{t^{1/\alpha}}{|x-z|}\bigg)^\lambda J(w,z) dz\nn\\
		&\le \frac{c_4}{\rho^{d+\alpha}}\sum_{n=0}^N 2^{n(d-\lambda)} \Phi\bigg( \frac{(\rho \wedge A_0)^2}{(\delta_\Sigma(w) \wedge A_0)(	2^{n} t^{1/\alpha} \wedge A_0)}\bigg)\nn\\
	&\le \frac{c_5}{\rho^{d+\alpha}}  \Phi\bigg( \frac{\rho \wedge A_0}{\delta_\Sigma(w) \wedge A_0}\bigg)\sum_{n=0}^N 2^{n(d-\lambda)} \bigg( \frac{\rho}{2^{n-1}t^{1/\alpha}}\bigg)^{\beta_1}\nn\\
	&\le \frac{c_62^{(N+1)(d-\lambda-\beta_1)}}{t^{\beta_1/\alpha}\rho^{d+\alpha-\beta_1}}  \Phi\bigg( \frac{\rho \wedge A_0}{\delta_\Sigma(w) \wedge A_0}\bigg) \le \frac{c_7t^{(\lambda-d)/\alpha}}{\rho^{\lambda+\alpha}}  \Phi\bigg( \frac{\rho \wedge A_0}{\delta_\Sigma(w) \wedge A_0}\bigg) .
	\end{align*}
	Combining this with \eqref{e:Meyer-integral-estimate-1}, we arrive at the result.
	\qed

	Recall that $p^{(\overline R)}(t,x,y)$ denotes $p(t,x,y)$, and the constant $\Lambda$ is defined in \eqref{e:def-Lambda}.

	\begin{lem}\label{l:UHK-without-Phi}
	If there exist  $\lambda\in [0,d- \beta_1)$ and $C>0$ such that for all $\rho \in (0,\overline R]$, $x,y\in {F}$ and $t\in (0,(\rho \wedge R_0)^\alpha)$,
		\begin{align}\label{e:UHK-without-Phi-ass}
			p^{(\rho)}(t,x,y) \le C t^{-d/\alpha} \bigg(1 \wedge \frac{t^{1/\alpha}}{|x-y|}\bigg)^\lambda,
		\end{align}
		then for any $\eps \in (0,\alpha-\beta_1)$,  there exists  $C'=C'(\lambda, C, \eps )>0$ such that for all $\rho \in (0,\overline R]$, $x,y\in {F}$ and $t\in (0,(\rho \wedge R_0)^\alpha)$,
		\begin{align}\label{e:UHK-without-Phi-result}
				p^{(\rho)}(t,x,y) \le C' t^{-d/\alpha} \bigg(1 \wedge \frac{t^{1/\alpha}}{|x-y|}\bigg)^{\lambda+\eps }.
		\end{align}
	\end{lem}
	\pf   Let $\rho \in (0,\overline R]$, $x,y\in {F}$,  $t\in (0,(\rho \wedge R_0)^\alpha)$ and let $\delta \in (0,1)$ satisfy
	 \begin{align}\label{e:choice-delta}
	 	(1-\delta)(\lambda + \alpha- \beta_1)   =\lambda +  \eps.
	 \end{align}
	  Set $r:=|x-y|$ and  $\sigma:= t^{\delta/\alpha}r^{1-\delta} $.
 If $r\le \Lambda^{1/(1-\delta)}t^{1/\alpha}$, then \eqref{e:UHK-without-Phi-result}  follows from Proposition \ref{p:NDU}, and if $\rho\le \sigma$, then \eqref{e:UHK-without-Phi-result} follows from Corollary  \ref{c:truncated-off-diagonal}.

Suppose that $r> \Lambda^{1/(1-\delta)}t^{1/\alpha}$  and $\rho>\sigma$.  Then, by the definition of $\sigma$, we have $t<(\sigma/\Lambda)^\alpha$.	Applying Corollary  \ref{c:truncated-off-diagonal} with $\lambda$ replaced by $\lambda+\eps$, we get 
\begin{align*}	p^{(\sigma)}(t,x,y)\le  \frac{c_1t^{(\lambda+\eps-d)/\alpha}}{ r^{\lambda+\eps}}.
	\end{align*} 
Thus, by Proposition \ref{p:comparison-Meyer}, to obtain \eqref{e:UHK-without-Phi-result},  it suffices to show that
\begin{align*}
	I:=\int_0^t \int_F p^{(\sigma)}(s, x, z) \int_{B_F(z,\rho)\setminus B(z,\sigma)} J(z,w)  p^{(\rho)}(t-s,w,y) \,dw dz ds\le  \frac{c_2t^{(\lambda+\eps-d)/\alpha}}{ r^{\lambda+\eps}}
\end{align*} 
for some constant $c_2>0$ independent of $\rho,x,y,t$.	

 Since $t<(\sigma \wedge R_0)^\alpha \le (\rho \wedge R_0)^\alpha$, using  \eqref{e:UHK-without-Phi-ass}, symmetry and a change of variables, we get 
	\begin{align*}
		I&\le c_3 \int_0^{t/2} \int_F \int_{B_F(z,\rho)\setminus B(z,\sigma)}  p^{(\sigma)}(s,x,z) t^{-d/\alpha} \bigg(1 \wedge \frac{t^{1/\alpha}}{|y-w|}\bigg)^\lambda J(z,w) \,dw\, dz\, ds\\
		&\quad + c_3 \int_{t/2}^{t} \int_F \int_{B_F(z,\rho)\setminus B(z,\sigma)} t^{-d/\alpha} \bigg(1 \wedge \frac{t^{1/\alpha}}{|x-z|}\bigg)^\lambda  J(z,w)  p^{(\rho)}(t-s,w,y) \,dw\, dz\, ds\nn\\
	&\le c_3 \int_0^{t/2} \int_F p^{(\sigma)}(s,x,z)\int_{F\setminus B(z,\sigma)}  t^{-d/\alpha} \bigg(1 \wedge \frac{t^{1/\alpha}}{|y-w|}\bigg)^\lambda J(z,w) \,dw\, dz\, ds\\
	&\quad + c_3 \int_{0}^{t/2} \int_F  p^{(\rho)}(s,y,w)  \int_{F\setminus B(w,\sigma)} t^{-d/\alpha} \bigg(1 \wedge \frac{t^{1/\alpha}}{|x-z|}\bigg)^\lambda J(w,z)  \,dz\, dw\, ds\nn\\
	&=:c_3 (I_1+I_2).
	\end{align*}
	Since $\lambda<d-\beta_1$ and $t<(\sigma/\Lambda)^\alpha$, applying Lemma \ref{l:Meyer-integral-estimate-1} and  Proposition \ref{p:blow-up-integrated}, we obtain
	\begin{align*}
		I_1&\le \frac{c_4t^{(\lambda-d)/\alpha}}{\sigma^{\lambda+\alpha}}  \int_0^{t/2} \int_F p^{(\sigma)}(s,x,z)	 \Phi\bigg( \frac{\sigma \wedge A_0}{\delta_\Sigma(z) \wedge A_0}\bigg)  dz\, ds\\
		&\le \frac{c_5t^{(\lambda-d)/\alpha}}{\sigma^{\lambda+\alpha}}  \bigg( \frac{\sigma}{t^{1/\alpha}}\bigg)^{\beta_1} \int_0^{t/2} \int_F p^{(\sigma)}(s,x,z)	 \Phi\bigg( \frac{t^{1/\alpha} \wedge A_0}{\delta_\Sigma(z) \wedge A_0}\bigg)  dz\, ds\\
		&\le \frac{c_6t^{1+(\lambda-d-\beta_1)/\alpha}}{\sigma^{\lambda+\alpha-\beta_1}}   = \frac{c_6 t^{(\lambda+\eps-d)/\alpha}}{r^{\lambda+\eps}},
	\end{align*}
	where we used \eqref{e:Phi-scaling} with \eqref{e:without-A0} in the second inequality, and \eqref{e:choice-delta} in the equality. Similarly, since $t<(\rho/\Lambda)^\alpha$, applying Lemma \ref{l:Meyer-integral-estimate-1} and Proposition \ref{p:blow-up-integrated}, we get that
	\begin{align*}
		I_2&\le  \frac{c_7t^{(\lambda-d)/\alpha}}{\sigma^{\lambda+\alpha}}  \int_0^{t/2} \int_F p^{(\rho)}(s,y,w)	 \Phi\bigg( \frac{\sigma \wedge A_0}{\delta_\Sigma(w) \wedge A_0}\bigg)  dw\, ds\\
		&\le \frac{c_8t^{(\lambda-d)/\alpha}}{\sigma^{\lambda+\alpha}}  \bigg( \frac{\sigma}{t^{1/\alpha}}\bigg)^{\beta_1} \int_0^{t/2} \int_F p^{(\rho)}(s,y,w)	 \Phi\bigg( \frac{t^{1/\alpha} \wedge A_0}{\delta_\Sigma(w) \wedge A_0}\bigg)  dw\, ds\\
		&\le \frac{c_{9}t^{1+(\lambda-d-\beta_1)/\alpha}}{\sigma^{\lambda+\alpha-\beta_1}}   = \frac{c_{9} t^{(\lambda+\eps-d)/\alpha}}{r^{\lambda+\eps}}.
	\end{align*} 
	The proof is complete. \qed 
	
	By Proposition \ref{p:NDU},  condition \eqref{e:UHK-without-Phi-ass} in Lemma \ref{l:UHK-without-Phi} holds with $\lambda=0$. Applying  Lemma \ref{l:UHK-without-Phi}  finitely many times, we obtain the following result.
	
	\begin{lem}\label{l:UHK-without-Phi-iterated}
	For any $\lambda \in [0, d+\alpha-2\beta_1)$, there exists $C=C(\lambda)>0$ such that for all $\rho \in (0,\overline R]$, $x,y\in {F}$ and $t\in (0,(\rho \wedge R_0)^\alpha)$,
	\begin{align*}
		p^{(\rho)}(t,x,y) \le C t^{-d/\alpha} \bigg(1 \wedge \frac{t^{1/\alpha}}{|x-y|}\bigg)^{\lambda }.
	\end{align*}
	\end{lem}

	\begin{lem}\label{l:blow-up-integrated-improve-1}	There exists $C>0$ such that for all $a>0$, $\rho \in (0,\overline R]$,  $x\in {F}$ and $t\in (0, (\rho \wedge R_0)^\alpha)$,
		\begin{align*}
			\int_F p^{(\rho)}(t, x, z) \Phi \bigg( \frac{a}{\delta_\Sigma(z)\wedge A_0}\bigg)   dz\,\le C\Phi \bigg( \frac{a}{t^{1/\alpha}\wedge A_0}\bigg) .
		\end{align*}
	\end{lem}
	\pf We have
	\begin{align*}
		&\int_F p^{(\rho)}(t, x, z) \Phi \bigg( \frac{a}{\delta_\Sigma(z)\wedge A_0}\bigg)   dz\\
		&\le  \Phi \bigg( \frac{a}{t^{1/\alpha}\wedge A_0}\bigg)	\int_{F:\delta_\Sigma(z)\ge t^{1/\alpha}} p^{(\rho)}(t, x, z)   dz + 	\int_{F: \delta_\Sigma(z)<t^{1/\alpha}} p^{(\rho)}(t, x, z) \Phi \bigg( \frac{a}{\delta_\Sigma(z)\wedge A_0}\bigg)   dz\\
		&\le \Phi \bigg( \frac{a}{t^{1/\alpha}\wedge A_0}\bigg)	 + 	\int_{F: \delta_\Sigma(z)<t^{1/\alpha}} p^{(\rho)}(t, x, z) \Phi \bigg( \frac{a}{\delta_\Sigma(z)\wedge A_0}\bigg) \, dz.
	\end{align*}
	Thus, it suffices to estimate $I:=\int_{F: \delta_\Sigma(z)<t^{1/\alpha}} p^{(\rho)}(t, x, z) \Phi \big( \frac{a}{\delta_\Sigma(z)\wedge A_0}\big) \, dz$.   Using Lemma \ref{l:UHK-without-Phi-iterated} with $\lambda=d-\beta_1$ in the first inequality below, \eqref{e:Phi-scaling}, \eqref{e:without-A0} and  Proposition \ref{p:Aikawa} with $q=(\beta_1+\gamma)/2$ in the second inequality, and $\beta_1<\gamma$ in the second equality, we get that
	\begin{align*}
		I
		&= \sum_{n=1}^\infty   \int_{B_F(x, t^{1/\alpha}): \,  \delta_\Sigma(z)\in [2^{-n} t^{1/\alpha}, 2^{1-n}t^{1/\alpha})} p^{(\rho)}(t, x, z) \Phi \bigg( \frac{a}{\delta_\Sigma(z)\wedge A_0}\bigg) \, dz\\
		&\quad +  \sum_{n=1}^\infty  \sum_{m=1}^\infty \int_{B_F(x, 2^{m}t^{1/\alpha}) \setminus B(x, 2^{m-1}t^{1/\alpha}) :  \delta_\Sigma(z)\in [2^{-n} t^{1/\alpha}, 2^{1-n}t^{1/\alpha})} p^{(\rho)}(t, x, z) \Phi \bigg( \frac{a}{\delta_\Sigma(z)\wedge A_0}\bigg) \, dz \\
		&\le c_1t^{-d/\alpha} \sum_{n=1}^\infty \Phi \bigg( \frac{a}{2^{-n}t^{1/\alpha}\wedge A_0}\bigg)   \int_{B_F(x, t^{1/\alpha}): \,  \delta_\Sigma(z)< 2^{1-n}t^{1/\alpha}} dz\\
		&\quad + c_1t^{-d/\alpha} \sum_{n=1}^\infty  \sum_{m=1}^\infty  2^{-(m-1)(d-\beta_1)} \Phi \bigg( \frac{a}{2^{-n}t^{1/\alpha}\wedge A_0}\bigg)  \int_{B_F(x, 2^{m}t^{1/\alpha}) :  \delta_\Sigma(z)< 2^{1-n}t^{1/\alpha}}   dz \\
		&\le c_2\Phi \bigg( \frac{a}{t^{1/\alpha}\wedge A_0}\bigg)  \left(  \sum_{n=1}^\infty 2^{n\beta_1-n(\beta_1+\gamma)/2}+ \sum_{n=1}^\infty  \sum_{m=1}^\infty  2^{n\beta_1 -n (\beta_1+\gamma)/2-m(d-\beta_1) + m(d-(\beta_1+\gamma)/2)}\right)  \\
		&= c_3\Phi \bigg( \frac{a}{t^{1/\alpha}\wedge A_0}\bigg) .
	\end{align*}
	The proof is complete.
	\qed
	
	\begin{prop}\label{p:UHK-tt}
		For any $\lambda \in (0, d+\alpha)$, there exists $C=C(\lambda)>0$ such that for all $\rho \in (0,\overline R]$, $x,y\in {F}$ and $t\in (0,(\rho \wedge R_0)^\alpha)$,
		\begin{align}\label{e:UHK-tt}
			p^{(\rho)}(t,x,y) \le C t^{-d/\alpha} \bigg(1 \wedge \frac{t^{1/\alpha}}{|x-y|}\bigg)^{\lambda} \Phi\bigg( \frac{(|x-y|\wedge A_0)^2}{(t^{1/\alpha}\wedge A_0)^2}\bigg).
		\end{align}
	\end{prop}
	\pf 
	Let $\rho \in (0,\overline R]$, $x,y\in F$ and $t\in (0,(\rho \wedge R_0)^\alpha)$.  Set $r:=|x-y|$, $\delta:=(d+\alpha-\lambda)/(d+\alpha)$,  and $\sigma:=t^{\delta/\alpha}r^{1-\delta}$. 		
	Following the arguments in the first paragraph of the  proof of Lemma \ref{l:UHK-without-Phi},  we see that, to obtain \eqref{e:UHK-tt}, it suffices to show that for the case $t^{1/\alpha}<\sigma<\rho$,
	\begin{align*}
	I:=	\int_0^t \int_F p^{(\sigma)}(s, x, z) \int_{B_F(z,\rho)\setminus B(z,\sigma)} J(z,w)  p^{(\rho)}(t-s,w,y) \,dw\, dz\, ds
	\end{align*}
	is bounded above by the right-hand side of \eqref{e:UHK-tt}.
	
	Suppose $t^{1/\alpha}<\sigma<\rho$.		
	 Using Lemma \ref{l:jump-density-monotonicty} and symmetry in the first inequality below,   Lemma \ref{l:blow-up-integrated-improve-1}  in the second and  third,   \eqref{e:Phi-scaling} and  \eqref{e:without-A0} in the fourth, and $\beta_1<\alpha$ and the monotonicity of $\Phi$ with $\sigma<r$ in the last, we obtain 
\begin{align*}
	I&\le \frac{c_1}{\sigma^{d+\alpha}}\int_0^t \int_F p^{(\sigma)}(s, x, z) \int_{B_F(z,\rho)\setminus B(z,\sigma)}  p^{(\rho)}(t-s,y,w) \Phi\bigg( \frac{(\sigma \wedge A_0)^2}{(\delta_\Sigma(z) \wedge A_0) (\delta_\Sigma(w) \wedge A_0)}\bigg)  dw dz ds\\
	&\le \frac{c_2}{\sigma^{d+\alpha}}\int_0^t \int_F p^{(\sigma)}(s, x, z)  \Phi\bigg( \frac{(\sigma \wedge A_0)^2}{(\delta_\Sigma(z) \wedge A_0) ((t-s)^{1/\alpha} \wedge A_0)}\bigg)   dz\, ds\\
	&\le  \frac{c_3}{\sigma^{d+\alpha}}\int_0^t  \Phi\bigg( \frac{(\sigma \wedge A_0)^2}{(s^{1/\alpha} \wedge A_0) ((t-s)^{1/\alpha} \wedge A_0)}\bigg)   ds\\
		&\le  \frac{c_4 t^{2\beta_1/\alpha}}{\sigma^{d+\alpha}} \Phi\bigg( \frac{(\sigma \wedge A_0)^2}{(t^{1/\alpha} \wedge A_0)^2}\bigg) \int_0^t s^{-\beta_1/\alpha} (t-s)^{-\beta_1/\alpha}   ds \le \frac{c_5t}{\sigma^{d+\alpha}} \Phi\bigg( \frac{(r \wedge A_0)^2}{(t^{1/\alpha} \wedge A_0)^2}\bigg) .
\end{align*}	
 Since $t^{1/\alpha}<\sigma<r$ and $t\sigma^{-d-\alpha}= t^{(\lambda-d)/\alpha} r^{-\lambda}$, the proof is complete. \qed

	Note that for any fixed $a>0$, 
	\begin{align}\label{e:r-r2-monotone}
	\text{the function }  r\mapsto	\frac{r^2}{a\wedge r} \text{ is increasing on $(0,\infty)$}.
	\end{align}

	Recall that $\delta_\Sigma(x,t)$ is defined in \eqref{e:def-deltaDt}. 
	Define for $x,y \in F$ and $t>0$,
	\begin{align*}
		\delta^\vee_\Sigma(x,y,t) := \delta_\Sigma(x,t)\vee \delta_\Sigma(y,t)=
		\delta_\Sigma(x) \vee \delta_\Sigma(y) \vee t^{1/\alpha}.
	\end{align*}

	\begin{lem}\label{l:Meyer-integral-estimate-with-Phi}
	For any   $\lambda \in [0, d+\beta_1)$, there exists $C=C(\lambda)>0$ such that for all $0<t^{1/\alpha}<\rho<r<\overline R$ and  $x,w\in F$ 
		\begin{align*}
			&\int_{B_F(x, r)\setminus B(w,\rho)}  t^{-d/\alpha} \bigg(1 \wedge \frac{t^{1/\alpha}}{|x-z|}\bigg)^\lambda  \Phi\bigg( \frac{(|x-z| \wedge A_0)^2}{(\delta^\vee_\Sigma(x,z,t) \wedge |x-z| \wedge A_0) (t^{1/\alpha}\wedge A_0)}\bigg)  J(w,z) dz\\
			&  \le \frac{Ct^{(\lambda-d-\beta_1)/\alpha}r^{d+\beta_1-\lambda}}{\rho^{d+\alpha}} \Phi\left( \frac{(\rho \wedge A_0)^2}{(\delta_\Sigma(w) \wedge A_0) (\delta_\Sigma(x,t) \wedge A_0) }\right).
		\end{align*}
	\end{lem}
	\pf 
	Let $0<t^{1/\alpha}<\rho<r<\overline R$ and $x, w\in F$, and define 
	$$
	f(z):=t^{-d/\alpha} \bigg(1 \wedge \frac{t^{1/\alpha}}{|x-z|}\bigg)^\lambda  \Phi\bigg( \frac{(|x-z| \wedge A_0)^2}{(\delta^\vee_\Sigma(x,z,t) \wedge |x-z| \wedge A_0) (t^{1/\alpha}\wedge A_0)}\bigg)  J(w,z), \quad z \in F.
	$$
	  By \eqref{e:r-r2-monotone}, we have for all $z\in B_F(x,t^{1/\alpha})$,
	\begin{align*}
f(z)	&\le 	t^{-d/\alpha}  \Phi\bigg( \frac{(t^{1/\alpha} \wedge A_0)^2}{(\delta^\vee_\Sigma(x,z,t) \wedge t^{1/\alpha} \wedge A_0) (t^{1/\alpha}\wedge A_0)}\bigg) = \Phi(1)t^{-d/\alpha}.
	\end{align*}
	Hence, using   Lemma \ref{l:tail-estimate-2}(i), we get that
	\begin{align}\label{e:Meyer-integral-estimate-with-Phi-1}
		&	
		 \int_{B_F(x,t^{1/\alpha}) \setminus B(w,\rho)}   f(z) dz \nn\\
		&\le \Phi(1)t^{-d/\alpha}	 \int_{B_F(x,t^{1/\alpha}) \setminus B(w,\rho)}    J(w,z) dz  \le \frac{c_1}{\rho^{d+\alpha}} \Phi\bigg( \frac{(\rho \wedge A_0)^2}{(\delta_\Sigma(w) \wedge A_0)(	\delta_\Sigma(x,t) \wedge A_0)}\bigg).
	\end{align}

	For all $n\ge 1$ and $z\in B_F(x,2^nt^{1/\alpha}) \setminus B(x,2^{n-1}t^{1/\alpha})$, by  \eqref{e:r-r2-monotone}, we have
	\begin{align}\label{e:Meyer-integral-estimate-with-Phi-monotone}
		f(z) \le 2^{-(n-1)\lambda}t^{-d/\alpha}	\Phi\bigg( \frac{( 2^{n}t^{1/\alpha}  \wedge A_0)^2}{(\delta^\vee_\Sigma(x,z,t) \wedge 2^{n}t^{1/\alpha} \wedge A_0) (t^{1/\alpha}\wedge A_0)}\bigg)   J(w,z).
	\end{align}
		For all $n\ge 1$ with $2^n t^{1/\alpha}<\delta_\Sigma(x,t)$,
	 using \eqref{e:Meyer-integral-estimate-with-Phi-monotone} in the first inequality below, and \eqref{e:Phi-scaling}, \eqref{e:without-A0} and Lemma \ref{l:tail-estimate-2}(i) in the second, we obtain
		\begin{align}\label{e:Meyer-integral-estimate-with-Phi-2-case1}
		&	 \int_{B_F(x,2^nt^{1/\alpha}) \setminus (B(x,2^{n-1}t^{1/\alpha}) \cup B(w,\rho))} f(z) dz\nn\\
		&\le  2^{-(n-1)\lambda}t^{-d/\alpha}	\Phi\bigg( \frac{ 2^{n}t^{1/\alpha}  \wedge A_0}{t^{1/\alpha}\wedge A_0}\bigg)  \int_{B_F(x,2^nt^{1/\alpha}) \setminus B(w,\rho)}     J(w,z) dz\nn\\
		&\le  \frac{c_22^{n(d+\beta_1-\lambda )} }{\rho^{d+\alpha}} \Phi\left( \frac{(\rho \wedge A_0)^2}{(\delta_\Sigma(w) \wedge A_0) (\delta_\Sigma(x,t) \wedge A_0) }\right).
	\end{align} 
On the other hand, 
	for all $n\ge 1$ with $2^n t^{1/\alpha}\ge \delta_\Sigma(x,t)$, 
	 using 
	 \eqref{e:Meyer-integral-estimate-with-Phi-monotone} in the first inequality below,   Lemma \ref{l:tail-estimate-2}(ii) with $u=\delta_\Sigma(x,t)$, 
 $s=2^nt^{1/\alpha}$   and $A=2^{n}t^{1/\alpha} \wedge A_0$ in the second, and \eqref{e:without-A0} in the third, we obtain 
	\begin{align}\label{e:Meyer-integral-estimate-with-Phi-2-case2}
	&	 \int_{B_F(x,2^nt^{1/\alpha}) \setminus (B(x,2^{n-1}t^{1/\alpha}) \cup B(w,\rho))} f(z) dz \nn\\
			&\le  2^{-(n-1)\lambda}t^{-d/\alpha}	 \int_{B_F(x,2^nt^{1/\alpha}) \setminus B(w,\rho)}   \Phi\bigg( \frac{( 2^{n}t^{1/\alpha}  \wedge A_0)^2}{(\delta^\vee_\Sigma(x,z,t) \wedge 2^{n}t^{1/\alpha} \wedge A_0) (t^{1/\alpha}\wedge A_0)}\bigg)   J(w,z) dz\nn\\
				&\le   \frac{c_32^{n(d-\lambda)}}{\rho^{d+\alpha}} \bigg[ 1 +\bigg( \frac{(2^n t^{1/\alpha} \wedge A_0)^2}{(2^nt^{1/\alpha}  \wedge A_0)(t^{1/\alpha} \wedge A_0)}\bigg)^{\beta_1} \bigg] \Phi\left( \frac{(\rho \wedge A_0)^2}{(\delta_\Sigma(w) \wedge A_0) (\delta_\Sigma(x,t) \wedge A_0) }\right) \nn\\
		& \le \frac{c_42^{n(d+\beta_1-\lambda)}}{\rho^{d+\alpha}} \Phi\left( \frac{(\rho \wedge A_0)^2}{(\delta_\Sigma(w) \wedge A_0) (\delta_\Sigma(x,t) \wedge A_0) }\right).
	\end{align} 
	Set  $N:=\min\{n\ge 1: 2^{n}t^{1/\alpha} \ge r\}$.  Combining \eqref{e:Meyer-integral-estimate-with-Phi-1},
	\eqref{e:Meyer-integral-estimate-with-Phi-2-case1} and \eqref{e:Meyer-integral-estimate-with-Phi-2-case2}, we arrive at  
	\begin{align*}
		&	\int_{B_F(x,r) \setminus  B(w,\rho)} f(z)dz	\le \bigg( \int_{B_F(x,t^{1/\alpha}) \setminus  B(w,\rho)}   +  \sum_{n=1}^N  \int_{B_F(x,2^nt^{1/\alpha}) \setminus (B(x,2^{n-1}t^{1/\alpha}) \cup B(w,\rho))}  \bigg) 
f(z)dz\nn\\
		&\le \frac{c_1}{\rho^{d+\alpha}} \Phi\bigg( \frac{(\rho \wedge A_0)^2}{(\delta_\Sigma(w) \wedge A_0)(	\delta_\Sigma(x,t) \wedge A_0)}\bigg) +  \sum_{n=1}^N \frac{c_52^{n(d+\beta_1-\lambda)}}{\rho^{d+\alpha}} \Phi\left( \frac{(\rho \wedge A_0)^2}{(\delta_\Sigma(w) \wedge A_0) (\delta_\Sigma(x,t) \wedge A_0) }\right)\nn\\
		&\le   \frac{c_62^{(N+1)(d+\beta_1-\lambda)}}{\rho^{d+\alpha}} \Phi\left( \frac{(\rho \wedge A_0)^2}{(\delta_\Sigma(w) \wedge A_0) (\delta_\Sigma(x,t) \wedge A_0) }\right) \\
		&\le \frac{c_6 (4r/ t^{1/\alpha})^{d+\beta_1-\lambda}}{\rho^{d+\alpha}} \Phi\left( \frac{(\rho \wedge A_0)^2}{(\delta_\Sigma(w) \wedge A_0) (\delta_\Sigma(x,t) \wedge A_0) }\right).
	\end{align*}
	\qed

	\begin{lem}\label{l:Meyer-integral-estimate-off-diagonal}
For any $\eps\in (0,1)$,	there exists $C=C(\eps)>0$ such that for all 
$x\in F$, $\rho \in (0,\overline R)$, 
$a>0$,   $t\in (0, (\rho \wedge R_0)^\alpha]$ and $r\ge t^{1/\alpha}$  satisfying $\rho \le t^{\eps/\alpha} r^{1-\eps}$,
		\begin{align*}
			\int_{F\setminus B(x,r)} p^{(\rho)}(t,x,z) \Phi\bigg( \frac{a}{\delta_\Sigma(z) \wedge A_0} \bigg) dz	&  \le C\Phi\left( \frac{a}{ (\delta_\Sigma(x) \vee r) \wedge A_0 }\right).
		\end{align*}
	\end{lem}
	\pf Applying Corollary \ref{c:truncated-off-diagonal} with $\lambda= d+1$ and using Lemma \ref{l:Aikawa-integrals}(ii), we obtain
	\begin{align*}
		&\int_{F\setminus B(x,r)} p^{(\rho)}(t,x,z) \Phi\bigg( \frac{a}{\delta_\Sigma(z) \wedge A_0} \bigg) dz \le c_1\int_{F\setminus B(x,r)} t^{-d/\alpha} \bigg(  \frac{t^{1/\alpha}}{|x-z|}\bigg)^{d+1} \Phi\bigg( \frac{a}{\delta_\Sigma(z) \wedge A_0} \bigg) dz\\
		&\le \frac{ c_1t^{1/\alpha}}{r^{d+1}} \sum_{n=1}^\infty 2^{-(n-1)(d+1)} \int_{B_F(x,2^nr)\setminus B(x,2^{n-1}r)}  \Phi\bigg( \frac{a}{\delta_\Sigma(z) \wedge A_0} \bigg) dz\\
			&\le \frac{ c_2t^{1/\alpha}}{r}\Phi\bigg( \frac{a}{(\delta_\Sigma(x) \vee r) \wedge A_0} \bigg)  \sum_{n=1}^\infty 2^{-n}  = \frac{ c_2t^{1/\alpha}}{r}\Phi\bigg( \frac{a}{(\delta_\Sigma(x) \vee r) \wedge A_0} \bigg)  . 
	\end{align*}
	Since $t^{1/\alpha}\le r$, the result follows. \qed

	\begin{lem}\label{l:UHK-onesided-Phi}
		If there exist  $\lambda\in [d,d+\beta_1)$ and $C>0$ such that for all $\rho \in (0,\overline R]$, $x,y\in {F}$ and $t\in (0,(\rho \wedge R_0)^\alpha)$,
		\begin{align}\label{e:UHK-onesided-Phi-ass}
			p^{(\rho)}(t,x,y) \le  C t^{-d/\alpha} \bigg(1 \wedge \frac{t^{1/\alpha}}{|x-y|}\bigg)^{\lambda } \Phi\bigg( \frac{(|x-y| \wedge  A_0)^2}{(\delta_\Sigma^\vee(x,y,t)  \wedge |x-y| \wedge  A_0) (t^{1/\alpha}\wedge A_0)}\bigg),
		\end{align}
		then for any $\eps \in (0,\alpha-\beta_1)$,  there exists $C'=C'(\lambda, C, \eps )>0$ such that for all $\rho \in (0,\overline R]$, $x,y\in {F}$ and $t\in (0,(\rho \wedge R_0)^\alpha)$,
		\begin{align}\label{e:UHK-onesided-Phi-result}
			p^{(\rho)}(t,x,y) \le C' t^{-d/\alpha} \bigg(1 \wedge \frac{t^{1/\alpha}}{|x-y|}\bigg)^{\lambda +\eps } \Phi\bigg( \frac{(|x-y| \wedge A_0)^2}{(\delta_\Sigma^\vee(x,y,t) \wedge |x-y| \wedge A_0) (t^{1/\alpha}\wedge A_0)}\bigg).
		\end{align}
	\end{lem}
	\pf   
	Let $\rho \in (0,\overline R]$, $x,y\in F$ and $t\in (0,(\rho \wedge R_0)^\alpha)$. By symmetry, without loss of generality, we   assume that $\delta_\Sigma(x) \ge \delta_\Sigma(y)$.	 Set $r:=|x-y|$, $\delta:=(\alpha-\beta_1-\eps)/(d+\alpha)$ and $\sigma:=t^{\delta/\alpha}r^{1-\delta}$.
	From the  first  paragraph of the proof of Lemma \ref{l:UHK-without-Phi},  
	we know that, to prove  \eqref{e:UHK-onesided-Phi-result}, 
	it suffices to show that for the case $t^{1/\alpha}<\sigma<\rho$,
	\begin{align*}
		 \int_0^t \int_F p^{(\sigma)}(s, x, z) \int_{B_F(z,\rho)\setminus B(z,\sigma)} J(z,w)  p^{(\rho)}(t-s,w,y) \,dw\, dz\, ds
	\end{align*}
	is bounded above by the right-hand side of \eqref{e:UHK-onesided-Phi-result}.

	Suppose $t^{1/\alpha}<\sigma<\rho$. If $\delta_\Sigma(x)\le t^{1/\alpha}$, then,   since $\lambda+\eps<d+\alpha$, 
	\eqref{e:UHK-onesided-Phi-result} follows from Proposition \ref{p:UHK-tt}. 
	Assume $\delta_\Sigma(x) \ge t^{1/\alpha}$. 
	 Set
	\begin{align*}
			I_1&:=	\int_0^t \int_{F\setminus B(x,r)} p^{(\sigma)}(s, x, z) \int_{B_F(z,\rho)\setminus B(z,\sigma)} J(z,w)  p^{(\rho)}(t-s,w,y) \,dw\, dz\, ds,\\
				I_2&:=	\int_0^t \int_{ B_F(x,r)} p^{(\sigma)}(s, x, z) \int_{B_F(z,\rho)\setminus B(z,\sigma)} J(z,w)  p^{(\rho)}(t-s,w,y) \,dw\, dz\, ds.
	\end{align*}	
	 Using symmetry and Lemma \ref{l:jump-density-monotonicty} in the first inequality, Lemma  \ref{l:Meyer-integral-estimate-off-diagonal} in the second, Lemma \ref{l:blow-up-integrated-improve-1} in the third, \eqref{e:Phi-scaling} and  \eqref{e:without-A0} in the fourth, and the monotonicity of $\Phi$ and    $\delta_\Sigma(x) \vee r = \delta_\Sigma(x,t) \vee r \ge \delta_\Sigma(x,t) \wedge r$    in the fifth, we obtain
	\begin{align*}
	I_1&\le \frac{c_1}{\sigma^{d+\alpha}}\int_0^t \int_F p^{(\rho)}(t-s,y,w) \nn\\
	&\qquad \times  \int_{F\setminus (B(x,r) \cup B(w, \sigma))} p^{(\sigma)}(s, x, z) \Phi \bigg( \frac{(\sigma \wedge A_0)^2}{(\delta_\Sigma(z) \wedge A_0)(\delta_\Sigma(w)\wedge A_0)}\bigg)   \,dz\, dw\, ds\nn\\
	&\le \frac{c_2}{\sigma^{d+\alpha}}\int_0^t \int_F p^{(\rho)}(t-s,y,w) \Phi \bigg( \frac{(\sigma \wedge A_0)^2}{((\delta_\Sigma(x) \vee r) \wedge A_0)(\delta_\Sigma(w)\wedge A_0)}\bigg)     dw\, ds\nn\\
		&\le \frac{c_3}{\sigma^{d+\alpha}}\int_0^t \Phi \bigg( \frac{(\sigma \wedge A_0)^2}{((\delta_\Sigma(x) \vee r) \wedge A_0)( (t-s)^{1/\alpha}\wedge A_0)}\bigg)      ds\nn\\
			&\le \frac{c_4t^{\beta_1/\alpha}}{\sigma^{d+\alpha}}\Phi \bigg( \frac{(\sigma \wedge A_0)^2}{((\delta_\Sigma(x)\vee r)\wedge A_0)( t^{1/\alpha}\wedge A_0)}\bigg)\int_0^t   (t-s)^{-\beta_1/\alpha}    ds \nn\\
				&= \frac{c_5t}{\sigma^{d+\alpha}}\Phi \bigg( \frac{(\sigma \wedge A_0)^2}{((\delta_\Sigma(x)\vee r)\wedge A_0)( t^{1/\alpha}\wedge A_0)}\bigg)\le \frac{c_5t}{\sigma^{d+\alpha}}\Phi \bigg( \frac{(r \wedge A_0)^2}{(\delta_\Sigma(x,t)\wedge r\wedge A_0)( t^{1/\alpha}\wedge A_0)}\bigg).
	\end{align*} 
	For $I_2$, using symmetry and \eqref{e:UHK-onesided-Phi-ass}  in the first inequality below,	Lemma \ref{l:Meyer-integral-estimate-with-Phi} in the second,	Lemma \ref{l:blow-up-integrated-improve-1}  in the third, \eqref{e:Phi-scaling},   \eqref{e:without-A0}, $\sigma<r$ and the fact that $\delta_\Sigma(x) = \delta_\Sigma(x,t) = \delta_\Sigma(x,s)$ for all $s\in (0,t)$ in the fourth, and  $\lambda-d-\beta_1\ge -\beta_1>-\alpha$ in the fifth, we get that 
	\begin{align*}
			I_2&\le c_6	\int_0^t \int_F  p^{(\rho)}(t-s,y,w)\int_{ B_F(x,r)\setminus B(w,\sigma)}  s^{-d/\alpha} \bigg(1 \wedge \frac{s^{1/\alpha}}{|x-z|}\bigg)^\lambda \\
			&\qquad\qquad\qquad\qquad\qquad\quad \times  \Phi\bigg( \frac{(|x-z| \wedge A_0)^2}{(\delta^\vee_\Sigma(x,z,s) \wedge |x-z| \wedge A_0) (s^{1/\alpha}\wedge A_0)}\bigg)   J(w,z)  \,dz \, dw\, ds\\
			&\le \frac{c_7 r^{d+\beta_1-\lambda}}{ \sigma^{d+\alpha}} 	\int_0^t  s^{(\lambda-d-\beta_1)/\alpha}\int_F  p^{(\rho)}(t-s,y,w)  \Phi\left( \frac{(\sigma \wedge A_0)^2}{(\delta_\Sigma(w) \wedge A_0) (\delta_\Sigma(x,s) \wedge A_0) }\right)  dw\,ds\\
				&\le \frac{c_8 r^{d+\beta_1-\lambda}}{ \sigma^{d+\alpha}} 	\int_0^t  s^{(\lambda-d-\beta_1)/\alpha} \Phi\left( \frac{(\sigma \wedge A_0)^2}{((t-s)^{1/\alpha} \wedge A_0) (\delta_\Sigma(x,s) \wedge A_0) }\right)  ds\\
					&\le  \frac{c_9 t^{\beta_1/\alpha} r^{d+\beta_1-\lambda}}{ \sigma^{d+\alpha}} \Phi\left( \frac{(r \wedge A_0)^2}{(t^{1/\alpha} \wedge A_0) (\delta_\Sigma(x,t) \wedge A_0) }\right) 	\int_0^t  s^{(\lambda-d-\beta_1)/\alpha}  (t-s)^{-\beta_1/\alpha}  ds\\
					&\le  \frac{c_{10}t^{(\lambda+\alpha-d-\beta_1)/\alpha} r^{d+\beta_1-\lambda}}{ \sigma^{d+\alpha}} \Phi\left( \frac{(r \wedge A_0)^2}{ (\delta_\Sigma(x,t) \wedge A_0)(t^{1/\alpha} \wedge A_0) }\right) .
	\end{align*}
Since 	 $t^{1/\alpha} \le r$ and $\lambda<d+\beta_1$, we have
\begin{align*}
	\frac{t}{\sigma^{d+\alpha}} + \frac{t^{(\lambda+\alpha-d-\beta_1)/\alpha} r^{d+\beta_1-\lambda}}{ \sigma^{d+\alpha}}=\frac{t^{(\beta_1+\eps)/\alpha}}{r^{d+\beta_1+\eps}}+\frac{t^{(\lambda+\eps-d)/\alpha}}{r^{\lambda+\eps}} \le \frac{2t^{(\lambda+\eps-d)/\alpha}}{r^{\lambda+\eps}}.
\end{align*} 
	Thus, combining the estimates for $I_1$ and $I_2$,
		we arrive at the desired result.\qed

	\begin{prop}\label{p:UHK-onesided-Phi-iterated}
		For any $\lambda \in [d,d+\alpha)$, there exists $C=C(\lambda)>0$ such that for all $\rho \in (0,\overline R]$, $x,y\in {F}$ and $t\in (0,(\rho \wedge R_0)^\alpha)$,
		\begin{align*}
			p^{(\rho)}(t,x,y) \le C t^{-d/\alpha} \bigg(1 \wedge \frac{t^{1/\alpha}}{|x-y|}\bigg)^{\lambda } \Phi\bigg( \frac{(|x-y| \wedge A_0)^2}{(\delta_\Sigma^\vee(x,y,t) \wedge |x-y| \wedge  A_0) (t^{1/\alpha}\wedge A_0)}\bigg).
		\end{align*}
	\end{prop}
	\pf 
	Applying Proposition \ref{p:UHK-tt} with $\lambda=d+\beta_1$, and using \eqref{e:Phi-scaling-monotone} and \eqref{e:without-A0}, we obtain  for all $\rho \in (0,\overline R]$, $x,y\in {F}$ and $t\in (0,(\rho \wedge R_0)^\alpha)$,
	\begin{align*}
		p^{(\rho)}(t,x,y)& \le c_1 t^{-d/\alpha} \bigg(1 \wedge \frac{t^{1/\alpha}}{|x-y|}\bigg)^{d+\beta_1} \Phi\bigg( \frac{(|x-y|\wedge A_0)^2}{(t^{1/\alpha}\wedge A_0)^2}\bigg)\\
		& \le c_2 t^{-d/\alpha} \bigg(1 \wedge \frac{t^{1/\alpha}}{|x-y|}\bigg)^{d+\beta_1} \bigg(1 + \frac{|x-y|}{t^{1/\alpha}}\bigg)^{\beta_1}\Phi\bigg( \frac{|x-y|\wedge A_0}{t^{1/\alpha}\wedge A_0}\bigg)\\
			& \le 2^{\beta_1}c_2 t^{-d/\alpha} \bigg(1 \wedge \frac{t^{1/\alpha}}{|x-y|}\bigg)^{d} \Phi\bigg( \frac{|x-y|\wedge A_0}{t^{1/\alpha}\wedge A_0}\bigg)\\
				& \le 2^{\beta_1}c_2 t^{-d/\alpha} \bigg(1 \wedge \frac{t^{1/\alpha}}{|x-y|}\bigg)^{d} \Phi\bigg( \frac{(|x-y|\wedge A_0)^2}{(\delta^\vee_\Sigma(x,y,t) \wedge |x-y| \wedge A_0)(t^{1/\alpha}\wedge A_0)}\bigg).
	\end{align*}
	 Here we used the fact that
	  $	1+ a \le 2(1\wedge (1/a))^{-1}$ for all $a>0$ in the third inequality.  Thus, \eqref{e:UHK-onesided-Phi-ass} holds for $\lambda=d$. Now, applying Lemma \ref{l:UHK-onesided-Phi} finitely many times, we get the  result. \qed

	 \begin{prop}\label{p:blow-up-integrated-improve-final}	There exists $C>0$ such that for all $a>0$, $\rho \in (0,\overline R]$,  $x\in {F}$ and $t\in (0, (\rho \wedge R_0)^\alpha)$,
	 	\begin{align*}
	 		\int_F p^{(\rho)}(t, x, z) \Phi \bigg( \frac{a}{\delta_\Sigma(z)\wedge A_0}\bigg)   dz\,\le C\Phi \bigg( \frac{a}{\delta_\Sigma(x,t)\wedge A_0}\bigg) .
	 	\end{align*}
	 \end{prop}
	 \pf If $\delta_\Sigma(x)\le t^{1/\alpha}$, then the result follows from Lemma \ref{l:blow-up-integrated-improve-1}. Suppose $\delta_\Sigma(x) >t^{1/\alpha}$. 
	Let $I:=\int_{F: \delta_\Sigma(z)<\delta_\Sigma(x)} p^{(\rho)}(t, x, z) \Phi \big( \frac{a}{\delta_\Sigma(z)\wedge A_0}\big) \, dz$.
	 As in the proof of Lemma \ref{l:blow-up-integrated-improve-1}, we have
	 \begin{align*}
	\int_F p^{(\rho)}(t, x, z) \Phi \bigg( \frac{a}{\delta_\Sigma(z)\wedge A_0}\bigg)   dz	\le \Phi \bigg( \frac{a}{\delta_\Sigma(x)\wedge A_0}\bigg)	 + I.
	 \end{align*}
	 Thus, it suffices to estimate $I$.  
	 Note that
	\begin{align*}
		I
		&=  \int_{B_F(x, \delta_\Sigma(x)/2): \,  \delta_\Sigma(z)<\delta_\Sigma(x)} p^{(\rho)}(t, x, z) \Phi \bigg( \frac{a}{\delta_\Sigma(z)\wedge A_0}\bigg) \, dz\\
		&\; +  \sum_{n=1}^\infty  \sum_{m=0}^\infty \int_{B_F(x, 2^{m}\delta_\Sigma(x)) \setminus B(x, 2^{m-1}\delta_\Sigma(x)) :  \delta_\Sigma(z)\in [2^{-n} \delta_\Sigma(x), 2^{1-n}\delta_\Sigma(x))} p^{(\rho)}(t, x, z) \Phi \bigg( \frac{a}{\delta_\Sigma(z)\wedge A_0}\bigg) \, dz \\
		&=:I_0 + \sum_{n=1}^\infty \sum_{m=0}^\infty I_{n,m}.
	\end{align*}
	 For all $z\in B_F(x,\delta_\Sigma(x)/2)$, we have $\delta_\Sigma(z)\ge \delta_\Sigma(x)/2$. Thus, by \eqref{e:Phi-scaling} and \eqref{e:without-A0}, we get that 
	 \begin{align*}
	 	I_0 \le c_1\Phi \bigg( \frac{a}{\delta_\Sigma(x)\wedge A_0}\bigg)  \int_{B_F(x, \delta_\Sigma(x)/2)} p^{(\rho)}(t, x, z) \, dz \le c_1\Phi \bigg( \frac{a}{\delta_\Sigma(x)\wedge A_0}\bigg).
	 \end{align*}
	 On the other hand, for all $n\ge 1$ and $m\ge 0$, using Proposition \ref{p:UHK-onesided-Phi-iterated} with $\lambda=d+\beta_1$ in the first inequality below, and \eqref{e:Phi-scaling}, \eqref{e:without-A0} and Proposition \ref{p:Aikawa} with $q=(\beta_1+\gamma)/2$ in the second, we obtain
	 \begin{align*}
	 	&I_{n,m}\le c_2t^{-d/\alpha} \bigg(\frac{t^{1/\alpha}}{2^{m-1}\delta_\Sigma(x)}\bigg)^{d+\beta_1}\Phi\bigg( \frac{(2^m \delta_\Sigma(x) \wedge A_0)^2}{(\delta_\Sigma(x) \wedge  A_0) (t^{1/\alpha}\wedge A_0)}\bigg) \Phi \bigg( \frac{a}{2^{-n} \delta_\Sigma(x)\wedge A_0}\bigg)\\
	 	&\qquad\quad \times \int_{B_F(x, 2^{m}\delta_\Sigma(x))  :  \delta_\Sigma(z)< 2^{1-n}\delta_\Sigma(x)}  dz \\
	 	&\le \frac{c_32^{(m+n)\beta_1}t^{\beta_1/\alpha}}{2^{m(d+\beta_1)}\delta_\Sigma(x)^{d+\beta_1}}  \Phi \bigg( \frac{a}{ \delta_\Sigma(x)\wedge A_0}\bigg)  \bigg( \frac{2^m \delta_\Sigma(x)}{t^{1/\alpha}}\bigg)^{\beta_1}  (2^m\delta_\Sigma(x))^{d-(\beta_1+\gamma)/2} (2^{1-n}\delta_\Sigma(x))^{(\beta_1+\gamma)/2} \\
	 	&= \frac{c_4}{2^{m(\gamma-\beta_1)/2}2^{n(\gamma-\beta_1)/2}}\Phi \bigg( \frac{a}{ \delta_\Sigma(x)\wedge A_0}\bigg). 
	 \end{align*}
	 Hence, since $\gamma>\beta_1$, we obtain $\sum_{n=1}^\infty \sum_{m=0}^\infty I_{n,m} \le c_5 \Phi \Big( \frac{a}{ \delta_\Sigma(x)\wedge A_0}\Big)$. The proof is complete.
	 \qed

	 \begin{lem}\label{l:EP}
	 	There exists $C>0$ such that for all  $\rho \in (0,\overline R)$,  $x\in {F}$, 
	 	$r>0$ and $t\in (0, (\rho \wedge R_0)^\alpha)$,
	 	\begin{align*}
	 		\int_{F\setminus B(x,r)} p^{(\rho)}(t,x,y) dy \le \frac{Ct^{(\alpha-\beta_1)/(2\alpha)}}{r^{(\alpha-\beta_1)/2}}.
	 	\end{align*}
	 \end{lem}
	 \pf Applying Proposition \ref{p:UHK-onesided-Phi-iterated} with $\lambda=d+(\alpha+\beta_1)/2$,  and using Lemma \ref{l:Aikawa-integrals}(ii), \eqref{e:Phi-scaling} and \eqref{e:without-A0},  since $\beta_1<\alpha$, we obtain
	 \begin{align*}
	 	&\int_{F\setminus B(x,r)} p^{(\rho)}(t,x,y) dy  \\
	 	&\le c_1t^{(\alpha+\beta_1)/(2\alpha)} \int_{F\setminus B(x,r)}  \frac{1}{|x-y|^{d+(\alpha+\beta_1)/2}} \Phi\bigg( \frac{(|x-y| \wedge A_0)^2}{(\delta_\Sigma(y) \wedge |x-y| \wedge  A_0) (t^{1/\alpha}\wedge A_0)}\bigg) dy \\
	 	&\le 	 \frac{c_1t^{(\alpha+\beta_1)/(2\alpha)}}{r^{d+(\alpha+\beta_1)/2}}\\
	 	&\quad \times \sum_{n=1}^\infty \int_{B_F(x,2^nr)\setminus B(x,2^{n-1}r)}   \frac{1}{2^{(n-1)(d+ (\alpha+\beta_1)/2)}}\Phi\bigg( \frac{(2^nr \wedge A_0)^2}{(\delta_\Sigma(y) \wedge 2^{n-1}r \wedge  A_0) (t^{1/\alpha}\wedge A_0)}\bigg)dy\\
	 	&\le 	  	 \frac{c_2t^{(\alpha+\beta_1)/(2\alpha)}}{r^{(\alpha+\beta_1)/2}} \sum_{n=1}^\infty   \frac{1}{2^{n(\alpha+\beta_1)/2}} \Phi\bigg( \frac{(2^nr \wedge A_0)^2}{((\delta_\Sigma(x) \vee 2^{n}r) \wedge 2^{n-1}r \wedge  A_0) (t^{1/\alpha}\wedge A_0)}\bigg)\\
	 	&= 	  	 \frac{c_2t^{(\alpha+\beta_1)/(2\alpha)}}{r^{(\alpha+\beta_1)/2}} \sum_{n=1}^\infty   \frac{1}{2^{n(\alpha+\beta_1)/2}} \Phi\bigg( \frac{(2^nr \wedge A_0)^2}{( 2^{n-1}r \wedge  A_0) (t^{1/\alpha}\wedge A_0)}\bigg)\\
	 	&\le	  	 \frac{c_3t^{(\alpha+\beta_1)/(2\alpha)}}{r^{(\alpha+\beta_1)/2}} \sum_{n=1}^\infty   \frac{1}{2^{n(\alpha+\beta_1)/2}} \bigg( \frac{2^nr}{t^{1/\alpha}}\bigg)^{\beta_1} = \frac{c_4t^{(\alpha-\beta_1)/(2\alpha)}}{r^{(\alpha-\beta_1)/2}}.
	 \end{align*}\qed
	 
	 \begin{prop}\label{p:truncated-improved}
	 	There exist  $C>0$ and $K_0>1$ such that for all  $\rho \in (0,\overline R)$,  $t\in (0, (\rho \wedge R_0)^\alpha)$ and  $x,y\in {F}$ with $|x-y|\ge K_0\rho$,
	 	\begin{align*}
	 		p^{(\rho)}(t,x,y)\le \frac{Ct}{\rho^{d+\alpha}}.
	 	\end{align*}
	 \end{prop}
	 \pf By Lemma \ref{l:EP}, there exists $k\ge 1$ such that   for all  $\rho \in (0,\overline R)$,  $x_0\in {F}$, $r\ge k\rho$ and $t\in (0, (\rho \wedge R_0)^\alpha)$,
	 \begin{align*}
	 	\sup_{x\in B_{F}(x_0,r/4)} \int_{F\setminus B(x_0,r)} p^{(\rho)}(t,x,y) dy  &\le 	\sup_{x\in F} \int_{F\setminus B(x,3r/4)} p^{(\rho)}(t,x,y) dy  	\\
	 	&\le  \frac{c_1t^{(\alpha-\beta_1)/(2\alpha)}}{r^{(\alpha-\beta_1)/2}} \le \frac{c_1\rho^{(\alpha-\beta_1)/2}}{(k\rho)^{(\alpha-\beta_1)/2}}\le\frac12.
	 \end{align*}
	 Hence, \eqref{e:truncated-self-improvement-reverse-ass-fixed-t} holds with $g(r,t):=( c_1t^{(\alpha-\beta_1)/(2\alpha)}r^{-(\alpha-\beta_1)/2})
	  \wedge (1/2)$.
	 Set $K_0:= 20k(1+ \frac{\alpha+d}{\alpha-\beta_1})$.	By 
	 Corollary \ref{c:truncated-self-improvement-reverse}, it follows that  for all $\rho \in (0,\overline R)$, $t\in (0,(\rho \wedge R_0)^\alpha)$,   and $x,y \in F$ with $|x-y|\ge K_0\rho$,
	 \begin{align*}
	 	p^{(\rho)} (t,x,y) \le  c_2 t^{-d/\alpha}\bigg( \frac{2c_1 (t/2)^{(\alpha-\beta_1)/(2\alpha)}}{(k\rho)^{(\alpha-\beta_1)/2}}\bigg)^{(K_0/(10k)) -2} = \frac{c_3t}{\rho^{d+\alpha}}.
	 \end{align*}
	 \qed

	 \begin{lem}\label{l:dominant-term}
			 		For any $t>0$ and $x,y \in F$ satisfying $|x-y|\le t^{1/\alpha}$, 
	 	\begin{align*}
	 		\frac{ t }{|x-y|^{d+\alpha}}\Phi\left( \frac{(|x-y| \wedge A_0)^2}{( \delta_\Sigma(x,t) \wedge A_0)(\delta_\Sigma(y,t)  \wedge A_0) }\right)\ge  t^{-d/\alpha}.
	 	\end{align*}
	Moreover, 	there exists $C>1$ such that for any $t>0$ and $x,y \in F$ satisfying $|x-y|> t^{1/\alpha}$, 
	 	\begin{align*}
	 		\frac{ t }{|x-y|^{d+\alpha}}\Phi\left( \frac{(|x-y| \wedge A_0)^2}{( \delta_\Sigma(x,t) \wedge A_0)(\delta_\Sigma(y,t)  \wedge A_0)}\right)\le Ct^{-d/\alpha}.
	 	\end{align*}
	 \end{lem}
	 \pf If $|x-y|\le t^{1/\alpha}$, then using $\Phi\ge 1$, we get
	 \begin{align*}
	 	\frac{ t }{|x-y|^{d+\alpha}}\Phi\left( \frac{(|x-y| \wedge A_0)^2}{( \delta_\Sigma(x,t) \wedge A_0)(\delta_\Sigma(y,t)  \wedge A_0) }\right)\ge \frac{t}{|x-y|^{d+\alpha}}\ge  t^{-d/\alpha}.
	 \end{align*}  If $|x-y|>t^{1/\alpha}$, then using Lemma \ref{l:jump-density-monotonicty}, we obtain
	 \begin{align*}
	 	&	\frac{ t }{|x-y|^{d+\alpha}}\Phi\left( \frac{(|x-y| \wedge A_0)^2}{( \delta_\Sigma(x,t) \wedge A_0)(\delta_\Sigma(y,t)  \wedge A_0) }\right)\\
	 	&\le 	\frac{ c_1t }{t^{(d+\alpha)/\alpha}}\Phi\bigg( \frac{(t^{1/\alpha} \wedge A_0)^2}{( \delta_\Sigma(x,t) \wedge A_0)(\delta_\Sigma(y,t)  \wedge A_0)}\bigg)\le c_1\Phi(1) t^{-d/\alpha}.
	 \end{align*}
	 \qed

	Define
	$$
\wt p(t,x,y) :=	t^{-d/\alpha} \wedge  \left( \frac{ t }{|x-y|^{d+\alpha}}\Phi\left( \frac{(|x-y| \wedge A_0)^2}{(\delta_\Sigma(x,t)  \wedge A_0)(\delta_\Sigma(y,t)  \wedge A_0) }\right) \right). 
$$
	\begin{prop}\label{p:main-HKE}
		There exists $C>0$ such that $p(t,x,y)  \le  C\wt p(t,x,y)$ for all $(t,x,y) \in (0,R_0^\alpha) \times F\times F$.
	\end{prop}
	\pf	 
	 Let $(t,x,y) \in (0,R_0^\alpha)\times F \times F$ and set $r:=|x-y|$.  
	 If $r\le K_0t^{1/\alpha}$, where $K_0\ge 1$ is the constant in Proposition \ref{p:truncated-improved},  then the result follows from Proposition \ref{p:NDU}. Suppose $r>K_0t^{1/\alpha}$.
	 By Lemma \ref{l:dominant-term}, to get the desired result, it suffices to show that
	 \begin{align}\label{e:HKE-upper-bound-1}
	 	p(t,x,y) \le \frac{ c_1t }{r^{d+\alpha}}\Phi\left( \frac{(r \wedge A_0)^2}{(\delta_\Sigma(x,t) \wedge A_0)(\delta_\Sigma(y,t) \wedge A_0) }\right).
	 \end{align}
	  By Proposition \ref{p:truncated-improved}, since $\Phi\ge 1$, we get \begin{align}\label{e:HKE-upper-bound-2} p^{(r/K_0)}(t,x,y) \le \frac{c_2t}{r^{d+\alpha}} \le \frac{ c_2t }{r^{d+\alpha}}\Phi\left( \frac{(r \wedge A_0)^2}{(\delta_\Sigma(x,t) \wedge A_0)(\delta_\Sigma(y,t) \wedge A_0) }\right).
	  	\end{align} 
	  On the other hand, using Lemma \ref{l:jump-density-monotonicty} and symmetry in the first inequality below, and Proposition \ref{p:blow-up-integrated-improve-final} twice in the second, \eqref{e:Phi-scaling} and \eqref{e:without-A0} in the third, and $\beta_1<\alpha$ in the equality, we get
	  \begin{align}\label{e:HKE-upper-bound-3}
	  	& \int_0^t \int_F p^{(r/K_0)}(s, x, z) \int_{F\setminus B(z,r/K_0)} J(z,w)  p(t-s,w,y) \,dw\, dz\, ds \nn\\
	  	 &\le \frac{c_3}{r^{d+\alpha}} \int_0^t \int_F p^{(r/K_0)}(s, x, z) \int_{F\setminus B(z,r/K_0)}   p(t-s,y,w) \Phi\bigg( \frac{((r/K_0) \wedge A_0)^2}{(\delta_\Sigma(z) \wedge A_0) (\delta_\Sigma(w) \wedge A_0)}\bigg) \,dw\, dz\, ds \nn\\
	  	  &\le \frac{c_4}{r^{d+\alpha}} \int_0^t  \Phi\bigg( \frac{((r/K_0) \wedge A_0)^2}{(\delta_\Sigma(x,s) \wedge A_0) (\delta_\Sigma(y, t-s) \wedge A_0)}\bigg)  ds\nn\\
	  	   &\le \frac{c_5t^{2\beta_1/\alpha}}{r^{d+\alpha}}\Phi\bigg( \frac{(r \wedge A_0)^2}{(\delta_\Sigma(x,t) \wedge A_0) (\delta_\Sigma(y, t) \wedge A_0)}\bigg) \int_0^t  s^{-\beta_1/\alpha}(t-s)^{-\beta_1/\alpha}  ds\nn\\
	  	   &= \frac{c_6t}{r^{d+\alpha}}\Phi\bigg( \frac{(r \wedge A_0)^2}{(\delta_\Sigma(x,t) \wedge A_0) (\delta_\Sigma(y, t) \wedge A_0)}\bigg).
	  \end{align}
	  Combining \eqref{e:HKE-upper-bound-2} and \eqref{e:HKE-upper-bound-3} with  Proposition \ref{p:comparison-Meyer},
	   we get \eqref{e:HKE-upper-bound-1}. 
	   The proof is complete. \qed

\subsection{Proof of Theorem \ref{t:main-HKE}}\label{ss:proof-main}
In this subsection, we complete the proof of Theorem \ref{t:main-HKE}.

\medskip

	 	 \noindent \textbf{Proof of Theorem \ref{t:main-HKE}.} 	
	Theorem \ref{t:main-HKE}  holds for $t\in (0, R_0^\alpha)$ by 
	Subsection \ref{ss:t-lhke} and Proposition \ref{p:main-HKE}. 
	Suppose $R_0 <\infty$ and $R_0^\alpha\le t <T$. Choose $n \in \N$ such that $t/n < R_0^\alpha$. 
	Since $	\wt p(t/n^2, z,w) \asymp \wt  p(t/n, z,w)\asymp \wt  p(t, z,w)$ for  all $z,w \in D$, 
	by what has been proven for $t\in (0, R_0^\alpha)$ and the semigroup property, 
	\begin{align*}
	p(t,x,y) &= \int_F \dots \int_F p(t/n, x, y_1)\cdots p(t/n, y_{n-1}, y) dy_1 \cdots dy_{n-1}\\
	&\asymp \int_F \dots \int_F\wt p(t/n, x, y_1)\cdots \wt p(t/n, y_{n-1}, y) dy_1 \cdots dy_{n-1}\\
		&\asymp \int_F \dots \int_F\wt p(t/n^2, x, y_1)\cdots \wt p(t/n^2, y_{n-1}, y) dy_1 \cdots dy_{n-1}\\
		&\asymp \int_F \dots \int_F p(t/n^2, x, y_1)\cdots p(t/n^2, y_{n-1}, y) dy_1 \cdots dy_{n-1}\\
		&=
		p(t/n,x,y) \asymp \wt{p}(t/n,x,y) \asymp \wt{p}(t,x,y). 
	\end{align*}
\qed	
		 
 \section{Examples}\label{s:appl}

	In this section, we present three 
	examples that fall  within our framework: 
	the process introduced in \cite{DRV17, Vondra} 
in connection 
	with the nonlocal Neumann problem,
	the trace of the $\alpha$-stable process,
and resurrected processes in the closed half-space, which are conservative versions of the processes in the open half-space studied in \cite{KSV-jfa25}.  
	Throughout this section, we 
	assume that $d\ge 2$ and  $\alpha \in (0,2)$.
	 Recall that the 
	 fractional Laplacian $\Delta^{\alpha/2}:=-(-\Delta)^{\alpha/2}$ on $\R^d$ can be defined as 
	\begin{align*}
	\Delta^{\alpha/2} u(x):= p.v. \int_{\R^d} (u(y)-u(x)) J^\alpha(x,y)dy:= \lim_{\eps \downarrow 0} \int_{\{y\in \R^d: |x-y|>\eps\}} (u(y)-u(x)) J^\alpha(x,y)dy,
	\end{align*}
	where
	\begin{align*}
		J^\alpha(x,y):=\frac{2^\alpha \Gamma((d+\alpha)/2)}{\pi^{d/2} |\Gamma(-\alpha/2)|} |x-y|^{-d-\alpha}, \quad x,y \in \R^d.
	\end{align*}
	
	\begin{defn}\label{d:Lip}
	\rm 	 (i)	 	 We say that an open set $D\subset \R^d$ is  \textit{Lipschitz} if there exist  $\wh R>0$ and $\Lambda_0>0$  such that   for any $Q \in \partial D$, there exist a Lipschitz function $\Psi:\R^{d-1}\to \R$ with $\Psi(\wt 0)=0$ and $|\Psi(\wt y)-\Psi(\wt z)| \le \Lambda_0|\wt y - \wt z|$ for all $\wt y, \wt z\in \R^{d-1}$, and an orthonormal coordinate system CS$_{Q}$ with origin at $Q$ such that
	 	\begin{align}\label{e:local-coordinate}
	 		B_D(Q, \wh R) = \left\{ y= (\wt y, y_d) \in B(0, \wh R) \text{ in CS$_{Q}$} : y_d>\Psi(\wt y) \right\}.
	 	\end{align}

	 	\noindent (ii)  For an increasing continuous function $\ell:[0,\infty) \to [0,\infty)$ with $\ell(0)=0$, we say that an open set $D\subset \R^d$ is  $C^{1,\ell}$ if there exist  $\wh R>0$ and $\Lambda>0$ such that for any $Q \in \partial D$, there exist a function $\Psi:\R^{d-1} \to \R$ with  $	\Psi(\wt 0)= |\nabla \Psi(\wt 0)|=0$ and $ |\nabla \Psi(\wt y)-\nabla\Psi(\wt z)| \le \Lambda \ell(|\wt y-\wt z|)$ for all $\wt y, \wt z \in \R^{d-1},$ 
	 	and an orthonormal coordinate system CS$_{Q}$ with origin at $Q$ such that \eqref{e:local-coordinate} holds.
	 	
	 	If $D$ is $C^{1,\ell}$ with $\ell(r)=r$, then $D$ is said to be $C^{1,1}$.  If $D$ is $C^{1,\ell}$ with $\ell$ satisfying the Dini condition $\int_0^1 u^{-1}\ell(u)du<\infty$, then $D$ is said to be $C^{1,\rm Dini}$.
	 \end{defn}

	 \begin{defn}\label{d:half-space-like}
	\rm 	We say that an open set $D\subset \R^d$ is \textit{half-space-like} if there are constants $a<b$ such that $\{(\wt x,x_d)\in \R^d: x_d>b\} \subset D \subset \{(\wt x,x_d)\in \R^d: x_d>a\}$.
	 \end{defn}
	 
	  \subsection{Process associated with the nonlocal Neumann problem}
	 Let $D\subset \R^d$ be a Lipschitz open set
	 and $F:=\overline D$. 
The elliptic boundary value problem
	\begin{align}\label{e:nonlocal-Neumann}
\begin{cases} 
		{\Delta^{\alpha/2} u=f} \quad  \text{ in $D$,}\\[1pt]
		{\sN^{\alpha/2}_D u =0}  \quad  \text{ in $\R^d \setminus 
			F$,} 
	\end{cases}    
	\end{align}
	 where the ``nonlocal normal derivative" $\sN^{\alpha/2}_D $ is defined by
	 \begin{align*}
	 	N^{\alpha/2}_D u (x) = \int_D (u(x)-u(y)) J^\alpha(x,y)dy,
	 \end{align*}
	 was introduced in \cite{DRV17} as a nonlocal analogue  of the classical Neumann problem. As noticed in \cite{Ab20}, \eqref{e:nonlocal-Neumann} can be reformulated as a problem in $D$ with 
	 the regional-type nonlocal operator
	 \begin{align*}
	 	\Delta^{\alpha/2,N}_D u(x):= p.v. \int_D (u(y)-u(x)) K^\alpha_D(x,y)dy,
	 \end{align*}
	 where
	 \begin{align*}
	 	K^\alpha_D(x,y):= J^\alpha(x,y) + \int_{D^c} \frac{J^\alpha(x,z)J^\alpha(z,y)}{ \int_D J^\alpha(z,w)dw } \, dz, \quad x,y \in D.
	 \end{align*}
	 We refer to \cite{Vondra} for a probabilistic interpretation of $\Delta_D^{\alpha/2,N}$. The associated  form is given by
	 \begin{align*}
	 	\sE^{\alpha/2,N}_D (u,v) &= \frac12 \int_{D\times D} (u(x)-u(y)) (v(x)-v(y)) K^\alpha_D(x,y)\,dxdy,\\
	 	\sF^{\alpha/2,N}_D&= \text{ $(\sE^{\alpha/2,N}_D)_1$-closure of ${\rm Lip}_c(F)$ in $L^2(F)$. }
	 \end{align*}
	 It was proved in \cite{Ab20} that if $D$ is a bounded $C^{1,1}$ open set, then $K_D^\alpha(x,y)$ blows up at the boundary. 
Assume that
	 
	 \medskip

	 \setlength{\leftskip}{0.18in}
	 \setlength{\rightskip}{0.18in}
	 
	\noindent  {\bf (B)} $D$ is Lipschitz, and is either bounded or half-space-like or the complement of a bounded set.

	 \setlength{\leftskip}{0in}
	 \setlength{\rightskip}{0in}
	 
	 \medskip
	 
	Under {\bf (B)}, $D$ is $\kappa$-fat with localization constant $R_0=\diam(D)$ and 
	dim$_A(\partial D)=d-1$, 
 see Remark \ref{r:some-examples}(ii).  
We also have the following estimates on $K_D^\alpha(x,y)$.
	 
	 \begin{prop}\label{p:kernelestimate-Neumann}
	 	If {\bf (B)} holds, then for all $x,y \in D$,
	 	\begin{align*}
	 		K_D^\alpha(x,y) \asymp \frac{1}{|x-y|^{d+\alpha}}  \log	\bigg(e+ \frac{(|x-y| \wedge \diam(D^c))^2}{(\delta_{\partial D}(x)\wedge \diam(D^c))(\delta_{\partial D}(y)\wedge \diam(D^c))}\bigg).
	 	\end{align*}
	 \end{prop}
	 
	 The proof of Proposition \ref{p:kernelestimate-Neumann} will be given in 
	 Subsection \ref{ss:7.1}.
	 
	 By Proposition \ref{p:kernelestimate-Neumann}, 
	 under {\bf (B)},  
	 the form $(	\sE^{\alpha/2,N}_D , 	\sF^{\alpha/2,N}_D)$ satisfies condition {\bf (A)} with $\Phi(r) = \log(e+r)$ and $A_0=\diam(D^c)$. Thus, applying Theorem \ref{t:main-HKE}, 
we get sharp two-sided estimates of the heat kernel $p(t, x, y)$ of $(\sE^{\alpha/2,N}_D , 	\sF^{\alpha/2,N}_D)$ for $t< \diam(D)^\alpha$. 
When $\diam(D)<\infty$, since  $(	\sE^{\alpha/2,N}_D , 	\sF^{\alpha/2,N}_D)$  is conservative,  it is easy to show that $p(t,x,y) \asymp 1$ for all $(t,x,y) \in [\diam(D)^\alpha,\infty) \times F\times F$.
 
	 \subsection{Trace of $\alpha$-stable processes}

	 Let $Y=(Y_t)_{t\ge 0}$ be an isotropic $\alpha$-stable process in $\R^d$ with generator $\Delta^{\alpha/2}$, and let
	  $D\subset \R^d$ be a
	  Lipschitz   open set  and $F:=\overline D$. Define 

	  $$A_t:=\int_0^t \1_{\{Y_s \in F\}}ds,\quad t\ge 0,$$ 
	  and let $\tau_t:=\inf\{s>0:A_s>t\}$ be its right-continuous inverse. Define a process $X=(X_t)_{t\ge 0}$ by $X_t:=Y_{\tau_t}$. Then $X$ is a 
	    symmetric Hunt process on $F$, called the
	  \textit{trace process}  of $Y$ on $F$.  
	  Let 
	  $g^\alpha_{\overline{D}^c}(z,w)$ 
	  be the Green function of the isotropic $\alpha$-stable process on 
	  $\overline{D}^c$. 
	  The Poisson kernel 
	  $P_{\overline{D}^c}^\alpha(z,x)$ of $\overline{D}^c$ 
	  for the isotropic $\alpha$-stable process
	  is defined by
	  \begin{align*}
	  	P^\alpha_{\overline{D}^c}(z,x) = \int_{D^c} g^\alpha_{\overline{D}^c}(z,w) 
	  	J^\alpha(w,x)\,dw, \quad z\in D^c, \, x\in D.
	  \end{align*}
	  Since $D$ is Lipschitz, by 
	  \cite[Theorem 1]{Wu02} and the Ikeda-Watanabe formula,  
	  for any non-negative bounded Borel function $f$ on $\overline D$ and any 
	  $z\in \overline{D}^c$, we have
	  \begin{align}\label{e:Ikeda-Watanabe}
	  	\E_z\left[ f(Y_{\sigma_{\overline D}}) :\sigma_{\overline D} <\infty  \right] & = 	\E_z\left[ f(Y_{\sigma_{\overline D}}) :\sigma_{\overline D} <\infty, \, Y_{\sigma_{\overline D}} \in D \right] = \int_D 
	  	P^\alpha_{\overline{D}^c}(z,x) f(x) \, dx,
	  \end{align}
	  where $\sigma_{\overline D}:=\inf\{t>0:Y_t \in \overline D\}$ is the first hitting time of $\overline D$. For $f:\overline D \to \R$, let 
	  $P^\alpha_{\overline{D}^c}f(z):=f(z)$ 
	  for $z\in \overline D$ and let 
	  $P^\alpha_{\overline{D}^c}f(z):=\int_D P^\alpha_{\overline{D}^c}(z,x)f(x) dx$ for $z\in\overline{D}^c$
	  if the integral is well-defined. 
By \eqref{e:Ikeda-Watanabe} and \cite[Corollary 5.6.1 and Theorem 5.7.6]{CF12}, 
	  the Dirichlet form of $X$ is given by 
	  \begin{align*}
	  	\wc \sE(f,f)
	  	  	&= \frac12 \int_{D\times D} (f(x)-f(y))^2 \left( J^\alpha(x,y) + N^\alpha_D(x,y) \right) dxdy	+ \int_D f(x)^2 V_D^\alpha(x) dx,
	  \end{align*}
	  where
	  \begin{align*}
	  	N^\alpha_D(x,y)&:=\int_{D^c} J^\alpha(x,z) 
	  	P^\alpha_{\overline{D}^c} (z,y) dz = \int_{D^c \times D^c} J^\alpha(x,z) g^\alpha_{\overline{D}^c}(z,w) 
	  	J^\alpha(w,y) \, dzdw,\\
	  	V^\alpha_D(x)&:= \int_{ D^c}  J^\alpha(x,z)   \P_z(\sigma_{\overline D} = \infty) dz.
	  	\end{align*}
	  		 Hence, the jump kernel  $	\wc J^\alpha_D(x,y)$ of  $\wc \sE$ is given by 	   
	  \begin{align*}
	  	\wc J^\alpha_D(x,y):= J^\alpha(x,y)+ N^\alpha_D(x,y).
	  	\end{align*}
	  	 See \cite[Chapter 5]{CF12} for more results on trace processes.   The subprocess   $X^{ D}$  of $X$ killed upon leaving  $D$, or equivalently, hitting $\partial D$, is referred to as the \textit{path-censored $\alpha$-stable process} on $D$; see \cite{KPW14}.

The next result gives two-sided estimates on  $\wc J_D^\alpha(x,y)$ under the assumption
	  
	  \medskip
	  
	  \setlength{\leftskip}{0.18in}
	  \setlength{\rightskip}{0.18in}
	  
\noindent 	  {\bf (C)} $D$ is $C^{1,\rm Dini}$,  and is either bounded or half-space-like or the complement of a bounded set.

	  \setlength{\leftskip}{0in}
	  \setlength{\rightskip}{0in}

	  \medskip

	 The assertion of the following proposition for bounded $C^{1,\rm Dini}$ open sets will not be used in this paper, but it is of independent interest. 
	 The assumption {\bf (C)} implies that 
	 $D$ is $\kappa$-fat with localization constant $R_0=\diam(D)$ and 
	 dim$_A(\partial D)=d-1$.
	 
	 \begin{prop}\label{p:kernelestimate-trace}
	 	If {\bf (C)} holds, then for all $x,y \in D$,
	 	\begin{align*}
	 	\wc J_D^\alpha(x,y) \asymp \frac{1}{|x-y|^{d+\alpha}}  \left[ 1+  \bigg( \frac{(|x-y|\wedge \diam(D^c))^2}{(\delta_{\partial D}(x)\wedge \diam(D^c))(\delta_{\partial D}(y)\wedge \diam(D^c))}\bigg)^{\alpha/2}\right].
	 	\end{align*}
	 \end{prop}
	 
	 The proof of Proposition \ref{p:kernelestimate-trace} will be given in 
	 Subsection \ref{ss:7.2}.

	 When $D$ is either a  half-space-like $C^{1,\rm Dini}$ open set or the complement of the closure of a bounded $C^{1,\rm Dini}$  open set, we have  $\P_z(\sigma_{\overline D} =\infty) =0$ for all $z\in D^c$. Consequently, 
	 the Dirichlet form of $X$ is of pure-jump type  and has no killing part.  
By Proposition \ref{p:kernelestimate-trace}, under {\bf (C)}, $\wc J_D^\alpha(x,y)$ satisfies condition {\bf (A)} with $\Phi(r) = 1+r^{\alpha/2}$ and $A_0=\diam(D^c)$.  
Therefore, if $D$ is either a  half-space-like $C^{1,\rm Dini}$ open set or an exterior $C^{1,\rm Dini}$ open set, applying Theorem \ref{t:main-HKE}, 
	 we get sharp two-sided estimates of the heat kernel $p(t, x, y)$ of $X$ for  
	 $t\in (0, \infty)$.

	We mention in passing  that,   when $D$ is a bounded $C^{1,\rm Dini}$ open set, the killing measure $V_D^\alpha(x)dx$ of $X$  is non-trivial. Consequently, the Dirichlet form associated with $X$ has a non-trivial killing part and does not fall into the framework of this paper. 	
Nevertheless,	the resurrected process $X^{\rm res}$ associated with the form $\wc \sE^{\rm res}(f,f):= \wc \sE (f,f)  - \int_D f(x)^2 V_D^\alpha(x)dx$ is covered by our framework.  Theorem \ref{t:main-HKE} together with Proposition \ref{p:kernelestimate-trace} therefore yields sharp two-sided estimates of the heat kernel $p^{\rm res}(t, x, y)$ of $X^{\rm res}$ for  
$t\in (0, T]$ for every $T>0$. 
Moreover, when $D$ is a bounded $C^{1,\rm Dini}$ open set, using the boundary Harnack principle for $Y$ (see, e.g., \cite[Corollary 1.9]{CS26}), one can obtain that  $\P_z(\sigma_{\overline D} = \infty) \asymp (\delta_{\partial D}(z) \wedge 1)^{\alpha/2}$ for $z\in \overline D^c$.  Using this, one can further show that $V_D^\alpha(x)\asymp \delta_{\partial D}^{-\alpha/2}(x)$ for $x\in D$.  It follows that 
$V_D^\alpha(x)dx$ belongs to the Kato class of $ X^{\rm res}$, and hence the trace process $X$ can be viewed as a Feynman-Kac perturbation of $X^{\rm res}$	 by a Kato class potential. Consequently, short-time sharp two-sided heat kernel estimates for $X$ follow  from standard perturbation arguments. We do not pursue this direction here.

\subsection{Resurrected processes in the closed half-space}\label{ss:res-proc}
We now revisit the construction  
of the resurrected process in the open upper half-space $\R^d_+=\{x=(\wt{x}, x_d): \wt{x}\in \R^{d-1}, x_d>0\}$
given in \cite[Section 2]{KSV-jfa25}, where all details can be found.
This construction yields a jump kernel that blows up at the boundary $\partial \R^d_+$. 
In this subsection, we use this  jump kernel to obtain the conservative process $X$ that fits in our framework.

Consider the $\alpha$-stable process starting from a point in the upper half-space $\R^d_+$. 
Upon exiting the half-space and landing at a point  
$z\in \R^d_-=\{x=(\wt{x}, x_d): \wt{x}\in \R^{d-1}, x_d<0\}$
(the lower half-space), the process is immediately returned to a point $y\in \R^d_+$ according to the probability distribution $p(z,y)dy$. For $x,y\in \R^d_+$ let
$$
q(x,y):=\int_{\R^d_-}|x-z|^{-d-\alpha}p(z,y)\, dz
$$
denote the density of the return point in $\R^d_+$, provided the process starts from $x$. The kernel $q(x,y)$, called the \emph{resurrection kernel}, introduces additional jumps from $x$ to $y$. The resurrected process has the jump kernel  
\begin{equation}\label{e:jump-kernel-res}
 J(x,y)=|x-y|^{-d-\alpha} +q(x,y), \qquad x,y\in \R^d_+.
\end{equation}
Depending on the value of $\alpha\in (0,2)$, this process may be conservative or killed upon  hitting the boundary
$\partial \R^d_+$ of $\R^d_+$.

Let $\Psi$ be a positive function on 
 $[1, \infty)$  
satisfying  the following weak scaling condition: 
there exist  constants $c_1, c_2>0$ and $-\infty<\gamma_1\le \gamma_2< 1\wedge\alpha$ such that
$$
c_1\bigg(\frac{R}{r}\bigg)^{\gamma_1}\leq\frac{\Psi(R)}{\Psi(r)}\leq c_2 \bigg(\frac{R}{r}\bigg)^{\gamma_2},\quad 
1\le r<R <\infty.
$$
Assume that $p(z,y)$, $z\in \R^d_-$, $y\in \R^d_+$,  is given by the normalized version of
$$
\wt p(z, y)
:={|z_d|^{\alpha}}
\Psi\left(\frac{|y-z|^2}{y_d|z_d|}\right)
|y-z|^{-d-\alpha}.
$$
It is shown in \cite[Proposition 2.2]{KSV-jfa25} that 
the resurrection kernel is symmetric for this kind of $p$.
We note in passing that when  $\Psi(r)=1\vee r^{\alpha/2}$,  the resurrected process coincides with the trace process 
(killed upon hitting $\partial \R^d_+$ when $\alpha\in (1,2)$),
while when $\Psi(r)=1$, we recover the process associated with nonlocal Neumann problem. 
 Define  $\Psi_1(u):=1$ for $0\le u<2$ and 
$$
 \Psi_1(u):=\left( \int_1^2 \frac{\Psi(v)}{v}\, dv \right)^{-1}\int_1^u \frac{\Psi(v)}{v}\, dv  \quad\text{for $u\ge 2$}.
$$  
It is shown in \cite[Theorem 2.4]{KSV-jfa25} that the jump kernel $J(x,y)$ satisfies assumption \As\  with the 
blow-up weight function $\Phi=\Psi_1$ and $A_0=\infty$.

Consider the bilinear form 
$$
\EE(u,v)=\frac{1}{2}\int_{\R^d_+\times \R^d_+} (u(x)-u(y))(v(x)-v(y))J(x,y)\, dx dy,
$$
where the jump kernel $J(x,y)$ is given by \eqref{e:jump-kernel-res}. Let $\FF$ be the $\EE_1$ closure of 
${\rm Lip}_c(\overline{\R^d_+})$ in $L^2(\overline{\R^d_+})$ and denote the corresponding conservative process in $\overline{\R^d_+}$ by $X$. Applying Theorem \ref{t:main-HKE} and Remark \ref{r:th}(iii), we get global sharp  two-sided estimates for the heat kernel of $X$.
	 
\section{Proofs of Propositions \ref{p:kernelestimate-Neumann} and \ref{p:kernelestimate-trace}}\label{s:proof-examples} 
Let $d\ge 2$, $\alpha \in (0,2)$ and let $D\subset \R^d$ be a Lipschitz open set.  Set $R_0:=\diam(D)$ and $A_0:=\diam(D^c)$. 

	We will use the following simple fact several times:
	 \begin{equation}\label{e:dist-D-Dc}
	 	|x-w| \ge \delta_{\partial D}(x) \vee \delta_{\partial D}(w) \quad \text{for all} \;\, x \in D \text{ and } w \in D^c.
	 \end{equation}
	\begin{lem}\label{l:dist-D}
		(i) If $D$ is bounded, then $|y-z| \le \delta_{\partial D}(z) +R_0$ for all $y \in D$ and $z\in D^c$.
\noindent (ii)	If $D^c$ is bounded, then $|x-y| \le 2\delta_{\partial D}(y) + A_0$ for all $x,y \in D$ with $\delta_{\partial D}(y) \ge \delta_{\partial D}(x)$.
 	\end{lem}  
	 \pf (i) 
	 Let $Q_z\in \partial D$ be such that $\delta_{\partial D}(z) = |z-Q_z|$. Then we have 
	 $|y-z| \le |y-Q_z| + |Q_z-z| \le R_0 + \delta_{\partial D}(z)$.
	 
	 \noindent (ii)  Let $Q_x, Q_y\in \partial D$ be such that $\delta_{\partial D}(x) = |x-Q_x|$ and $\delta_{\partial D}(y) = |y-Q_y|$. Then we have $	|x-y| \le \delta_{\partial D}(x) + |Q_x-Q_y| + \delta_{\partial D}(y) \le 2\delta_{\partial D}(y)+A_0.$ \qed 
	 
For $Q\in \partial D$, $\theta>0$ and $r\in (0,\infty]$, define
\begin{align*}
	\sC_Q(r,\theta):=\{z=(\wt z, z_d) \text{ \rm in  CS$_{Q}$}: 
	-r<z_d<-\theta |\wt z|\}, 
\end{align*}
where CS$_Q$ is the coordinate system in Definition \ref{d:Lip}(i). Since $D$ is Lipschitz, there exist constants $r_0>0$  and $\theta>0$ such that
\begin{align}\label{e:exterior-cone}
	\sC_Q(r_0,\theta) \subset D^c \quad \text{for all $Q\in \partial D$.} 
\end{align}
Moreover, if $D$ is half-space-like, then there exists  $\theta>0$ such that \eqref{e:exterior-cone} holds with $r_0=\infty$.

For the remainder of this section, we let $r_0$ and $\theta$ be the constants in \eqref{e:exterior-cone}.

	 \begin{lem}\label{l:reflected-point}
	 	Suppose that $D$ is either bounded or half-space-like.
	 Then	 there exists   $\eps \in (0,1)$ depending only on $D$  such that for 
		 any  $y\in D$,  there exists $y^* \in D^c$ satisfying
	 	\begin{align}\label{e:reflect-point}
	 				\delta_{\partial D}(y^*)  \ge \eps \delta_{\partial D}(y) \quad \text{and} \quad |y-y^*| \le 2\delta_{\partial D}(y).
	 	\end{align} 
	 Moreover,	for any $y \in D$ and  
	  $w \in B(y^*, \delta_{\partial D}(y^*)/2)$, 
	 it holds that
	 	\begin{gather}
	 		\delta_{\partial D}(w)  \ge \eps \delta_{\partial D}(y)/2,\quad\;\; 
	 		|w-y| \le 3\delta_{\partial D}(y),\label{e:reflect-point1}
	 	\end{gather}
	 	and 
	 	\begin{align}\label{e:reflect-point2}
	 			|z-w| \le 4|y-z| \quad \text{for all $z\in D^c$}.
	 	\end{align} 	
	 \end{lem}
	 \pf We first prove \eqref{e:reflect-point}.  Let $y\in D$ and $Q_y\in \partial D$ be such that $|y-Q_y|=\delta_{\partial D}(y)$. By \eqref{e:exterior-cone}, there exist constants  $r_0\in (0,\infty]$ and $\theta>0$ independent of $Q_y$ such that $\sC_{Q_y}(r_0,\theta) \subset D^c$. Further, we have $r_0=\infty$ when $D$ is half-space-like. 
	 Let $y^* \in D^c$ 
with coordinates $ y^* = (\wt 0 , -(\delta_{\partial D}(y) \wedge (r_0/2)))$ in CS$_{Q_y}$. 
	Observe that
	 $$
	 \delta_{\partial D}(y^*) \ge \delta_{\partial \sC_{Q_y}(r_0,\theta)}(y^*) \ge c_1 (\delta_{\partial D}(y) \wedge (r_0/2)) \quad \text{and} \quad  |y-y^*| \le |y-Q_y| + |y^*- Q_y| \le 2\delta_{\partial D}(y).
	 $$
	 Since $\delta_{\partial D}(y) \wedge (r_0/2) \ge (1 \wedge (r_0/(2R_0)))\delta_{\partial D}(y)$ if $D$ is bounded, and $r_0=\infty$ if $D$ is half-space-like, we obtain \eqref{e:reflect-point}.

	For \eqref{e:reflect-point1}, using   \eqref{e:reflect-point} and  \eqref{e:dist-D-Dc}, we obtain 
		$\delta_{\partial D}(w) \ge \delta_{\partial D}(y^*)/2 \ge \eps \delta_{\partial D}(y)/2$ and $	|w-y| \le |y-y^*| + \delta_{\partial D}(y^*)/2 \le (3/2)|y-y^*| \le 3\delta_{\partial D}(y)$. 
For \eqref{e:reflect-point2}, we get from   \eqref{e:reflect-point1} and  \eqref{e:dist-D-Dc} that $|z-w| \le |y-z| +|w-y|  \le |y-z| + 3\delta_{\partial D}(y) \le 4|y-z|.$
	 The proof is complete. \qed 
	 
	 	For $x,y\in D$, define
	 \begin{align*}
	 	\delta_{\partial D}^\wedge(x,y):=\delta_{\partial D}(x) \wedge \delta_{\partial D}(y) \quad \text{ and } \quad \delta_{\partial D}^\vee(x,y):=\delta_{\partial D}(x) \vee \delta_{\partial D}(y).
	 \end{align*}

	 For all $\gamma \in \R$ and $y \in D$, define a function $f_{\gamma,y}:D \to (0,\infty)$ by
	 \begin{align*}
	 	f_{\gamma,y} (x) = \begin{cases}
	 		(|x-y| + \delta_{\partial D}^\vee(x,y))^{\gamma} &\mbox{ if } \gamma>0,\\
	 		\log(e+|x-y|/\delta_{\partial D}(x))&\mbox{ if } \gamma=0,\\
	 		\delta_{\partial D}(x)^{\gamma}&\mbox{ if } \gamma<0.
	 	\end{cases} 
	 \end{align*}
	 Observe that  for all $\gamma \in \R$ and $y,x\in D$, 
	 \begin{align}\label{e:f-gamma-bound}
	 	f_{\gamma,y} (x) \ge (|x-y| + \delta_{\partial D}^\vee(x,y))^{\gamma}.
	 \end{align}

	 \begin{lem}\label{l:technical-1}
	 	Let $a_1\ge -d$, $a_2\ge -d$ and $a_3\ge 0$ be  constants such that  $a_3<d+a_1+a_2$. Then  there exists $C>0$ such that for all $x,y \in D$,
	 	\begin{align*}
	 		&	\int_{D^c} \frac{\delta_{\partial D}(z)^{a_3}}{|x-z|^{d+a_1}|z-y|^{d+a_2}}dz \le  C\left( \frac{ f_{a_3-a_1,y}(x)}{(|x-y| + \delta_{\partial D}^\vee(x,y))^{d+a_2} } +  \frac{ f_{a_3-a_2,x}(y)}{(|x-y| + \delta_{\partial D}^\vee(x,y))^{d+a_1} }  \right).	\end{align*}
	 \end{lem}
	 \pf  Let $x,y \in D$ and 
	 we use the notation $\delta^\vee:=\delta_{\partial D}(x)\vee \delta_{\partial D}(y)$ in this proof. Define
	 \begin{align*}
	 	E_1&:=\left\{z\in D^c: |x-z|\le |y-z|, \, \delta_{\partial D}(z) \le 2|x-y|\right\},\\
	 	E_2&:=\left\{z\in D^c: |x-z|> |y-z|, \, \delta_{\partial D}(z) \le 2|x-y|\right\},\\
	 	E_3&:=\left\{z\in D^c:  \delta_{\partial D}(z)> 2|x-y|\right\} = D^c \setminus (E_1\cup E_2).
	 \end{align*}
	 For all $z \in E_1$, we have $|x-y|\le |y-z|+|x-z|\le 2|y-z|$. Combining this with \eqref{e:dist-D-Dc}, we get
	  $
	 |y-z|= |y-z| \vee |x-z| \ge (|x-y|/2) \vee \delta^\vee \ge (|x-y| + \delta^\vee)/4$ for all $ z \in E_1.$
	 Using this in the first inequality below and \eqref{e:dist-D-Dc}  in the second, we obtain
	 \begin{align*}
	 &	\int_{E_1} \frac{\delta_{\partial D}(z)^{a_3}}{|x-z|^{d+a_1}|z-y|^{d+a_2}}dz\le \frac{4^{d+a_2}}{(|x-y| + \delta^\vee)^{d+a_2} } \int_{E_1} \frac{\delta_{\partial D}(z)^{a_3} }{|x-z|^{d+a_1}} dz\\
	 	&\le \frac{4^{d+a_2}}{(|x-y| + \delta^\vee)^{d+a_2} } \int_{B(x, \delta_{\partial D}(x) + 2|x-y|) \setminus B(x,\delta_{\partial D}(x))} \frac{dz}{|x-z|^{d+a_1-a_3}} \le  \frac{c_1 f_{a_3-a_1,y}(x)}{(|x-y| + \delta^\vee)^{d+a_2} }.
	 \end{align*} 
	Similarly, for all $z \in E_2$, we have  
	 $
	 |x-z| \ge (|x-y|/2) \vee \delta^\vee \ge (|x-y| + \delta^\vee)/4.
	 $
	 Using this in the first inequality below and  \eqref{e:dist-D-Dc} in the second, we get
	 \begin{align*}
	 	&\int_{E_2} \frac{\delta_{\partial D}(z)^{a_3}}{|x-z|^{d+a_1}|z-y|^{d+a_2}}dz\le \frac{4^{d+a_1}}{(|x-y| + \delta^\vee)^{d+a_1} } \int_{E_2} \frac{\delta_{\partial D}(z)^{a_3} }{|z-y|^{d+a_2}} dz\\
	 	&\le \frac{4^{d+a_1}}{(|x-y| + \delta^\vee)^{d+a_1} } \int_{
		B(y, \delta_{\partial D}(y) + 2|x-y|) 
		\setminus B(y,\delta_{\partial D}(y))} \frac{dz}{|y-z|^{d+a_2-a_3}} \le  \frac{c_2 f_{a_3-a_2,x}(y)}{(|x-y| + \delta^\vee)^{d+a_1} }.
	 \end{align*} 
	  Finally, for all $z\in E_3$, we have $	|x-z| \wedge |y-z| \ge \delta_{\partial D}(z) >2|x-y|$. Thus $\frac12 |x-z| \le |y-z| \le 2 |x-z|$
	  and, by \eqref{e:dist-D-Dc},  \begin{align}\label{e:technical-2}
	 	|x-z| \ge  2^{-1} (|x-z| \vee |y-z|)  \ge |x-y| \vee (\delta^\vee/2) \ge (|x-y|+\delta^\vee)/4 \quad \text{for all} \;\, z \in E_3.
	 \end{align}
	 Using $|x-z| \asymp |y-z|$ and \eqref{e:technical-2} in the first inequality below,  $d+a_1+a_2-a_3>0$ in the equality and  \eqref{e:f-gamma-bound} in the second inequality, we obtain
	 \begin{align*}
	 	&\int_{E_3} \frac{\delta_{\partial D}(z)^{a_3}}{|x-z|^{d+a_1}|z-y|^{d+a_2}}dz \\
	 	&\le c_3	\int_{B(x, (|x-y|+\delta^\vee)/4)^c} \frac{dz}{|x-z|^{2d+a_1+a_2-a_3}} = \frac{c_4}{(|x-y|+\delta^\vee )^{d+a_1+a_2-a_3}} \le  \frac{c_4 f_{a_3-a_1,y}(x)}{(|x-y| + \delta^\vee)^{d+a_2} }.
	 \end{align*}
	 The proof is complete. \qed

	 \subsection{Proof of Proposition \ref{p:kernelestimate-Neumann}} \label{ss:7.1}

Throughout this subsection,	 we assume that $D$ satisfies {\bf (B)}.	 By the argument leading to  \cite[Proposition 2.1]{AFR23}, we have that, for any $ R>0$, there are comparison constants depending on $R$  such that for all $Q \in \partial D$ and  $x,y \in B_D(Q, R)$,
	 \begin{align}\label{e:interaction-Neumann-bounded}
	 	K_D^\alpha(x,y)  \asymp  \frac{1}{|x-y|^{d+\alpha}}  \log	\bigg(e+ \frac{|x-y|}{\delta_{\partial D}^\wedge(x,y)}\bigg).
	 \end{align}

	 \begin{lem}\label{l:reformulation-1}
	 		 For all $x,y \in D$,
	 	\begin{equation}\label{e:reformulation-1}
	 		\log	\bigg(e+ \frac{|x-y|}{\delta^\wedge_{\partial D}(x,y)}\bigg) \asymp 	\log	\bigg(e+ \frac{|x-y|^2}{\delta_{\partial D}(x)\delta_{\partial D}(y)}\bigg).
	 	\end{equation} 
	 \end{lem} 
	 \pf Without loss of generality, by symmetry, it suffices to prove \eqref{e:reformulation-1}  for all $x,y \in D$ with $\delta_{\partial D}(x)\le \delta_{\partial D}(y)$. 
	 If $\delta_{\partial D}(y) \ge 2|x-y|$, then $\delta_{\partial D}(x)\ge \delta_{\partial D}(y) - |x-y| \ge \delta_{\partial D}(y)/2$ so that $\delta_{\partial D}(x) \asymp \delta_{\partial D}(y)$ and $\delta_{\partial D}(x) \ge |x-y|$. Hence, the both hand sides of \eqref{e:reformulation-1} are comparable with  $1$. If  $\delta_{\partial D}(y)<2|x-y|$, then we have
	 \begin{align*}
	 	\log	\bigg(e+ \frac{|x-y|}{\delta_{\partial D}(x)}\bigg) \le  	\log	\bigg(e+ \frac{2|x-y|^2}{\delta_{\partial D}(x)\delta_{\partial D}(y)}\bigg) \le 2	\log	\bigg(e+ \frac{|x-y|^2}{\delta_{\partial D}(x)\delta_{\partial D}(y)}\bigg)
	 \end{align*}
	 and 
	 \begin{align*}
	 	2\log	\bigg(e+ \frac{|x-y|}{\delta_{\partial D}(x)}\bigg) \ge  	\log	\bigg(e+ \frac{|x-y|^2}{\delta_{\partial D}(x)^2}\bigg) \ge 	\log	\bigg(e+ \frac{|x-y|^2}{\delta_{\partial D}(x)\delta_{\partial D}(y)}\bigg).
	 \end{align*} The proof is complete.\qed

	 \begin{lem}\label{l:Poisson-2}
	 If $D$ is  either  half-space-like or the complement of a bounded set, then   
	  for $z\in \overline{D}^c$,
	 	\begin{align*}
	 \int_D J^\alpha(z,x)dx \asymp \delta_{\partial D}(z)^{-\alpha}.
	 	\end{align*}  
	 \end{lem}
	 \pf 
	 For all $z \in \overline{D}^c$, we have $
	 \int_D J^\alpha(z,x) dx\le 	c_1 \int_{B(z, \delta_{\partial D}(z))^c} |z-x|^{-d-\alpha}dx  = c_2 \delta_{\partial D}(z)^{-\alpha}.$
	Applying Lemma \ref{l:reflected-point} to 
	$\overline{D}^c$, we get that for any $z\in \overline{D}^c$, 
	there exists 
	$z^* \in D$ satisfying $\delta_{\partial D}(z^*) \ge \eps \delta_{\partial D}(z)$ and $|z-z^*| \le 2\delta_{\partial D}(z)$ for some constant $\eps \in (0,1)$. 
	Hence, for all $z\in D^c$, we obtain
	 \begin{align*}
	 	\int_D J^\alpha(z,x) dx \ge c_3	
		\int_{B(z^*, \delta_{\partial D}(z^*))} \frac{dx}{|z-x|^{d+\alpha}} \ge  \frac{c_4\delta_{\partial D}(z^*)^d}{ (|z- z^*| + \delta_{\partial D}(z^*))^{d+\alpha}}
		\ge c_5 \delta_{\partial D}(z)^{-\alpha}.
	 \end{align*} \qed

	 \noindent \textbf{Proof of Proposition \ref{p:kernelestimate-Neumann}.}  
	   If $D$ is bounded, then  the result follows  from \eqref{e:interaction-Neumann-bounded} with  $R=2R_0$	together with Lemma \ref{l:reformulation-1}.

	Suppose that $D$ is either half-space-like or the complement of a bounded set. 
	Let $x,y \in D$. By symmetry, without loss of generality, we assume that $\delta_{\partial D}(x)\le \delta_{\partial D}(y)$. Let  $Q_x\in \partial D$ be such that $\delta_{\partial D}(x) = |x-Q_x|$. By Lemma \ref{l:Poisson-2}, we have
	 \begin{align}\label{e:interaction-Neumann-unbounded}
	 	K^\alpha_D(x,y)\asymp \frac{1}{|x-y|^{d+\alpha}}+ \int_{D^c} \frac{\delta_{\partial D}(z)^\alpha}{|x-z|^{d+\alpha} |z-y|^{d+\alpha}}dz.
	 \end{align} 
Hence,	 applying Lemma \ref{l:technical-1} with $a_1=a_2=a_3=\alpha$,  we obtain
	 \begin{align}\label{e:interaction-Neumann-unbounded-upper}
	 	K^\alpha_D(x,y) &\le \frac{c_1}{|x-y|^{d+\alpha}} + \frac{c_2}{(|x-y|+\delta_{\partial D}(y))^{d+\alpha}} \log \bigg(e + \frac{|x-y|}{\delta_{\partial D}(x)}\bigg).
	 \end{align}

	 (i) Assume first that $D$ is half-space-like.  By \eqref{e:interaction-Neumann-unbounded-upper} and Lemma \ref{l:reformulation-1}, we obtain
	 \begin{align*}
	 		K^\alpha_D(x,y)\le \frac{c_1+c_2}{|x-y|^{d+\alpha}} \log \bigg(e + \frac{|x-y|}{\delta_{\partial D}(x)}\bigg) \le  \frac{c_3}{|x-y|^{d+\alpha}} \log \bigg(e + \frac{|x-y|^2}{\delta_{\partial D}(x)\delta_{\partial D}(y)}\bigg).
	 \end{align*}
	 
For the lower bound, by Lemma \ref{l:reformulation-1} and \eqref{e:interaction-Neumann-unbounded}, it suffices to consider the case $|x-y|/\delta_{\partial D}(x) \ge e$.
Suppose $|x-y|/\delta_{\partial D}(x) \ge e$.  
By \eqref{e:exterior-cone},  there exists $\theta >0$ independent of $x$ such that $\sC_{Q_x}(\infty, \theta) \subset D^c$. Set $E:= \sC_{Q_x} (|x-y|,2\theta)$.
	 Then $E\subset \sC_{Q_x}(\infty, \theta) \subset  D^c$. Moreover, for all $z\in E$ 
	 with coordinates 
	 $z=(\wt z,z_d) $ in CS$_{Q_x}$, since $ 2\theta |\wt z|<-z_d<|x-y|$, we get that
 \begin{gather*}
	 	\delta_{\partial D}(z) \ge \delta_{\partial \sC_{Q_x}(\infty,\theta)}(z) = \frac{-z_d-\theta |\wt z|}{\sqrt{\theta^2+1}} \ge \frac{-z_d}{2\sqrt{\theta^2+1}},\\
	 		|x-z| \le |x| + |z| \le \delta_{\partial D}(x) + (-z_d) \sqrt{1+ (2\theta)^{-2}}, \\
	 	|z-y| 
	  	\le
	  	 |x-y|+\delta_{\partial D}(x) + (-z_d) \sqrt{1+ (2\theta)^{-2}}
	  	  \le (1+e^{-1}+ \sqrt{1+ (2\theta)^{-2}})|x-y|.
 \end{gather*}
	 Combining these inequalities with \eqref{e:interaction-Neumann-unbounded}, we arrive at
	 \begin{align*}
	 	&K^\alpha_D(x,y) \ge  c_4\int_{E} \frac{\delta_{\partial D}(z)^\alpha}{|x-z|^{d+\alpha} |z-y|^{d+\alpha}}dz \ge \frac{c_5}{|x-y|^{d+\alpha}}\int_{E} \frac{|z_d|^{\alpha}}{(\delta_{\partial D}(x) + |z_d|)^{d+\alpha}}dz\\
	 	&\ge \frac{c_6}{|x-y|^{d+\alpha}}\int_{0}^{|x-y|} \frac{s^{\alpha}}{(\delta_{\partial D}(x)+s)^{d+\alpha}} \int_0^{s/(2\theta)} \ell^{d-2} d\ell \, ds= \frac{c_7}{|x-y|^{d+\alpha}}\int_{0}^{|x-y|} \frac{s^{d-1+\alpha}}{(\delta_{\partial D}(x)+s)^{d+\alpha}} ds\\
	 	&\ge \frac{c_7}{|x-y|^{d+\alpha}}\int_{\delta_{\partial D}(x)}^{|x-y|} \frac{s^{d-1+\alpha}}{(2s)^{d+\alpha}} ds=\frac{c_8}{|x-y|^{d+\alpha}}\log \bigg( \frac{|x-y|}{\delta_{\partial D}(x)}\bigg) \ge \frac{c_8}{2|x-y|^{d+\alpha}} \log \bigg(e+ \frac{|x-y|}{\delta_{\partial D}(x)}\bigg) ,
	 \end{align*} 
	 where we used the fact that $\log r \ge 2^{-1}\log(e+r)$ for all $r \ge e$ in the last inequality.  Combining this with Lemma \ref{l:reformulation-1}, we get the desired lower bound, and the proof for half-space-like $D$ is complete. 
	 
	  (ii) Suppose that $D$ is the complement of a bounded set.  We consider the following three cases separately.

	 \smallskip
	 
	\noindent \textit{Case 1:}  $\delta_{\partial D}(y) < 2A_0$.  By Lemma \ref{l:dist-D}(ii), we have $|Q_x-y| \le  \delta_{\partial D}(x)+|x-y|  \le \delta_{\partial D}(x) +2\delta_{\partial D}(y)+A_0<
	7A_0$, and thus, $x,y \in B(Q_x, 7A_0)$. Applying \eqref{e:interaction-Neumann-bounded} with $Q=Q_x$ and $R=7A_0$, 
	and using Lemma \ref{l:reformulation-1}, we obtain the result in this case.
	 
	 \smallskip
	 
\noindent \textit{Case 2:} $\delta_{\partial D}(x) <A_0$ and $\delta_{\partial D}(y)\ge 2A_0$. 
In this case, by Lemma \ref{l:dist-D}(ii) we have   $5\delta_{\partial D}(y)/2 \ge  |x-y| \ge\delta_{\partial D}(y) - \delta_{\partial D}(x)\ge\delta_{\partial D}(y)/2\ge A_0 \ge \delta_{\partial D}(x)$ and,
 so by Lemma \ref{l:dist-D}(i) and \eqref{e:dist-D-Dc}, $|y-z| \asymp |y-x|\asymp \delta_{\partial D}(y)$ for all $z \in D^c$. 
  Hence, by \eqref{e:interaction-Neumann-unbounded}, it holds that
	 \begin{align}\label{e:interaction-Neumann-unbounded-case2-1}
	 	K^\alpha_D(x,y)&\asymp \frac{1}{|x-y|^{d+\alpha}} \left( 1+  \int_{D^c} \frac{\delta_{\partial D}(z)^\alpha}{|x-z|^{d+\alpha}}dz\right).   
	 \end{align}
Observe  that $D^c \subset B(x,\delta_{\partial D}(x) + 2A_0) \setminus B(x,\delta_{\partial D}(x))$.  Using \eqref{e:dist-D-Dc}, we get that
\begin{align}\label{e:interaction-Neumann-unbounded-case2-2}  \int_{D^c} \frac{\delta_{\partial D}(z)^\alpha}{|x-z|^{d+\alpha}}dz &\le \int_{B(x,\delta_{\partial D}(x)  + 2A_0) \setminus B(x,\delta_{\partial D}(x))} \frac{dz}{|x-z|^{d}} \nn\\  &= c_{9}  	\log\left(1+\frac{ 2A_0}{\delta_{\partial D}(x)}\right)  \le 2c_9\log\left(e+\frac{A_0}{\delta_{\partial D}(x)}\right).	 \end{align}
Recall from \eqref{e:exterior-cone} that $\sC_{Q_x}(r_0, \theta) \subset D^c$.  Note that  for all $z\in \sC_{Q_x} (r_0/2,2\theta)$ 
with coordinates  $z=(\wt z,z_d) $ in CS$_{Q_x}$, since $ 2\theta |\wt z|<-z_d<r_0/2$, we have 
	\begin{gather*}
		\delta_{\partial D}(z) \ge \delta_{\partial \sC_{Q_x}(r_0,\theta)}(z) = \left( \frac{-z_d-\theta |\wt z|}{\sqrt{\theta^2+1}}\right) \wedge (r_0+z_d) \ge \frac{-z_d}{2\sqrt{\theta^2+1}}
	\end{gather*}
	and $	|x-z| \le |x| + |z| \le \delta_{\partial D}(x) + (-z_d) \sqrt{1+ (2\theta)^{-2}}$.	Thus, if $\delta_{\partial D}(x)<r_0$, then we get
	\begin{align*}
		 &\int_{D^c} \frac{\delta_{\partial D}(z)^\alpha}{|x-z|^{d+\alpha}}dz \ge  \int_{\sC_{Q_x}(r_0/2, 2\theta) } \frac{\delta_{\partial D}(z)^\alpha}{|x-z|^{d+\alpha}}dz 
		\ge c_{10}  \int_0^{r_0/2}  \frac{s^\alpha}{(\delta_{\partial D}(x) + s)^{d+\alpha}}\, \int^{s/(2\theta)}_{s/(4\theta)} \ell^{d-2} d\ell \,ds \nn\\
		  &\ge c_{10}  \int^{r_0/2}_{\delta_{\partial D}(x)/4}  \frac{s^\alpha}{(5 s)^{d+\alpha}}\, 
		  \int_{s/(4\theta)}^{s/(2\theta)} \ell^{d-2} d\ell \,ds  = c_{11} \int^{r_0/2}_{\delta_{\partial D}(x)/4} \frac{ds}{s}  = c_{11} \log \bigg(\frac{2r_0}{\delta_{\partial D}(x)}\bigg).
	\end{align*}
	Consequently, there exists $c_{12}=c_{12}(r_0,A_0)>0$ such that 
	\begin{align}\label{e:interaction-Neumann-unbounded-case2-3}
		1+ \int_{D^c} \frac{\delta_{\partial D}(z)^\alpha}{|x-z|^{d+\alpha}}dz \ge 1 + c_{11} \1_{\{\delta_{\partial D}(x)<r_0\}}\log \bigg(\frac{2r_0}{\delta_{\partial D}(x)}\bigg) \ge  c_{12}\log\left(e+\frac{A_0}{\delta_{\partial D}(x)}\right).
	\end{align}
Using $\delta_{\partial D}(x) <A_0$ and $ |y-x| \asymp \delta_{\partial D}(y) \ge 2A_0$, we see that
	 \begin{align}\label{e:interaction-Neumann-unbounded-case2-4}
	 	\log \bigg(e+ \frac{(|x-y|\wedge A_0)^2}{( \delta_{\partial D}(x)\wedge A_0)(\delta_{\partial D}(y) \wedge A_0) }\bigg) \asymp 	\log \bigg(e+ \frac{A_0}{\delta_{\partial D}(x)}\bigg) \asymp \log \bigg(e+ \frac{2A_0}{\delta_{\partial D}(x)}\bigg) .
	 \end{align}
	 Combining \eqref{e:interaction-Neumann-unbounded-case2-1} with \eqref{e:interaction-Neumann-unbounded-case2-2},  \eqref{e:interaction-Neumann-unbounded-case2-3}  and \eqref{e:interaction-Neumann-unbounded-case2-4}, we arrive at the result.
	
	 \smallskip

	 \noindent \textit{Case 3:} $\delta_{\partial D}(x)\ge A_0$. In this case, we have
	by \eqref{e:dist-D-Dc} and Lemma \ref{l:dist-D}(i), 
	 $\delta_{\partial D}(x) < |x-z| \le \delta_{\partial D}(x) + A_0 \le 2\delta_{\partial D}(x)$ and   $\delta_{\partial D}(y) < |y-z| \le 2\delta_{\partial D}(y)$ for all $z\in D^c$. Further, by Lemma \ref{l:dist-D}(ii), $|x-y| \le 3\delta_{\partial D}(y)$. Consequently, we get
	 \begin{align*}
	 	 \int_{D^c} \frac{\delta_{\partial D}(z)^\alpha}{|x-z|^{d+\alpha} |z-y|^{d+\alpha}}dz \le  \frac{4^{d+\alpha} \int_{D^c}\delta_{\partial D}(z)^\alpha dz }{\delta_{\partial D}(x)^{d+\alpha}\delta_{\partial D}(y)^{d+\alpha}} \le \frac{c_{13}}{A_0^{d+\alpha}\delta_{\partial D}(y)^{d+\alpha}} \le \frac{c_{14}}{|x-y|^{d+\alpha}}.
	 \end{align*} 
	 Thus, by  \eqref{e:interaction-Neumann-unbounded}, $K_D^\alpha(x,y) \asymp |x-y|^{-d-\alpha}$. Since 
	 \begin{align*}
	 1\le \log \bigg( e+ \frac{(|x-y|\wedge A_0)^2}{(\delta_{\partial D}(x) \wedge A_0)(\delta_{\partial D}(y) \wedge A_0)  }\bigg)  \le \log(e+ 1)
	 \end{align*}
	 in this case, 
	 the desired assertion is valid.

	 The proof is complete.
	 \qed

	  \subsection{Proof of Proposition \ref{p:kernelestimate-trace}} \label{ss:7.2}
	  
	 Throughout this subsection,	 we assume that $D$ satisfies {\bf (C)}.

	 \begin{lem}\label{l:reformulation-2}
	 	For all $x,y \in D$, 
	 	\begin{equation}\label{e:reformulation-2}
	 		\bigg( 1\wedge  \frac{|x-y|^2}{\delta_{\partial D}(x)\delta_{\partial D}(y)}\bigg)^{1/2} \asymp   1 \wedge \frac{|x-y|}{\delta_{\partial D}^\vee(x,y)} \asymp   \frac{|x-y|}{|x-y|+\delta_{\partial D}(x) + \delta_{\partial D}(y)}.
	 	\end{equation}
	 \end{lem} 
	 \pf 
	 By symmetry, it suffices to prove  \eqref{e:reformulation-2}  for  $x,y \in D$ with $\delta_{\partial D}(x)\le \delta_{\partial D}(y)$. 
	 If $\delta_{\partial D}(y) \ge 2|x-y|$, then $\delta_{\partial D}(x)\ge \delta_{\partial D}(y) - |x-y| \ge \delta_{\partial D}(y)/2$ so that $\delta_{\partial D}(x) \asymp \delta_{\partial D}(y)$ and $\delta_{\partial D}(x) \ge |x-y|$. Hence,  all  
	 the terms
	 in \eqref{e:reformulation-2} are comparable with  $|x-y|/\delta_{\partial D}(x)$. If $\delta_{\partial D}(y) < 2|x-y|$, then    $|x-y|^2 \ge \delta_{\partial D}(x)\delta_{\partial D}(y)/4$ and $|x-y| \ge (\delta_{\partial D}(x)+\delta_{\partial D}(y))/4 \ge \delta_{\partial D}(x)/2$. Thus, 
	 all  
	 the terms
	 in \eqref{e:reformulation-2} are comparable with $1$. The proof is complete.\qed

	 \begin{lem}\label{l:path-censor-pre}
	(i) 	There are comparison constants such that for all $y \in D$ and $z\in D^c$,
	 	\begin{align}\label{e:path-censor-pre-1}  \bigg( 1 \wedge \frac{\delta_{\partial D}(y)}{|y-z| \wedge R_0}\bigg)^{\alpha/2} \bigg( 1 \wedge \frac{\delta_{\partial D}(z)}{|y-z|\wedge R_0}\bigg)^{\alpha/2}\asymp \frac{\delta_{\partial D}(y)^{\alpha/2}\delta_{\partial D}(z)^{\alpha/2}}{|y-z|^{\alpha} } \bigg(1+ \frac{\delta_{\partial D}(z)}{R_0} \bigg)^{\alpha/2}.
	 	\end{align}
	 	
	 	\noindent (ii) For all $y\in D$ and $z\in D^c$, it holds that 
	 	\begin{align}\label{e:path-censor-pre-2}
	 	1+ \frac{|y-z|}{R_0} \le 
		3 \bigg(1+ \frac{\delta_{\partial D}(z)}{R_0} \bigg).
	 	\end{align}
	 \end{lem}
	 \pf (i) If $R_0=\infty$, then by \eqref{e:dist-D-Dc}, both sides of \eqref{e:path-censor-pre-1} are equal to $\delta_{\partial D}(y)^{\alpha/2} \delta_{\partial D}(z)^{\alpha/2}|y-z|^{-\alpha}$.
	 
	 Assume $R_0<\infty$. If $\delta_{\partial D}(z) \le R_0$, then by Lemma \ref{l:dist-D}(i), we get 
	 $|y-z| \le \delta_{\partial D}(z)   +R_0\le 2R_0$. 
	 	 Thus,  by \eqref{e:dist-D-Dc}, the left-hand side of \eqref{e:path-censor-pre-1}
	  is comparable to
	 \begin{align*}
	  \frac{\delta_{\partial D}(y)^{\alpha/2} \delta_{\partial D}(z)^{\alpha/2}}{|y-z|^\alpha} \asymp \frac{\delta_{\partial D}(y)^{\alpha/2}\delta_{\partial D}(z)^{\alpha/2}}{|y-z|^{\alpha} } \bigg(1+ \frac{\delta_{\partial D}(z)}{R_0} \bigg)^{\alpha/2}.
	 \end{align*}
 Suppose $\delta_{\partial D}(z) > R_0$.  By \eqref{e:dist-D-Dc} and Lemma \ref{l:dist-D}(i),  we see that 
 $|y-z|   \asymp \delta_{\partial D}(z)>R_0.$ 
	 Thus,  the left-hand side of \eqref{e:path-censor-pre-1} equals to $(\delta_{\partial D}(y)/R_0)^{\alpha/2}$ and the right-hand side of \eqref{e:path-censor-pre-1} is comparable to
$({\delta_{\partial D}(y)}/{R_0})^{\alpha/2}.$
	 The proof for (i) is complete.

	 \noindent (ii) If $R_0=\infty$, then \eqref{e:path-censor-pre-2} holds trivially. Assume $R_0<\infty$.  If $\delta_{\partial D}(z) < R_0$, then since 
	 $|y-z| < 
	 2R_0$ by Lemma \ref{l:dist-D}(i), the left-hand side of \eqref{e:path-censor-pre-2} is bounded above by $3$, and hence  \eqref{e:path-censor-pre-2} holds. If $\delta_{\partial D}(z) \ge R_0$, then since $|y-z| \le \delta_{\partial D}(z) +
	 R_0\le 2\delta_{\partial D}(z)$ by Lemma \ref{l:dist-D}(i), we get $1+ |y-z|/R_0 \le 
	 3\delta_{\partial D}(z)/R_0$. This completes the proof.
	 \qed

	 \begin{lem}\label{l:path-censor}
	 	 Suppose that $D$ is either bounded or half-space-like. Then there are comparison constants such that for all $y \in D$ and $z\in D^c$,
		\begin{align}\label{e:lemma7.8}
	 		&	\int_{D^c} \bigg( 1 \wedge \frac{\delta_{\partial D}(z)}{|z-w| \wedge R_0}\bigg)^{\alpha/2} \bigg( 1 \wedge \frac{\delta_{\partial D}(w)}{|z-w| \wedge R_0}\bigg)^{\alpha/2} \frac{dw}{|z-w|^{d-\alpha}|w-y|^{d+\alpha}}
			\nn\\
	 		&\asymp \frac{\delta_{\partial D}(z)^{\alpha/2}}{|y-z|^{d}  \, \delta_{\partial D}(y)^{\alpha/2}} \bigg(1+ \frac{\delta_{\partial D}(z)}{R_0} \bigg)^{\alpha/2}.
\end{align}
	 \end{lem}
	 \pf  
	 (Lower bound) We will use the notation $y^*$ from Lemma \ref{l:reflected-point}.
	 By 
	 \eqref{e:reflect-point}--\eqref{e:reflect-point2},
	 we have
	 \begin{align*}
	 	&\int_{D^c} \bigg( 1 \wedge \frac{\delta_{\partial D}(z)}{|z-w| \wedge R_0}\bigg)^{\alpha/2} \bigg( 1 \wedge \frac{\delta_{\partial D}(w)}{|z-w| \wedge R_0}\bigg)^{\alpha/2} \frac{dw}{|z-w|^{d-\alpha}|w-y|^{d+\alpha}}\\
	 	&\ge  c_1\bigg( 1 \wedge \frac{\delta_{\partial D}(z)}{|y-z| \wedge R_0}\bigg)^{\alpha/2} \bigg( 1 \wedge \frac{\delta_{\partial D}(y)}{|y-z| \wedge R_0}\bigg)^{\alpha/2}\frac{1}{|y-z|^{d-\alpha}\delta_{\partial D}(y)^{d+\alpha}} \int_{B(y^*, \delta_{\partial D}(y^*)/2)} dw\\
	 	&\ge c_2 \bigg( 1 \wedge \frac{\delta_{\partial D}(z)}{|y-z|\wedge R_0}\bigg)^{\alpha/2} \bigg( 1 \wedge \frac{\delta_{\partial D}(y)}{|y-z| \wedge R_0}\bigg)^{\alpha/2} \frac{1}{|y-z|^{d-\alpha}\delta_{\partial D}(y)^{\alpha}} .
	 \end{align*}
	 Combining this with Lemma \ref{l:path-censor-pre}(i), we get the  lower bound.
	 
	 \smallskip
	 
	 (Upper bound) 	 
             We rewrite the integral on the right-hand side of \eqref{e:lemma7.8}  as $I_1+I_2+I_3$, with
	 \begin{align*}
	 	I_i:=\int_{A_i} \bigg( 1 \wedge \frac{\delta_{\partial D}(z)}{|z-w| \wedge R_0}\bigg)^{\alpha/2} \bigg( 1 \wedge \frac{\delta_{\partial D}(w)}{|z-w| \wedge R_0}\bigg)^{\alpha/2} \frac{dw}{|z-w|^{d-\alpha}|w-y|^{d+\alpha}}, \quad i=1,  \dots,3,
	 \end{align*}
 where $A_1:=B(z,\delta_{\partial D}(z)/2)$, $A_2:=\{w \in D^c \setminus B(z,\delta_{\partial D}(z)/2): |w-y| \ge 2|z-w|\}$  and 	$A_3:=D^c \setminus (A_1 \cup A_2).$

	 For all $w\in A_1$, by \eqref{e:dist-D-Dc}, we have $|w-y| \ge |y-z|-|z-w| \ge |y-z|/2$. Using this, we get
	 \begin{align*}
	 	I_1&\le \int_{A_1} \frac{dw}{|z-w|^{d-\alpha}|w-y|^{d+\alpha}}	\le \frac{2^{d+\alpha}}{|y-z|^{d+\alpha}} \int_{B(z, \delta_{\partial D}(z)/2)} \frac{dw}{|z-w|^{d-\alpha}} \\
	 	&= \frac{c_3 \delta_{\partial D}(z)^{\alpha}}{|y-z|^{d+\alpha}} \le 
	 	\frac{c_3 \delta_{\partial D}(z)^{\alpha/2}}{|y-z|^{d}\delta_{\partial D}(y)^{\alpha/2}},
	 \end{align*}
	 where we used \eqref{e:dist-D-Dc} in the last inequality. 
	 Next,	 for all $w \in A_2$, we have $|y-z|\le |w-y| + |z-w| \le (3/2)|w-y|$ and $2|z-w|\le |w-y|\le |y-z|+|z-w|$ so that $|z-w|\le |y-z|$.  
	 Using the fact $a/(a \wedge  b) \le 1+(a/b)$ in the first inequality below, the above facts   in the second inequality,  Lemma \ref{l:path-censor-pre}(ii) and \eqref{e:dist-D-Dc} in the third, we obtain
	 \begin{align*}
	I_2	 	&\le \delta_{\partial D}(z)^{\alpha/2}\!	\int_{A_2} \!  \frac{(1+|z-w|/R_0)^{\alpha/2} dw}{|w-y|^{d+\alpha}|z-w|^{d-\alpha/2} }\le \frac{(3/2)^{d+\alpha}\delta_{\partial D}(z)^{\alpha/2}}{ |y-z|^{d+\alpha}}\!\int_{B(z, |y-z|)}\!\!
		\frac{(1+|z-w|/R_0)^{\alpha/2} dw}{|z-w|^{d-\alpha/2} }\nn\\
	 	&\le \frac{c_4\delta_{\partial D}(z)^{\alpha/2}}{|y-z|^{d+\alpha/2} }   \bigg(1+ \frac{|y-z|}{R_0} \bigg)^{\alpha/2}\le \frac{3^{\alpha/2} c_4\delta_{\partial D}(z)^{\alpha/2}}{|y-z|^{d}  \, \delta_{\partial D}(y)^{\alpha/2}} \bigg(1+ \frac{\delta_{\partial D}(z)}{R_0} \bigg)^{\alpha/2}.
	 \end{align*}
	 Thirdly, for all $w\in A_3$, we have	$|y-z|\le |w-y|+|z-w|< 3|z-w|$.   Using this  in the first inequality below and \eqref{e:dist-D-Dc} in the second, 
	 we obtain
	 \begin{align*}
	 	I_3&\le \frac{3^d}{|y-z|^{d-\alpha}}\bigg( 1 \wedge \frac{\delta_{\partial D}(z)}{|y-z| \wedge R_0}\bigg)^{\alpha/2}  \int_{A_3} \bigg( 1 \wedge \frac{\delta_{\partial D}(w)}{|y-z| \wedge R_0}\bigg)^{\alpha/2} \frac{dw}{|w-y|^{d+\alpha}}\\
	 	&\le \frac{3^d}{|y-z|^{d-\alpha}}\bigg( 1 \wedge \frac{\delta_{\partial D}(z)}{|y-z| \wedge R_0}\bigg)^{\alpha/2}  \int_{B(y, \delta_{\partial D}(y))^c} \bigg( 1 \wedge \frac{|w-y|}{|y-z| \wedge R_0}\bigg)^{\alpha/2} \frac{dw}{|w-y|^{d+\alpha}}\\
	 	&\le \frac{c_5}{|y-z|^{d-\alpha}}\bigg( 1 \wedge \frac{\delta_{\partial D}(z)}{|y-z| \wedge R_0}\bigg)^{\alpha/2} \\
	 	&\quad \times \bigg[ \frac{1}{(|y-z| \wedge R_0)^{\alpha/2}} \int_{\delta_{\partial D}(y)}^{\delta_{\partial D}(y) \vee ( |y-z| \wedge R_0)} \frac{d\ell}{\ell^{1+\alpha/2}} + \int_{\delta_{\partial D}(y) \vee ( |y-z| \wedge R_0)}^\infty \frac{d\ell}{\ell^{1+\alpha}}\bigg]\\
	 	&\le \frac{c_6}{|y-z|^{d-\alpha}\delta_{\partial D}(y)^{\alpha}}\bigg( 1 \wedge \frac{\delta_{\partial D}(z)}{|y-z| \wedge R_0}\bigg)^{\alpha/2} \bigg( 1 \wedge \frac{\delta_{\partial D}(y)}{|y-z| \wedge R_0}\bigg)^{\alpha/2} .
	 \end{align*}
	 Hence,  Lemma \ref{l:path-censor-pre}(i) yields the desired upper bound for $I_3$ and the proof is complete.
		 \qed

	 We recall the Green function estimates on $\overline{D}^c$.
	 
	 \begin{prop}\label{p:alpha-stable-Green}
	There are comparison constants such that  for all 
	 $z,w \in \overline{D}^c$,
	\begin{align*}
		g^{\alpha}_{\overline{D}^c}(z,w)
		\asymp \bigg( 1 \wedge \frac{\delta_{\partial D}(z)}{|z-w| \wedge R_0}\bigg)^{\alpha/2} \bigg( 1 \wedge \frac{\delta_{\partial D}(w)}{|z-w| \wedge R_0}\bigg)^{\alpha/2} \frac{1}{|z-w|^{d-\alpha}}. 
	\end{align*}
	 \end{prop}
	\pf If 
$D$ 
	is the complement of a bounded set, then the result follows from \cite[Theorem 5.7(ii)]{CS25+}, originally established in \cite[Corollary 1.3]{CS98}  for $C^{1,1}$ open sets.

 If $D$ is either bounded or half-space-like, then 
 the results follow from  the arguments in  \cite[Corollaries 1.4 and 1.5]{CT11} using the small-time heat kernel estimates (for $C^{1,\mathrm{Dini}}$ open sets) in \cite[Theorem 5.7(i)]{CS25+}.
 Indeed, although \cite{CT11} treats only $C^{1,1}$ open sets, the same argument extends to $C^{1,\mathrm{Dini}}$ open sets by using \cite[Theorem 5.7(i)]{CS25+}. We omit the details. \qed

	 \begin{prop}\label{p:alpha-stable-Poisson}
	 (i) If  $D$ is either bounded or half-space-like, then for all $y\in D$ and  
	 $z\in \overline{D}^c$,
	 	\begin{align*}
				P^\alpha_{\overline{D}^c}(z,y)
				 \asymp  \frac{\delta_{\partial D}(z)^{\alpha/2}}{|y-z|^{d}  \, \delta_{\partial D}(y)^{\alpha/2}} \bigg(1+ \frac{\delta_{\partial D}(z)}{R_0} \bigg)^{\alpha/2}.
	 	\end{align*}

	 	\noindent (ii) If $D$ is the complement of a bounded set,  then for all $y\in D$ and 
		$z\in \overline{D}^c$,
	 	\begin{align*}
			 P^\alpha_{\overline{D}^c}(z,y)
			\asymp \frac{\delta_{\partial D}(z)^{\alpha/2}}{|z-y|^d \, \delta_{\partial D}(y)^{\alpha/2} (\delta_{\partial D}(y)+A_0)^{\alpha/2}}.
	 	\end{align*} 
	 \end{prop}
	 \pf 
 (i) The result follows by combining Proposition \ref{p:alpha-stable-Green} with Lemma \ref{l:path-censor}.

	 \noindent (ii) Following the argument leading to \cite[Theorem 1.5]{CS98} (with $D$ replaced 
	 by $\overline{D}^c$),
	 the result follows from Proposition \ref{p:alpha-stable-Green}. 
	 \qed

	 \begin{lem}\label{l:trace-ND} 
	 	There are comparison constants such that  for all $x,y \in D$,
	 	\begin{align}\label{e:trace-ND} 
	 		N^\alpha_D(x,y) &\asymp \bigg( 1\wedge  \frac{|x-y|^2}{\delta_{\partial D}(x)\delta_{\partial D}(y)} \bigg)^{d/2} \bigg( 1 + \frac{\delta_{\partial D}^\wedge(x,y)}{A_0} \bigg)^{-d-\alpha/2}   \bigg( 1 + \frac{\delta_{\partial D}^\vee(x,y)}{A_0} \bigg)^{-\alpha/2}\nn\\
	 		&\quad\;\; \times \frac{1}{|x-y|^{d} \delta_{\partial D}(x)^{\alpha/2} \delta_{\partial D}(y)^{\alpha/2}}.
	 	\end{align}
	 \end{lem}
	\pf  Let $x,y \in D$. By symmetry, without loss of generality, we assume that $\delta_{\partial D}(x)\le \delta_{\partial D}(y)$.
	 
	  If $D$ is the complement of a bounded set,  then 
	  combining Propositions \ref{p:alpha-stable-Green} and \ref{p:alpha-stable-Poisson} with the proof of  \cite[Theorem 2.6]{BGPKR} (with $D$ replaced by $\overline{D}^c$), 
	  we obtain
	 \begin{align*}
	 	N^\alpha_D(x,y) &\asymp \begin{cases} 
	 		\displaystyle	\delta_{\partial D}(x)^{-d-\alpha} \delta_{\partial D}(y)^{-d-\alpha}	&\mbox{ if } \delta_{\partial D}(x) \ge A_0,\\[3pt]
	 		\displaystyle \delta_{\partial D}(x)^{-\alpha/2} \delta_{\partial D}(y)^{-d-\alpha}	&\mbox{ if } \delta_{\partial D}(x)<A_0\le \delta_{\partial D}(y),\\[3pt]
	 		\displaystyle  \delta_{\partial D}(x)^{-\alpha/2}\delta_{\partial D}(y)^{-\alpha/2}(|x-y|+\delta_{\partial D}(x)+\delta_{\partial D}(y))^{-d}	&\mbox{ if  }  \delta_{\partial D}(y)<A_0.
	 	\end{cases}
	 \end{align*}
	 On the other hand, by Lemma \ref{l:reformulation-2}, the right-hand side of \eqref{e:trace-ND} is comparable to
	 \begin{align*}
	 	& \begin{cases} 
	 		\displaystyle	\bigg( 1\wedge  \frac{|x-y|}{\delta_{\partial D}(y)} \bigg)^{d} |x-y|^{-d} \delta_{\partial D}(x)^{-d-\alpha}\delta_{\partial D}(y)^{-\alpha}&\mbox{ if } \delta_{\partial D}(x) \ge A_0,\\[3pt]
	 		\displaystyle	\bigg( 1\wedge  \frac{|x-y|}{\delta_{\partial D}(y)} \bigg)^{d} |x-y|^{-d} \delta_{\partial D}(x)^{-\alpha/2}\delta_{\partial D}(y)^{-\alpha}	&\mbox{ if } \delta_{\partial D}(x)<A_0\le \delta_{\partial D}(y),\\[3pt]
	 		\displaystyle	\delta_{\partial D}(x)^{-\alpha/2}\delta_{\partial D}(y)^{-\alpha/2}(|x-y|+\delta_{\partial D}(x)+\delta_{\partial D}(y))^{-d} &\mbox{ if } \delta_{\partial D}(y)<A_0.
	 	\end{cases} 
	 \end{align*}
	 By Lemma \ref{l:dist-D}(ii), if $\delta_{\partial D}(y) \ge A_0$, then $|x-y| \le 2\delta_{\partial D}(y) + A_0 \le 3\delta_{\partial D}(y)$ so that $1 \wedge \frac{|x-y|}{\delta_{\partial D}(y)} \asymp \frac{|x-y|}{\delta_{\partial D}(y)} $. Consequently, \eqref{e:trace-ND} holds in this case.

	We now assume that $D$ is either bounded or half-space-like. 
	 First note that, by Lemma \ref{l:reformulation-2},   
	 \begin{align}\label{e:path-censored-reformulation}
	 	  \bigg( 1\wedge  \frac{|x-y|^2}{\delta_{\partial D}(x)\delta_{\partial D}(y)} \bigg)^{d/2} \frac{1}{|x-y|^d}\asymp \bigg( 1\wedge  \frac{|x-y|}{\delta_{\partial D}(y)} \bigg)^{d} \frac{1}{|x-y|^d} \asymp \frac{1}{(|x-y| + \delta_{\partial D}(y))^{d}}  \ge (2R_0)^{-d}.
	 \end{align}
	 By Proposition \ref{p:alpha-stable-Poisson}(i),  we have
	 \begin{align*}
	 	N^\alpha_D(x,y)=N^\alpha_D(y,x)&\asymp\frac{1}{\delta_{\partial D}(x)^{\alpha/2}} \int_{D^c} \frac{\delta_{\partial D}(w)^{\alpha/2}}{|y-w|^{d+\alpha}|x-w|^{d}  }\bigg(1+ \frac{\delta_{\partial D}(w)}{R_0} \bigg)^{\alpha/2} dw=:I. 
	 \end{align*}
	 Let  $y^* \in D^c$ be the point satisfying \eqref{e:reflect-point}. By  \eqref{e:reflect-point1} and \eqref{e:reflect-point2}, it holds that
	 \begin{align}\label{e:path-censored-1}
	 I
	 	&\ge 	\frac{1}{\delta_{\partial D}(x)^{\alpha/2}} \int_{B(y^*, \delta_{\partial D}(y^*)/2)} \frac{\delta_{\partial D}(w)^{\alpha/2}}{|y-w|^{d+\alpha}(|x-y|+|y-w|)^{d}  } dw\nn\\
	 	&\ge 	  \frac{c_1}{\delta_{\partial D}(x)^{\alpha/2}\delta_{\partial D}(y)^{d+\alpha/2}(|x-y| + \delta_{\partial D}(y))^{d}  }\int_{B(y^*, \delta_{\partial D}(y^*)/2)}dw\nn\\
	 	&\ge  \frac{c_2}{\delta_{\partial D}(x)^{\alpha/2}\delta_{\partial D}(y)^{\alpha/2}(|x-y| + \delta_{\partial D}(y))^{d}  }.
	 \end{align}
We now obtain the desired lower bound from \eqref{e:path-censored-1}.
	 
	 For the upper bound,  we have 
\begin{align*}	I&\le \frac{2^{\alpha/2}}{\delta_{\partial D}(x)^{\alpha/2}} \bigg(  \int_{D^c} \frac{\delta_{\partial D}(w)^{\alpha/2}}{|y-w|^{d+\alpha}|x-w|^{d}  } \,dw + R_0^{-\alpha/2} \int_{w\in D^c: \delta_{\partial D}(w) \ge R_0} \frac{\delta_{\partial D}(w)^{\alpha}}{|y-w|^{d+\alpha}|x-w|^{d}  } \,dw \bigg) \\	 	&=:\frac{2^{\alpha/2}}{\delta_{\partial D}(x)^{\alpha/2}} \left(I_1+I_2\right).	 \end{align*}
	 Applying Lemma \ref{l:technical-1} with $a_1 = 0$, $a_2=\alpha$ and $a_3=\alpha/2$, we get 
	\begin{align*}	I_1&\le c_3\bigg( \frac{1}{(|x-y|+\delta_{\partial D}(y))^{d+\alpha/2}}  + \frac{\delta_{\partial D}(y)^{-\alpha/2}}{(|x-y|+\delta_{\partial D}(y))^{d}} \bigg) \nn\\	 	&\le\frac{2c_3}{\delta_{\partial D}(y)^{\alpha/2}(|x-y|+\delta_{\partial D}(y))^{d}} \le c_4\bigg( 1\wedge  \frac{|x-y|^2}{\delta_{\partial D}(x)\delta_{\partial D}(y)} \bigg)^{d/2}\frac{1}{|x-y|^{d}  \delta_{\partial D}(y)^{\alpha/2}}, \end{align*}
	 where we used \eqref{e:path-censored-reformulation} in the last inequality. For $I_2$, let us assume that $R_0<\infty$. For all $w\in D^c$ with $\delta_{\partial D}(w) \ge R_0$, 
	   by \eqref{e:dist-D-Dc} and Lemma \ref{l:dist-D}(i), we have $|x-w| \wedge |y-w| \ge \delta_{\partial D}(w) \ge R_0$ and 
	 $$
	 |y-w| \le |x-w| + |x-y| \le |x-w| + R_0 \le 2|x-w|.
	 $$It follows that
 \begin{equation*} 
	 	I_2	  \le \frac{c_5}{R_0^{\alpha/2}}\int_{B(y,R_0)^c} \frac{dw}{|y-w|^{2d} } =\frac{c_6}{R_0^{d+\alpha/2}} \le \frac{c_6 }{R_0^d\delta_{\partial D}(y)^{\alpha/2}}.
	 \end{equation*}
 Combining this with \eqref{e:path-censored-reformulation}, we get the desired upper bound for $I_2$. The proof is complete.
	 \qed

	 \noindent \textbf{Proof of Proposition \ref{p:kernelestimate-trace}.} Let $x,y \in D$. By symmetry, we assume $\delta_{\partial D}(x)\le \delta_{\partial D}(y)$. Recall that $\wc J^\alpha_D(x,y)=J^\alpha(x,y) + N_D^\alpha(x,y)$. 
	 Note that $\wc J^\alpha_D(x,y)\ge J^\alpha(x,y)\ge c_1|x-y|^{-d-\alpha}$. If $|x-y|\le \delta_{\partial D}(x)$, then by Lemma \ref{l:trace-ND},
	 \begin{align*}
	 	\wc J^\alpha_D(x,y) \le \frac{c_2}{|x-y|^{d+\alpha}}  +   \frac{c_3}{|x-y|^{d} \delta_{\partial D}(x)^{\alpha/2} \delta_{\partial D}(y)^{\alpha/2}} \le  \frac{c_2+c_3}{|x-y|^{d+\alpha}}.
	 \end{align*}
	 Thus, the result holds in this case.

	 Suppose $|x-y|>\delta_{\partial D}(x)$. To get the result, it suffices to show that
	 \begin{align}\label{e:trace-proof-claim} 
	 	N_D^\alpha(x,y)  +\frac{1}{|x-y|^{d+\alpha}} \asymp \left( \frac{(|x-y| \wedge A_0)^2}{(\delta_{\partial D}(x) \wedge A_0)(\delta_{\partial D}(y) \wedge A_0)} \right)^{\alpha/2} \frac{1}{|x-y|^{d+\alpha}}.
	 \end{align}
	 Note that $\delta_{\partial D}(y) \le 
	 \delta_{\partial D}(x) 
	 +|x-y| <2|x-y|$, and thus,
	 \begin{align}\label{e:trace-case-blowup}
	 	\frac{|x-y|^2}{\delta_{\partial D}(x)\delta_{\partial D}(y)} \wedge  \frac{(|x-y| \wedge A_0)^2}{(\delta_{\partial D}(x) \wedge A_0)(\delta_{\partial D}(y) \wedge A_0)} \ge \frac12.
	 \end{align} Hence,  by Lemma \ref{l:trace-ND}, it holds that 
	 \begin{align}\label{e:trace-proof-1} 
	 		N^\alpha_D(x,y) & \asymp \bigg( 1 + \frac{\delta_{\partial D}(x)}{A_0} \bigg)^{-d-\alpha/2}   \bigg( 1 + \frac{\delta_{\partial D}(y)}{A_0} \bigg)^{-\alpha/2} \bigg( \frac{|x-y|^2}{\delta_{\partial D}(x)\delta_{\partial D}(y)}\bigg)^{\alpha/2}\frac{1}{|x-y|^{d +\alpha}}.
	 \end{align}
	 If $A_0=\infty$, then, by  \eqref{e:trace-case-blowup},  \eqref{e:trace-proof-claim}  follows from \eqref{e:trace-proof-1}. Assume that $A_0<\infty$.		   We consider the following three cases separately.

	 \smallskip
	 
	 \noindent \textit{Case 1:}  $\delta_{\partial D}(y) < 2A_0$. By Lemma \ref{l:dist-D}(ii),  $|x-y| <5A_0$. Thus, we get from \eqref{e:trace-proof-1}  that 
	 \begin{align*}
	 	N^\alpha_D(x,y) & \asymp \bigg( \frac{|x-y|^2}{\delta_{\partial D}(x)\delta_{\partial D}(y)}\bigg)^{\alpha/2}\frac{1}{|x-y|^{d +\alpha}}\asymp \left( \frac{(|x-y| \wedge A_0)^2}{(\delta_{\partial D}(x) \wedge A_0)(\delta_{\partial D}(y) \wedge A_0)} \right)^{\alpha/2} \frac{1}{|x-y|^{d+\alpha}}.
	\end{align*} 
By \eqref{e:trace-case-blowup}, this implies  \eqref{e:trace-proof-claim}.

	 \noindent \textit{Case 2:} $\delta_{\partial D}(x)<A_0$ and $\delta_{\partial D}(y) \ge 2A_0$. In this case,  we have $|x-y| \ge\delta_{\partial D}(y) - \delta_{\partial D}(x) \ge \delta_{\partial D}(y)/2 \ge  A_0$. Further, by Lemma \ref{l:dist-D}(ii), it holds that $|x-y|\le (5/2)\delta_{\partial D}(y)$. Using \eqref{e:trace-proof-1} and $A_0\le |x-y| \asymp \delta_{\partial D}(y)$,  we obtain
	 \begin{align*}
	 	N^\alpha_D(x,y) & \asymp  \bigg(  \frac{\delta_{\partial D}(y)}{A_0} \bigg)^{-\alpha/2} \bigg( \frac{\delta_{\partial D}(y)}{\delta_{\partial D}(x)}\bigg)^{\alpha/2}\frac{1}{|x-y|^{d +\alpha}} \\
	 	&= \left( \frac{(|x-y| \wedge A_0)^2}{(\delta_{\partial D}(x) \wedge A_0)(\delta_{\partial D}(y) \wedge A_0)} \right)^{\alpha/2} \frac{1}{|x-y|^{d+\alpha}},
	 \end{align*}
	implying	 \eqref{e:trace-proof-claim} by \eqref{e:trace-case-blowup}.
	 \smallskip
	 
	 \noindent \textit{Case 3:} $\delta_{\partial D}(x)\ge A_0$. By Lemma \ref{l:dist-D}(ii), $|x-y|\le 3\delta_{\partial D}(y)$. Thus, by \eqref{e:trace-proof-1},
	\begin{align*}
		N^\alpha_D(x,y) \le c_4 \bigg( \frac{\delta_{\partial D}(y)}{A_0} \bigg)^{-\alpha/2} \bigg( \frac{\delta_{\partial D}(y)}{\delta_{\partial D}(x)}\bigg)^{\alpha/2}\frac{1}{|x-y|^{d +\alpha}} \le  \frac{c_4}{|x-y|^{d+\alpha}}.
	\end{align*}
	Since the right-hand side of \eqref{e:trace-proof-claim} equals to $|x-y|^{-d-\alpha}$ in this case, this implies \eqref{e:trace-proof-claim}.

	\smallskip
	
	The proof is complete. \qed

\end{document}